\pdfoutput=1
\RequirePackage{ifpdf}
\ifpdf 
\documentclass[pdftex]{sigma}
\else
\documentclass{sigma}
\fi

\DeclareMathOperator{\Pol}{Pol}

\usepackage{stmaryrd}
\numberwithin{equation}{section}
\newtheorem{The1}{Theorem}
\newtheorem{The}{Theorem}[section]

\newtheorem{Pro}{Proposition}[section]

\newtheorem{Lem}{Lemma}[section]
{\theoremstyle{definition}
\newtheorem{Def}{Definition}[section]
\newtheorem{Rem}{Remark}[section]
\newtheorem{Not}{Notations}[section]

}

\begin{document}

\allowdisplaybreaks

\renewcommand{\PaperNumber}{059}

\newcommand{\arXivNumber}{1205.2992}

\FirstPageHeading

\ShortArticleName{Conf\/igurations of an Articulated Arm and Singularities of Special Multi-Flags}

\ArticleName{Conf\/igurations of an Articulated Arm\\and Singularities of Special Multi-Flags}

\Author{Fernand PELLETIER~$^\dag$ and Mayada SLAYMAN~$^\ddag$}

\AuthorNameForHeading{F.~Pelletier and M.~Slayman}

\Address{$^\dag$~Universit\'e de Savoie, Laboratoire de Math\'ematiques (LAMA),\\
\hphantom{$^\dag$}~Campus Scientifique, 73376 Le Bourget-du-Lac Cedex, France}
\EmailD{\href{mailto:pelletier@univ-savoie.fr}{pelletier@univ-savoie.fr}}

\Address{$^\ddag$~Department of Mathematical Sciences, Faculty of Sciences II, Lebanese University, Lebanon}
\EmailD{\href{mailto:mslayman@ul.edu.lb}{mslayman@ul.edu.lb}}

\ArticleDates{Received January 29, 2013, in f\/inal form May 18, 2014; Published online June 05, 2014}

\Abstract{P.~Mormul has classif\/ied the singularities of special multi-f\/lags in terms of ``EKR class''
encoded by sequences $j_1,\dots, j_k$ of integers (see~[Singularity Theory Seminar, Warsaw University of Technology, Vol.~8, 2003, 87--100] and [\textit{Banach Center Publ.}, Vol.~65, Polish
  Acad. Sci., Warsaw, 2004, 157--178]).
However, A.L.~Castro and R.~Montgomery have proposed in~[\textit{Israel~J. Math.} \textbf{192} (2012),
  381--427] a~codif\/ication of singularities of
multi-f\/lags by {\bf RC} and {\bf RVT} codes.
The main results of this paper describe a~decomposition of each ``EKR'' set of depth $1$ in terms of {\bf RVT}
 codes as well as characterize such a~set in terms of conf\/igurations of an articulated arm.
Indeed, an analogue description for some ``EKR'' sets of depth~$2$ is provided.
All these results give rise to a~complete characterization of all ``EKR'' sets for $1\leq k\leq 4$.}


\Classification{53C17; 58K99; 70B15; 70Q05; 93A30}

\section{Introduction and results}\label{intro}

A {\it special multi-flag} of step $m\geq 1$ and length $k\geq 1$ is a~sequence
\begin{gather*}
\mathbb{D}: \ D=D_{k}\subset\hfill D_{k-1}\hfill\subset\dots\subset D_{j}\subset\dots\subset D_{1} \subset D_{0}= TM
\end{gather*}
of distributions of constant rank on a~manifold~$M$ of dimension $ (k+1)m+1$ which satisf\/ies the following conditions
(see~\cite{M1}):
\begin{enumerate}\itemsep=0pt
\item[(i)] $D_{j-1}=[D_j,D_j]$ is the distribution generated by all Lie brackets of sections of $D_{j}$;
\item[(ii)] $D_j$ is a~distribution of constant rank $(k-j+1)m+1$;
\item[(iii)] each Cauchy characteristic subdistribution $L(D_j)$ of $D_j$ is a~subdistribution of constant corank one in
each $D_{j+1}$ for $j=1,\dots, k-1$, and $L(D_k)=0$;
\item[(iv)] there exists a~completely integrable subdistribution $F\subset D_1$ of corank one in $D_1$.
\end{enumerate}
(See Section~\ref{multiflag} for a~more precise def\/inition.)

The notion of special multi-f\/lags is described in~\cite{M2, PR1}.
Furthermore, for $m\geq 2$, the existence of a~completely integrable subdistribution~$F$ of $D_1$ implies property
(iii).
This result was f\/irst proved in~\cite{LR} for regular points, and in~\cite{Ad2, SY} for the general case.
When such a~distribution~$F$ exists, it is then unique (see Remark~\ref{cov}).
For $m=1$, a~special multi-f\/lag is a~Goursat f\/lag, and in this case conditions (iii) and (iv) are automatically
satisf\/ied but such a~distribution~$F$ is not unique.
One fundamental result on Goursat f\/lags is the existence of locally universal Goursat distributions proved by
R.~Montgomery and M.~Zhitomirskii in~\cite{MZ1}.
More precisely, they def\/ine the monster Goursat manifold which is constructed by applying~$k$ successive Cartan
prolongations.
On the other hand, the kinematic system of a~car with $k-1$ trailers can be described by an appropriate Goursat
distribu\-tion~$\Delta_k$ on~$\mathbb{R}^2
\times
({{\mathbb{S}}}^1)^k$.
Moreover, this Goursat distribution $ \Delta_k$ is dif\/feomorphic to the Cartan prolongation of the distribution~$\Delta_{k-1}$ on~$\mathbb{R}^2\times({\mathbb{S}}^1)^{k-1}$ (see Appendix~D of~\cite{MZ1} or Theorem~3.3 of~\cite{PLR}).

A special multi-f\/lag can be considered as a~generalization of the notion of Goursat f\/lags and the fundamental result
of~\cite{Ad2} and~\cite{SY} is again obtained by Cartan prolongation (see also~\cite{M2}).
Consequently, in this situation, we can build a~monster tower by successive Cartan prolongations of~$T\mathbb{R}^{m+1}$
(see~\cite{Ad2, CH, CMH, SY}):
\begin{gather}\label{proj}
\dots\rightarrow P^k(m)\rightarrow P^{k-1}(m)\rightarrow\dots\rightarrow P^j(m)\rightarrow\dots\rightarrow
P^{1}(m)\rightarrow P^0(m):=\mathbb{R}^{m+1},
\end{gather}
where each manifold $P^j(m)$ is endowed with a~{\it typical} distribution $\Delta_j$, the Cartan prolongation
of $\Delta_{j-1}$ for $1\leq l\leq k$.
In a~similar way we can def\/ine a~natural notion of spherical prolongation, which also gives rise to a~tower of sphere
bundles (see Section~\ref{cartsph}):
\begin{gather}\label{sph}
\dots\rightarrow\hat{P}^k(m)\rightarrow \hat{P}^{k-1}(m)\rightarrow\dots\rightarrow \hat{P}^j(m)\rightarrow
\dots\rightarrow\dots\rightarrow \hat{P}^{1}(m)\rightarrow \hat{P}^0(m):=\mathbb{R}^{m+1}.\!\!
\end{gather}
Again, each manifold $ \hat{P}^j(m)$ is endowed with a~{\it typical} distribution $\hat{\Delta}_j$ which represents the
spherical prolongation of $\hat{\Delta}_{j-1}$ for $j\geq 1$. Notice that we have a~canonical $2$-fold covering:
\begin{gather*}
\hat{P}^j(m)\rightarrow {P}^j(m)
\end{gather*}
for any $j\geq 1$ and $m\geq 2$.

An articulated arm def\/ined in~\cite{SP1} or a~system of rigid bars def\/ined in~\cite{LR} is a~kinematic system which can
be described by a~special multi-f\/lag.
More precisely, the conf\/iguration space ${\mathcal C}^k(m)$ of such a~kinematic system is dif\/feomorphic
to $\mathbb{R}^{m+1}\times({\mathbb{S}}^m)^k$, and this system is characterized by a~distribution ${\mathcal D}_k$ which generates a~special
multi-f\/lag of length~$k$ (see Section~\ref{arm}).
Thus we obtain a~natural tower of sphere bundles
\begin{gather}\label{congtower}
{\mathcal C}^k(m)\rightarrow {\mathcal C}^{k-1}(m)\rightarrow\dots\rightarrow{\mathcal C}^j(m)\rightarrow {\mathcal
C}^{j-1}(m)\rightarrow \dots\rightarrow {\mathcal C}^{1}(m)\rightarrow {\mathcal C}^0(m):=\mathbb{R}^{m+1},
\end{gather}
where each map ${\mathcal C}^j(m)\rightarrow {\mathcal C}^{j-1}(m)$ is a~sphere bundle, and each manifold ${\mathcal
C}^j(m)$ is endowed with a~{\it typical} distribution ${\mathcal D}_j$ associated with the corresponding articulated arm
of length~$j$ on~$\mathbb{R}^{m+1}$ for $1\leq l\leq k$. Note that by convention, ${\mathcal C}^0(m)$ represents the space
$\mathbb{R}^{m+1}$ endowed with the distribution ${\mathcal D}_0=T\mathbb{R}^{m+1}$ (see Section~\ref{armsph}).

In this context, we have the following result\footnote{The reader can f\/ind this result in~\cite{P} with a~summarized proof.
As all the arguments used to show this theorem are also essential for proving Theorem~\ref{2} and~\ref{3} below,
thus in this paper we give a~complete proof of this result.}:

\begin{The1}[see Theorems~\ref{artmod} and~\ref{distribhyper}(2)]\label{1}\quad
\begin{enumerate}\itemsep=0pt
\item[$1.$] Let $\hat{\Delta}_j$ be the canonical distribution on $\hat{P}^j(m)$ obtained after~$j$ successive spherical
prolongations.
Then for each $m\geq 2$ and $1\leq j\leq k$, there exists a~diffeomorphism $F^j$ from $\hat{P}^j(m)$ to
$\mathcal{C}^j(m)$ such that:
\begin{itemize}\itemsep=0pt
\item[$(i)$] $\rho^j\circ F^j =F^{j-1} \circ \hat{\pi}^j$, where $\hat{\pi}^j: \hat{P}^j(m)\rightarrow \hat{P}^{j-1}(m)$
and
$\rho^j: \mathcal{C}^j(m)\rightarrow\mathcal{C}^{j-1}(m)$ are the canonical projections,
\item[$(ii)$] $F^j_*(\hat{\Delta}_j)=\mathcal{D}_j$.
\end{itemize}
\item[$2.$] The commutative diagrams
\begin{gather*}
\begin{matrix}
\hat{P}^k(m)\rightarrow \hat{P}^{k-1}(m)\rightarrow\cdots\rightarrow
\hat{P}^{1}(m)\rightarrow\hat{P}^0(m):=\mathbb{R}^{m+1}
\\
\hspace*{-11mm}\downarrow
\hspace{19mm}\downarrow\hspace{6.5mm}\cdots
\hspace{11mm}
\downarrow\hspace{15mm}\downarrow
\\
{\mathcal C}^k(m)\rightarrow \mathcal{C}^{k-1}(m)\rightarrow\dots\rightarrow \mathcal{C}^1(m)\rightarrow
\mathcal{C}^0(m):=\mathbb{R}^{m+1}
\\
\hspace*{-11mm}\downarrow
\hspace{19mm}\downarrow\hspace{6.5mm}\cdots
\hspace{11mm}
\downarrow\hspace{15mm}\downarrow
\\
{P}^k(m)\rightarrow {P}^{k-1}(m)\rightarrow\dots\rightarrow {P}^{1}(m)\rightarrow{P}^0(m):=\mathbb{R}^{m+1}
\end{matrix}
\end{gather*}
have the following properties:
\begin{itemize}
\item[$(i)$] in each horizontal tower, the horizontal map between the space number~$j$ and the space number $j-1$ $(1\leq j\leq k)$
is a~sphere fibration for the first two lines and a~projective space fibration for the last line;
\item[$(ii)$] in each column number~$j$ $(1\leq j\leq k)$, each vertical map between the first two lines is a~diffeomorphism
which sends the typical distribution over the source space to the typical distribution over the image space, and each
vertical map between the last two lines is a~two-fold covering which has the same property.
\end{itemize}
\end{enumerate}
\end{The1}

The singularities of special multi-f\/lags were f\/irst described by P.~Mormul in~\cite{M1, M2}.
This classif\/ication was based on a~generalization of Cartan prolongation, and on some
``operation'' denoted ${\bf j}$ which produces a~new ($m+1$)-distribution from the previous one.
In this way, P.~Mormul constructs a~coding system which labels the singularity classes of germs of special multi-f\/lags
which he calls ``extended Kumpera--Ruiz singularity classes of multi-f\/lags'' $-$ in short ``EKR classes'' $-$ (for more
details see Section~\ref{EKR}).
An EKR class is coded by a~sequence $j_1,\dots, j_k$ such that $j_{l+1}\leq 1+ \max \{j_1,\dots, j_l\}$.
The integer \mbox{$\max\{j_1,\dots,j_k\}-1$} is called the depth of the EKR class.

Recently, A.L.~Castro and R.~Montgomery  proposed, in~\cite{CMH}, a~codif\/ication of singularities of
multi-f\/lags founded on the tower of projective bundles~\eqref{proj} using {\bf RC} and {\bf RVT} codes.
This codif\/ication gives rise to a~new classif\/ication of the singularities of special multi-f\/lags in terms of {\bf RVT}
classes.
More precisely, in tower~\eqref{proj} one can def\/ine sub-towers by taking the tower of Cartan prolongation of any f\/iber
of $P^j(m)\rightarrow P^{j-1}(m)$.
Therefore, we obtain the ``baby monsters'' (see~\cite{CMH}).
It follows that in each vector space $\Delta_k(p)\subset T_pP^k(m)$ we have a~family of ``critical'' hyperplanes, coming
from these sub-towers.
We can note, that one of these hyperplanes is the vertical space $V_pP^k(m)$, i.e.\ the tangent space on a~f\/iber of
$P^k(m)\rightarrow P^{k-1}(m)$.
A point $p\in P^k(m)$ can be written as $p=(p_{k-1},z)$, where $p_{k-1}\in P^{k-1}(m)$ and~$z$ is a~line in $\Delta_{k-1}(p_{k-1})$.
Therefore,~$p$ is called {\it vertical} if~$z$ is tangent to the f\/iber at $p_{k-1}$, and {\it tangency} if~$z$ is not
vertical but belongs to one critical hyperplane.
Otherwise,~$p$ is called {\it regular}.
Thus, we can af\/fect to~$p$ a~word composed of letters $\{R,V,T\}$ such that the letter of rank~$l$ is either~$R$,~$V$ or~$T$,
depending on whether the projection of~$p$ onto $P^l(m)$ is regular, vertical, or tangency, respectively.

The main results of this paper are to give a~complete description of some EKR sets in terms of {\bf RVT} codes, as well as an interpretation of such EKR classes and {\bf RVT} classes, in terms of the conf\/igurations of an articulated arm.
To make it clear, we need to consider further def\/initions and notations.

Let $\omega$ be any word in {\bf RVT} code. We denote by $R^h$ or $T^h$ a~sub-word of $\omega$
which is a~sequence of~$h$ consecutive letters~$R$ or~$T$ if $h>0$, and no letter~$R$ or~$T$, if $h=0$ respectively.
Consider now the multi-f\/lag $\mathbb{D}$ on the conf\/iguration space ${\mathcal C}^k(m)$ associated to an
articulated arm $(M_{0},\dots,M_{k})$ on $\mathbb{R}^{m+1}$ of length $k\geq 1$ (see Section~\ref{arm}).
The EKR set $\Sigma_{j_{1}\dots j_{k}}$ is the set of conf\/igurations $q\in {\mathcal C}^k(m)$ such that the germ of
$\mathbb{D}$ at~$q$ belongs to the EKR class coded by the sequence $j_1,\dots,j_k$.
In the same way, the {\bf RVT} set ${\mathcal C}_{\omega}$ is the set of conf\/igurations $q\in {\mathcal C}^k(m)$ whose
{\bf RVT} code is $\omega$.
The depth of $\Sigma_{j_{1}\dots j_{k}}$ is the depth of the EKR class ${j_{1}\dots j_{k}}$.
Finally, for any EKR class of $1$-depth we will denote by $\{i_1,\dots, i_\nu\}$ the set $\{i\in\{1,\dots,k\}\,|\,  j_i=2\}$.
We then have:
\begin{The1}[see Theorem~\ref{stratEKR}]\label{2} Let $(M_0,\dots,M_k)$ be an articulated arm.
\begin{enumerate}\itemsep=0pt
\item[$1.$] Each EKR set $\Sigma_{j_1\dots j_k}$ of depth $1$ is an analytic manifold of codimension $\nu$ of ${\mathcal
C}^k(m)$.
\item[$2.$] A {\bf RVT} set ${\mathcal C}_\omega$ is contained in $\Sigma_{j_1\dots j_k}$ if and only if $\omega$ is of type
$R^{h_0}VT^{l_1}R^{h_1}\dots V T^{l_\nu}R^{h_{\nu}}$ and each letter~$V$ is exactly at rank $i_1,\dots, i_\nu$.
Such set is an analytic submanifold of $\Sigma_{j_1\dots j_k}$ of codimension $l_1+\dots+l_\nu$.
\item[$3.$] The EKR set $\Sigma_{j_1\dots j_k}$ is the disjoint union of the RVT sets ${\mathcal C}_\omega$, where $\omega$ is
any word of type ${R^{h_0}VR^{h}R^{h_0}VT^{l_1}R^{h_1}\dots V T^{l_\nu}R^{h_{\nu}}}$.
\end{enumerate}
\end{The1}

The following result gives an interpretation of EKR sets of depth $1$ in terms of orthogonality properties of an
articulated arm:
\begin{The1}[see Theorems~\ref{stratEKR}(2) and~\ref{RVTconf1}] \label{3} Let $(M_0,\dots,M_k)$ be an articulated arm.
\begin{enumerate}\itemsep=0pt
\item[$1.$] A configuration $q\in {\mathcal C}^k(m)$ of the articulated arm belongs to the EKR set $\Sigma_{j_1\dots j_k}$ of
depth~$1$ if and only if in this configuration the segments $[M_{i-2},M_{i-1}]$ and $[M_{i-1},M_{i}]$ are
orthogonal at $M_{i-1}$ for all $i=i_1,\dots,i_\nu$.
\item[$2.$] A configuration $q\in \Sigma_{j_1\dots j_k}$ belongs to the RVT set ${\mathcal C}_{R^{h_0}VT^{l_1}R^{h_1}\dots V
T^{l_\nu}R^{h_{\nu}}}\subset \Sigma_{j_1\dots j_k}$ if and only if, at~$q$, the only orthogonality constraint is that
each segment $[M_{{i_\lambda}+l-1},M_{{i_\lambda}+l}]$ is orthogonal to the direction on $\mathbb{R}^{m+1}$ generated by
$\overrightarrow{M_{i_\lambda-2}M_{i_\lambda-1}}$ for all $l=0,\dots,l_\lambda$ and $\lambda=1,\dots,\nu$.
\end{enumerate}
\end{The1}

This paper is self-contained and organized as follows.

We f\/irst recall, in Section~\ref{multiflag}, the context and the essential results about
special multi-f\/lags which will be used in this paper.
We present a~summary on Cartan prolongation and tower of projective bundles in Section~\ref{projtower}.
Spherical prolongations, tower of sphere bundles and their properties are developed in the last Section~\ref{cartsph}.

Section~\ref{tower} is devoted to the conf\/igurations of an articulated arm of length $k\geq 1$ in $\mathbb{R}^{m+1}$.
The space ${\mathcal C}^k(m)$ of such conf\/igurations is presented in Section~\ref{arm}.
The relation between the tower of sphere bundles~\eqref{sph} and the tower~\eqref{congtower} of conf\/iguration spaces
${\mathcal C}^k(m)$ is given in Section~\ref{armsph}. Finally we present the hyperspherical coordinates on ${\mathcal
C}^k(m)$ in Section~\ref{hypersphcoord}.
The reader can f\/ind the proof of Theorem~\ref{1}(1) in Section~\ref{armsph} and Theorem~\ref{1}(2) in
Section~\ref{hypersphcoord}.

In Section~\ref{code}, we present a~summary of the {\bf RC} and {\bf RVT} codes def\/ined in~\cite{CMH}, and we adapt
these codes to the context of tower of sphere bundles.
Section~\ref{vertarm} gives some interpretations of the property of {\it verticality} in terms of conf\/igurations of
an articulated arm.
In an analogous maner, some interpretations of the property of {\it tangency} are given in the last Section~\ref{tangearm}.

Section~\ref{relatiomEKRRVT} is devoted to the relation between EKR sets of depth $1$ and {\bf RVT} sets.
In Section~\ref{EKR}, we summarize the def\/inition and the results concerning EKR classes based on~\cite{M1, M2}.
Section~\ref{stratEKR1} gives a~global description of EKR sets in terms of {\bf RVT} sets.
Section~\ref{depth1} presents an interpretation of EKR sets (of
depth at most $1$) and {\bf RVT} sets in terms of the configurations of an articulated arm.
Finally, in Section~\ref{2k4} a~characterization for some EKR sets of depth $2$ in terms of articulated arms
is given.
In this paragraph, we also give, for $1\leq k\leq 4$, the decomposition of EKR sets of depth at most $2$ in {\bf RVT}
sets and the corresponding interpretation in terms of conf\/igurations of an articulated arm.
The reader can f\/ind the proof of Theorems~\ref{2} and~\ref{3} in Sections~\ref{stratEKR1} and~\ref{depth1}
respectively.

\section{Preliminaries}\label{prelim}

\subsection{Special multi-f\/lags}\label{multiflag}

A distribution~$D$ on a~manifold~$M$ is an assignment $D: x\in M\mapsto D_x \subset TM$ where $D_x$ is a~linear subspace
of the tangent space $T_xM$.
A local vector f\/ield~$X$ on~$M$ is tangent to~$D$ if~$X(x)$ belongs to~$D_x$ for all~$x$ in the open set on which~$X$ is
def\/ined.
The distribution~$D$ is smooth if there exists a~set ${\mathcal X}$ of local vector f\/ields such that~$D_x$ is generated
by the set $\{X(x),X\in {\mathcal X}\}$ for all~$x$ in some open set~$U$.
We then say that~$D$ is generated by ${\mathcal X}$ on~$U$.

In this paper, all distributions are smooth and we denote by $\Gamma(D)$ the set of all local vector f\/ields
which are tangent to~$D$.
A~distribution will be called a~distribution of constant rank if~$D$ def\/ines a~subbundle of $TM$.
According to~\cite{Ad2} and~\cite{SY}, any pair $(M,D)$ where~$D$ is a~distribution of constant rank on a~smooth
manifold~$M$ is called a~{\it differential system}.
Given two dif\/ferential systems $(M,D)$ and $(N,\Delta)$ and two points $x\in M$ and $y\in N$, we will say that $(M,D,
x)$ is {\it locally equivalent} to $(N,\Delta,y)$ if there exists a~dif\/feomorphism $\phi$ from an open neighborhood~$U$
of~$x$ in~$M$ to a~neighborhood~$V$ of~$y$ in~$N$ such that $y=\phi(x)$ and $\phi_*(D_{| U})=\Delta_{| V}$.

Given a~distribution $D'$ on~$M$ such that $D'_{x}\subset D_{x}$ for all $x\in M$, we denote by $[D',D]$ the
distribution generated by the sets $\Gamma(D)$ and $\{[X,Y] \,|\,   X\in \Gamma(D'),\,  Y\in \Gamma(D)\}$.
The {\it Lie square} of a~distribution~$D$ is the distribution $D^2:=[D,D]$.
The {\it Cauchy characteristic distribution} $L(D)$ of a~distribution~$D$ is the distribution generated by the set
of vector f\/ields $\{X\in \Gamma(D)\,|\,  [X,Y]\in D, \,
 \forall\, Y\in \Gamma(D)\}$. If $L(D)$ def\/ines a~distribution of constant rank, then it is an integrable distribution.

A {\it special multi-flag} of step $m\geq 2$ and length $k\geq 1$ is a~sequence of distributions
\begin{gather*}
\mathbb{D}:\ D=D_{k}\subset\hfill D_{k-1}\hfill\subset\dots\subset D_{j}\subset\dots\subset D_{1} \subset D_{0}= TM
\end{gather*}
all of constant rank on a~manifold~$M$ of dimension $ (k+1)m+1$ which fulf\/ills the following conditions
(see~\cite{M1}):
\begin{enumerate}\itemsep=0pt
\item[(i)] $D_{j-1}=(D_j)^2$,
\item[(ii)] $D_j$ is a~distribution of constant rank $(k-j+1)m+1$,
\item[(iii)] each Cauchy characteristic subdistribution $L(D_j)$ of $D_j$ is of constant corank one in each~$D_{j+1}$,
for $j=1,\dots,k-1$, and $L(D_k)=0$,
\item[(iv)] there exists a~completely integrable subdistribution $F\subset D_1$ of corank one in $D_1$.
\end{enumerate}
In the sequel, a~f\/lag $\mathbb{D} $ which satisf\/ies conditions (i) and (ii) without conditions (iii) and (iv) will just
be called a~{\it multi-flag} of step~$m$ and we say that $\mathbb{D} $ is generated by~$D$.

The necessary and suf\/f\/icient condition for a~multi-f\/lag to be a~special multi-f\/lag is given by the following result
(see~\cite[Proposition 1.3]{Ad2}  and~\cite[Theorem 6.2]{SY}):
\begin{The}[see~\protect{\cite{Ad2,SY}}]
For $k\geq 2$ and $m\geq 1$ consider a~multi-flag of step~$m$:
\begin{gather*}
\mathbb{D}: \ D=D_{k}\subset\hfill D_{k-1}\hfill\subset\dots\subset D_{j}\subset\dots\subset D_{1} \subset D_{0}= TM.
\end{gather*}
$\mathbb{D}$ is a~special multi-flag if and only if there exists a~completely integrable subbundle $F $ of $D_1$ of corank~$1$.
Moreover, if the subbundle~$F$ exists, then it is unique.
\end{The}
\begin{Rem}
\label{cov}
The existence of the subbundle~$F$ in the previous theorem is crucial and is uniquely determined by the distribution $D_1$.
In fact, for any subbundle ${\mathcal D}$ of $TM$, the subbundle~$F$ of $\mathcal D$ was f\/irstly def\/ined by Kumpera
and Rubin in~\cite{KuRu}.
Let us give some details about this fundamental fact (all the following af\/f\/irmations are proved in~\cite{Ad1, KuRu}).

Let ${\mathcal D}^\perp\subset T^*M$ be the annihilating Pfaf\/f\/ian system of ${\mathcal D}$, the {\it polar system}
$\Pol({\mathcal D}^\perp)$ of the Pfaf\/f\/ian system $(\mathcal D)^\perp$ is def\/ined by
\begin{gather*}
\Pol\big({\mathcal D}^\perp\big)(x) =\big\{\alpha\in T_x^*M/{\mathcal D}^\perp(x) \,|\,   \alpha\wedge\delta\omega=0, \, \forall\,
\omega\in {\mathcal D}^\perp\big\},
\end{gather*}
where $\delta:{\mathcal D}^\perp\rightarrow\Lambda^2(T^*M/{\mathcal D}^\perp)$ is the Martinet tensor characterized by
$\delta\omega=d\omega$ $({\rm mod}\, {\mathcal D}^\perp)$.
Then, the {\it covariant system} $\widehat{{\mathcal D}^\perp}$ associated with ${\mathcal D}^\perp$ is
$q^{-1}(\Pol({\mathcal D}^\perp))\subset T^*M$ where $q:T^*M\rightarrow T^*M/{\mathcal D}^\perp$ is the
canonical projection.
It can be proved that, the announced subdistribution $F\subset {\mathcal D}$ is the distribution whose annihilating
Pfaf\/f\/ian system is $\widehat{{\mathcal D}^\perp}$ (see~\cite{Ad1}).
When the corank of $\mathcal D$ is at most $2$ and if ${\mathcal D}^2=TM$, then any such distribution~$F$ has corank $1$
in $\mathcal D$ and~$F$ is completely integrable (see~\cite{Ad1, KuRu}).
\end{Rem}

According to the previous def\/inition of a~special multi-f\/lag, we obtain the following sandwich f\/lag:
\begin{gather*}
\begin{matrix}
D_{k}&\subset\hfill D_{k-1}\hfill&\subset\dots\subset &D_{j}&\subset\dots\subset&D_{1} &\subset & D_{0}= TM\hfill
\\
\hfill\cup\hfill&\hfill
\qquad
\cup\hfill& \dots &\cup &\dots &\cup\hfill & &
\\
\hfill L(D_{k-1})&\subset L(D_{k-2})&\subset\dots\subset &
L(D_{j-1})&\subset\dots\subset& F\hfill& &
\end{matrix}
\end{gather*}

All vertical inclusions in this diagram are of codimension one, while all horizontal inclusions are of codimension~$m$.
From these inclusions, we can extract the following ``squares subdiagrams'' called ``sandwiches'', indexed by number~$j$
which corresponds to the upper left vertices $D_{j}$:
\begin{gather*}
\begin{matrix}
D_{j}&\subset&D_{j-1}
\\
\cup&&\cup
\\
L(D_{j-1})&\subset&L(D_{j-2})
\end{matrix}
\end{gather*}
This subdiagram is called a~{\it sandwich of rank~$j$}.

Given a~sandwich of rank~$j$, and a~point $x\in M$, we can look for the relative positions of the~$m$ dimensional
subspace $L(D_{j-2})/L(D_{j-1})(x)$ and the $1$-dimensional subspace $D_{j}/L(D_{j-1})(x)$ in the $(m+1)$-dimensional
vector space $D_{j-1}/L(D_{j-1})(x)$.
One (and only one) of the following situations then occurs:
\begin{itemize}\itemsep=0pt
\item[(i)] $L(D_{j-2})/L(D_{j-1})(x)\oplus D_{j}/L(D_{j-1})(x)=D_{j-1}/L(D_{j-1})(x)$,
\item[(ii)] $ D_{j}/L(D_{j-1})(x)\subset L(D_{j-2})/L(D_{j-1})(x)$.
\end{itemize}
\begin{Def}
\label{Cart}
A point $x\in M$ is called a~{\it Cartan point} if $k=1$ or if, for $k\geq 2$, the previous situation (i) is true in
each sandwich of rank~$j$, for $j=2,\dots,k$.
Otherwise~$x$ is called a~{\it singular point}.
\end{Def}

\subsection{Cartan prolongation and tower of projective bundles}\label{projtower}

Consider a~distribution~$D$ of constant rank $m+1$ on a~manifold~$M$ of dimension~$n$.
Classically the Grassmannian bundle $G(D,1,M)$ over~$M$ is the set
\begin{gather*}
G(D,1,M):= \bigcup_{x\in M} P(D(x),1),
\end{gather*}
where $P(D(x),1)$ is the projective space of the vector space $D(x)$.
Thus we have a~bundle $\pi:G(D,1,M)\rightarrow M$ whose f\/iber $\pi^{-1}(x)$ is dif\/feomorphic to the projective space $\mathbb{R} P^m$.
The {\it rank one Cartan prolongation} of the distribution~$D$ is the distribution $D^{(1)}$ def\/ined as follows: given a~point $(x,\lambda)\in G(D,1,M)$, we set
\begin{gather*}
D^{(1)}_{(x,\lambda)}:= d\pi^{-1}(\lambda)\subset T_{(x,\lambda)} G(D,1,M),
\end{gather*}
where $\lambda$ is a~direction of $D(x)$.
Then $D^{(1)}$ is a~distribution on $G(D,1,M)$ of constant rank $m+1$.
Let~$M$ be a~manifold of dimension $m+1$.
According to~\cite{SY}, for any $m\geq 2$ and $k\geq 1$ starting with $D=TM$, we obtain inductively a~tower of bundles
\begin{gather}\label{towerM}
P^{k}(M)\rightarrow P^{k-1}(M)\rightarrow\dots\rightarrow P^{1}(M)\rightarrow P^0(M):=M,
\end{gather}
where the sequences $(P^j(M))_{j=0,\dots,k}$ and $(\Delta_{j})_{j=0,\dots,k}$ are def\/ined inductively by
\begin{gather*}
P^0(M)=M
\qquad
\text{and}
\qquad
\Delta_0=TM,
\\
P^{j}(M)=G(\Delta_{j-1},1,P^{j-1}(M))
\qquad
\text{and}
\qquad
\Delta_{j}=(\Delta_{j-1})^{(1)}
\quad
\text{for}
\quad
j=1,\dots,k.
\end{gather*}

Notice that $P^{j}(M)$ is a~manifold of dimension $(j+1)m+1$ for $j=0,\dots,k$.

In the particular case of $M=\mathbb{R}^{m+1}$, we denote by $P^j(m)$ the manifold $P^j(\mathbb{R}^{m+1})$ for
$j=0,\dots,k$, and we obtain the corresponding tower of bundles:
\begin{gather*}
P^k(m)\rightarrow P^{k-1}(m)\rightarrow\dots\rightarrow P^{1}(m)\rightarrow P^0(m):=\mathbb{R}^{m+1}.
\end{gather*}

Then we have the following result:
\begin{The}[see~\protect{\cite{SY}}]\label{modelP}\quad
\begin{enumerate}\itemsep=0pt
\item[$1.$] On $ P^k(m)$, the distribution $\Delta_k$ generates a~special multi-flag of step~$m$ and length~$k$.
\item[$2.$] Given $x\in M$ and a~special multi-flag $\mathbb{D}:D=D_{k}\subset D_{k-1}\subset\dots\subset
D_{j}\subset\dots\subset D_{1} \subset D_{0}= TM$ of step $m\geq 2$ and length $k\geq 1$, there exists $y\in P^k(m)$ for
which the differential systems $(P^k(m),\Delta_k,y)$ and $(M,D,x) $ are locally equivalent.
\end{enumerate}
\end{The}

\begin{Rem}
 Theorem~\ref{modelP}(2) can be found precisely in~\cite{SY} and is called the ``Drapeau theorem''.
However, according to the def\/inition of a~special multi-flag, we can easily deduce this result from the following
theorem of~\cite{M2}:
\begin{The}[see \protect{\cite{M2}}]
Suppose that~$D$ is a~$(m+1)$-dimensional distribution on a~$(s+m)$-dimensional manifold~$M$ such that the two following
conditions are satisfied:
\begin{enumerate}\itemsep=0pt
\item[$1)$] $D_1=[D,D]$ is a~$(2m+1)$-dimensional distribution on~$M$,
\item[$2)$] there exists a~$1$-codimensional involutive subdistribution $E \subset D$ that preserves $D_1$, i.e.\
$[E,D_{1}] \subset D_{1}$.
\end{enumerate}
Then~$D$ is locally equivalent to the Cartan prolongation $(D_{1}/E)^{(1)}$ of the reduction $(D_{1}/E)$ of~$D_{1}$
modulo~$E$.
\end{The}
\end{Rem}

\subsection{Spherical prolongation, Cartan prolongation and tower of sphere bundles}\label{cartsph}

\begin{figure}[t]\centering
\includegraphics{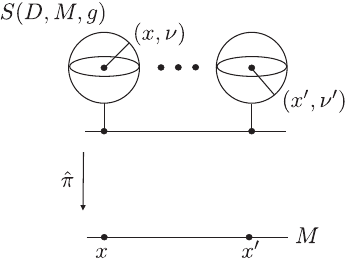}
\caption{Spherical united bundle.}\label{unitsphere}
\end{figure}

Consider a~distribution~$D$ of constant rank $m+1$ on a~manifold~$M$ of dimension~$n$.
Choose any Riemannian metric~$g$ on~$M$, and denote by $S(D,M,g)$ the unit sphere bundle of~$D$ associated with the
induced Riemannian metric on~$D$.
Then we obtain a~bundle $\hat{\pi}: S(D,M,g)\rightarrow M$ (see Fig.~\ref{unitsphere}).

Consider the antipodal action of $\mathbb{Z}_2$ on $S(D,M,g)$.
Clearly, the quotient of $S(D,M,g)$ by this action can be identif\/ied with $G(D,1,M)$ and the associated projection
$\tau:S(D,M,g)\rightarrow G(D,1,M)$ is both a~bundle morphism over~$M$ and a~two-fold covering.
In particular, $\tau$ is a~local dif\/feomorphism.
Consider now the distribution $D^{[1]}$ on $S(D,M,g)$ def\/ined by
\begin{gather*}
D^{[1]}_{(x,\nu)}:= \{v\in T_{(x,\nu)}S(D,M,g) \,|\,   d\hat{\pi}(v)=\lambda \nu~\text{for some}~\lambda\in\mathbb{R}\},
\end{gather*}
where $\nu$ is a~unit vector in $D(x)$.

\begin{figure}[t]\centering
\includegraphics{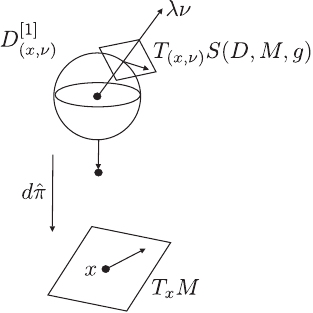}
\caption{Spherical prolongation on sphere bundle.}\label{prolong}
\end{figure}

The distribution $D^{[1]}$ is called the {\it rank one spherical prolongation} of $(M,D,g)$ (see Fig.~\ref{prolong}).

\begin{Rem}
The unit sphere bundle associated with~$D$ is def\/ined as soon as we f\/ix some Riemannian metric $g_D$ on~$D$.
In this case, the distribution $D^{[1]}$ is  well def\/ined and depends only on the Riemannian metric $g_D$ on~$D$.
The spherical prolongation $D^{[1]}$ depends only on the sub-Riemannian structure $(M,D,g_D)$.
However, for the sake of simplicity, we always consider Riemannian metrics on~$M$.
\end{Rem}
\begin{Lem}\label{spprolong}
With the previous notations we have
\begin{enumerate}\itemsep=0pt
\item[$(i)$] $\tau_*D^{[1]}=D^{(1)}$, and
\item[$(ii)$] there exists a~canonical Riemannian metric $\hat{g}$ on $S(D,M,g)$ which is uniquely defined from the
Riemannian metric~$g$ on~$M$.
\end{enumerate}
\end{Lem}
\begin{proof}
First of all we show (i) locally.
Choose a~chart domain~$U$ over which~$D$ is trivial.
Fix an orthonormal frame $\{e_0,\dots, e_m\}$ of~$D$ over~$U$.
Without loss of generality, we can assume that $D|_{U}\equiv \mathbb{R}^n
\times
\mathbb{R}^{m+1}$ therefore the bundle $S(D,M,g)|_{U}$ is isomorphic to $\mathbb{R}^n
\times
{\mathbb{S}}^m$, and $G(D,1,M)|_{U}$ is isomorphic to $\mathbb{R}^n
\times
\mathbb{R} P^m$.
Then, locally, $\tau: \mathbb{R}^n
\times
{\mathbb{S}}^m\rightarrow \mathbb{R}^n
\times
\mathbb{R} P^m$ is the map $(x,\nu)\rightarrow (x,[\nu])$, where $[\nu]$ is the line bundle generated by $\nu$.
According to the def\/inition of $D^{[1]}_{(x,\nu)}$ and $D^{(1)}_{\tau(x,[\nu])}$, we have
$\tau_*\big(D^{[1]}_{(x,\nu)}\big)=D^{(1)}_{\tau(x,[\nu])}$.
Since $\tau$ is a~local dif\/feomorphism then~(i) is proved  locally.
The map $\hat{\alpha}:S(D,M,g)\rightarrow S(D,M,g)$ given by $\hat{\alpha}(x,\nu)=(x,-\nu)$ is a~dif\/feomorphism which
commutes with $\tau$.
According to the def\/inition of $D^{[1]}$, we obtain
\begin{gather*}
\hat{\alpha}_*\big(D^{[1]}_{(x,\nu)}\big)=D^{[1]}_{(x,-\nu)}.
\end{gather*}
This ends the proof of~(i).

For~(ii), let $\bar{g}$ be the canonical Riemannian metric on $TM$ associated with~$g$.
Since $S(D,M,g)$ can be considered as a~submanifold of $TM$, the metric $\bar{g}$ induces a~Riemannian metric $\hat{g}$
on $S(D,M,g)$.
\end{proof}

Consider two Riemannian metrics $g_0$ and $g_1$ on~$M$.
We denote by $S_i(D,M)$ the sphere bundle of~$D$ associated with the metric $g_i$, and $D^{[1]}_i$ the spherical
prolongation of $(M,D,g_i)$ for $i=0,1$.
\begin{Lem}
\label{changeg}
There exists a~canonical isomorphism of sphere bundles $\psi: S_0(D,M)\rightarrow S_1(D,M)$ such that
$\psi_*\big(D^{[1]}_0\big)=D^{[1]}_1$
\end{Lem}
\begin{proof}
The set $D^{0}: = \bigcup\limits_{x\in M}[D_x\setminus\{0\}]$ is an open submanifold of $D\subset TM$ on which we
consider the map $\Psi:D^{0}\rightarrow D^{0} $ def\/ined by
\begin{gather*}
\Psi(x,u):=\left(x, \frac{u}{[g_1(u,u)]^{1/2}}\right).
\end{gather*}
If $\Pi: D\rightarrow M$ is the projection bundle, then for any $(x,u)\in D$ there exists a~neighborhood
$\hat{U}=\Pi^{-1}(U)\cap D^{0}$ around $(x,u)$ in $D^{0}$ such that $TD^{0}_{| \hat{U}}$ can be identif\/ied with $\hat{U}\times T_xM\times D_x$.
Then we have:
\begin{gather}\label{difPsi}
d\Psi(v,w)=\left(v, - \frac{g_1(u,w)}{2[g_1(u,u)]^{3/2}}\right).
\end{gather}
$\Psi$ is a~dif\/feomorphism from $D^{0}$ onto itself that commutes with $\Pi$ and sends $S_0(D,M)$
to $S_1(D,M)$.
It follows that the restriction $\psi$ of $\Psi$ to $S_0(D,M)$ is a~dif\/feomorphism onto $S_1(D,M)$.
Moreover, equation~\eqref{difPsi} shows that for any~$u$ in the f\/iber $D^{0}_x$ over~$x$, $d\Psi$ maps the linear span
$\mathbb{R} u$ onto itself.
Thus we have
\begin{gather*}
\psi_*\big(D^{[1]}_0\big)=D^{[1]}_1.\tag*{\qed}
\end{gather*}
\renewcommand{\qed}{}
\end{proof}

Consider a~dif\/ferential system $(M',D')$ and $\phi:M\rightarrow M'$ an injetive immersion such that $\phi_*(D_x)\subset
D'_{\phi(x)}$ for any $x\in M$.
Any Riemanian metric $g'$ on $M'$, induces, via $\phi$, a~Riemannian metric~$g$ on~$M$, and we can consider the
associated spherical prolongation.
This generates the following lemma:
\begin{Lem}\label{prolongmap}
With the above notations, the map $\hat{\phi}: S(D,M,g)\rightarrow S(D',M',g')$ defined by
\begin{gather*}
\hat{\phi}(x,\nu)=(\phi(x),d_x\phi(\nu))
\end{gather*}
is a~bundle morphism over $\phi$, which is an injective immersion, and $\hat{\phi}$ satisfies the following properties:
\begin{itemize}\itemsep=0pt
\item[$(i)$] $\hat{\phi}(S(D,M,g))=S(\phi_*(D),\phi(M),g')$, and
\item[$(ii)$] $\hat{\phi}_*\big(D^{[1]}\big)=(\phi_*(D))^{[1]}\subset (D')^{[1]}$.
\end{itemize}
Moreover, if $\phi$ is a~diffeomorphism such that $\phi_*(D)=D'$, then $\hat{\phi}$ is also a~diffeomorphism and we have
$\hat{\phi}_*\big(D^{[1]}\big)=(D')^{[1]}$.
On the other hand, the Riemannian metric $\hat{\phi}_*\hat{g'}$ is nothing else but the canonical metric~$\hat{g}$
naturally associated with~$g$ on~$M$.
\end{Lem}

\begin{proof}
First of all it is clear that $\hat{\phi}$ is smooth and is a~bundle morphism over $\phi$.
Moreover, $\hat{\phi}$ is injective since $\phi$ is an injective immersion.

Note that the tangent space $T_{(x,\nu)}S_x$ of the f\/iber $S_x$ over~$x$ of $S(D,M,g)$ can be identif\/ied with
\begin{gather*}
\{v \in D_x \,|\,  g(\nu,v)=0 \}.
\end{gather*}
Now, any $V\in T_{(x,\nu)}S(D,M,g)$ can be written as $V=(u,v)$ with $u\in T_xM$ and $v\in T_{(x,\nu)}S_x$.
Consequently we get
\begin{gather}\label{immers}
d_{(x,\nu)}\hat{\phi}(u,v)= (d_x\phi(u), d_x\phi(v)).
\end{gather}
It follows that $ \hat{\phi}$ is an immersion by equation~\eqref{immers}.

Indeed, since $\phi^*g'=g$, the dif\/ferential $d_x\phi$ is an isometry on its range, and then $d_x\phi(S_x)$
is the f\/iber over $\phi(x)$ of $S(\phi_*(D),\phi(M),g')$.
Thus~(i) is proved.

Let $\hat{\pi}:S(D,M,g)\rightarrow M$ and $\hat{\pi}':S(D',M',g')\rightarrow M'$ be the natural projections.
Then we have
\begin{gather*}
d\hat{\pi}'\circ d\hat{\phi}=d\phi\circ d\hat{\pi},
\end{gather*}
which yields
\begin{gather*}
\big\{\hat{\phi}_*(D^{[1]})\big\}_{\hat{\phi}(x,\nu)}=\big\{d\hat{\phi}(u,v);\,
(u,v)\in T_{(x,\nu)}S(D,M,g),
\,
d\hat{\pi}(u,v)=\lambda \nu
~\text{for some}~
\lambda\in \mathbb{R}\big\}
\\
\qquad
=\big\{d\hat{\phi}(u,v); (u,v)\in T_{(x,\nu)}S(D,M,g),
\;
d\phi\circ d\hat{\pi}(u)=d\hat{\pi}'\circ d\hat{\phi}(u,v)=\lambda d\phi(\nu)
, \,
\lambda\in \mathbb{R}\big\}
\\
\qquad
=\big\{(\phi_*(D))^{[1]}\big\}_{\hat{\phi}(x,\nu)}.
\end{gather*}
This ends the proof of~(ii).

Assume now that $\phi$ is a~dif\/feomorphism such that $\phi_*(D)=D'$ and let $\psi=\phi^{-1}$.
According to the def\/inition of $\hat{\phi}$ and $\hat{\psi}$, it follows trivially that $\hat{\psi}\circ\hat{\phi}={\rm Id}$.
Besides, based on the def\/inition of $[\phi_*(D)]^{[1]}$, and since $d_x\phi $ is an isomorphism, we must have $
\big\{(\phi_*(D))^{[1]}\big\}_{\hat{\phi}(x,\nu)}=\big\{(D')^{[1]}\big\}_{\hat{\phi}(x,\nu)}$.
Finally, since $\phi$ is an isometry from $(M,g)$ to $(M',g')$, $d\phi$ is also an isometry for $(TM, \bar{g})$ and
$(TM',\bar{g}')$ if $\bar{g}$ and $\bar{g}'$ are the canonical Riemannian metrics on the tangent bundles induced by~$g$
and $g'$, respectively.
As $\hat{g}$ and $\hat{g}'$ are the restrictions of $\bar{g}$ and $\bar{g}'$ to $S(D,M,g)\subset TM$ and
$S(D',M',g')\subset TM'$, respectively, we obtain the last property and the proof of the lemma.
\end{proof}

Consequently, as in the context of Cartan prolongation, for any $m\geq 2$ and $k\geq 1$ we can inductively def\/ine a~tower of
sphere bundles (for a~f\/ixed choice of the metric~$g$ on a~manifold~$M$) as
\begin{gather}\label{towerS}
\hat{P}^k(M)\rightarrow \hat{P}^{k-1}(M)\rightarrow\dots\rightarrow\hat{P}^j(M)\rightarrow
\hat{P}^{j-1}(M)\rightarrow\dots\rightarrow \hat{P}^{1}(M)\rightarrow \hat{P}^0(M):=M,\!\!\!
\end{gather}
where $\hat{P}^j(M)$ is a~manifold of dimension $(j+1)m+1$ for any $j=0,\dots, k$, and on each $\hat{P}^j(M)$ we have
a~canonical distribution $\hat{\Delta}_j$ and a~Riemannian metric $g_j$. All these data are inductively def\/ined by:
\begin{itemize}\itemsep=0pt
\item $g_0=g$ is a~given Riemannian metric on $\hat{P}^0(M)=M$, $\hat{\Delta}_0=TM$,
\item for $1\leq j\leq k$:
\begin{itemize}\itemsep=0pt
\item[$\circ$] $\hat{P}^j(M)=S\big(\hat{\Delta}_{j-1},\hat{P}^{j-1}(M),g_{j-1}\big)$,
\item[$\circ$] $\hat{\Delta}_j=\big(\hat{\Delta}_{j-1}\big)^{[1]}$,
\item[$\circ$] $g_j$ is the Riemannian metric $\hat{g}_{j-1}$ on $S\big(\hat{\Delta}_{j-1},\hat{P}^{j-1}(M),g_{j-1}\big)$ associated
with $g_{j-1}$.
\end{itemize}
\end{itemize}

Note that, if $g'$ is another Riemannian metric on~$M$, according to Lemmas~\ref{changeg} and~\ref{prolongmap}, we construct, by
induction, a~family of dif\/feomorphisms $\psi^j$ such that, if
\begin{gather*}
\hat{P}'^k(M)\rightarrow \hat{P}'^{k-1}(M)\rightarrow\dots\rightarrow \hat{P}'^{1}(M)\rightarrow \hat{P}'^0(M):=M
\end{gather*}
is the tower of sphere bundles associated with the chosen metric $g'$ on~$M$ we obtain, for all $j=0,\dots, k$:
\begin{itemize}\itemsep=0pt
\item $\psi^j(\hat{P}^j(M))=\hat{P}'^j(M)$,
\item $\psi^j$ is f\/iber preserving,
\item $\psi^j_*(\hat{\Delta}_j)=\hat{\Delta}'_j$.
\end{itemize}
Therefore, {\it the properties which characterize the tower~\eqref{towerS} are independent of the choice of
the Riemannian metric~$g$ on~$M$}.

For the sake of simplicity we write $\hat{P}^j(m): =\hat{P}^j(\mathbb{R}^{m+1})$ for any $j\in \mathbb{N}$.
From Theorem~\ref{modelP} and Lemma~\ref{spprolong} we can deduce the following result:
\begin{The}\label{modelS}
Consider the tower of sphere bundles
\begin{gather}\label{towerS0}
\hat{P}^k(m)\rightarrow \hat{P}^{k-1}(m)\rightarrow\dots\rightarrow \hat{P}^j(m)\rightarrow
\hat{P}^{j-1}(m)\rightarrow\dots\rightarrow\hat{P}^{1}(m)\rightarrow\hat{P}^0(m):=\mathbb{R}^{m+1}
\end{gather}
associated with the canonical metric on $\mathbb{R}^{m+1}$.
Then the following properties hold:
\begin{enumerate}\itemsep=0pt
\item[$1.$] There exists a~canonical two-fold covering $\tau^j: \hat{P}^j(m)\rightarrow P^j(m)$ such that
$\tau^j(\hat{\Delta}_j)=\Delta_j$ for all $j=1,\dots,k$.
\item[$2.$] On $ \hat{P}^j(m)$, the distribution $\hat{\Delta}_j$ generates a~special multi-flag of step~$m$ and length~$j$,
for all $j=1,\dots,k$.
\item[$3.$] Let $\mathbb{D}:D=D_{k}\subset D_{k-1}\subset\dots\subset D_{j}\subset\dots\subset D_{1} \subset D_{0}= TM$ be
a~special multi-flag of step $m\geq 2$ and length $k\geq 1$.
Then, for any $x\in M$, there exists $y\in \hat{P}^k(m)$ for which the differential system
$(\hat{P}^k(m),\hat{\Delta}_k,y)$ is locally equivalent to the differential system $(M,D,x) $.
\end{enumerate}
\end{The}
This tower~\eqref{towerS0} will be called the spherical tower of special multi-f\/lags of step~$k$.

\section{Tower of sphere bundles associated with a~kinematic system}\label{tower}

\subsection{A kinematic system for special multi-f\/lags}\label{arm}

We set ourselves in the context of~\cite{LR, SP1}.
Consider, in $\mathbb{R}^{m+1}$, a~family of~$k$ segments $[M_{i},M_{i+1}]$, where $i=0,\dots,k-1$ and $m\geq 2$, keeping
a~constant length $l_i=1$ between $M_i$ and $M_{i+1}$, with articulation at points $M_i$, for $i=1,\dots,k-1$.

\begin{figure}[t]\centering
\includegraphics{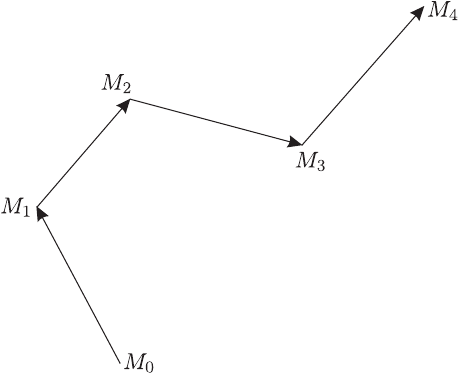}
\caption{Articulated arm of length $4$.}\label{arm4}
\end{figure}

Such a~system is called a~``$k$-bar system'' in~\cite{LR} and an ``articulated arm of length~$k$'' in~\cite{SP1}.
The kinematic evolution of the extremity $M_0$, {\it under the constraint that the velocity of each point $M_i$,
$i=0,\dots,k-1$, is colinear with the segment $[M_i,M_{i+1}]$}, is completely described in terms of hyperspherical
coordinates in~\cite{SP1}, whereas results of f\/latness and controllability for such a~system are proved in~\cite{LR} (see
Fig.~\ref{arm4}).

A~special multi-f\/lag of step $m\geq 2$ and length $k\geq 1$ is associated with this kinematic system as explained in the following.
We can decompose $(\mathbb{R}^{m+1})^{k+1}$
into the product $\mathbb{R}^{m+1}_0\times\dots\times\mathbb{R}^{m+1}_i\times\dots\times\mathbb{R}^{m+1}_{k}$.
Let $x_i=(x_i^1,\dots,x_i^{m+1})$ be the canonical coordinates in the space $\mathbb{R}^{m+1}_i$ which is equipped with
its canonical scalar product $\langle\cdot,\cdot\rangle$.
The space $(\mathbb{R}^{m+1})^{k+1}$ is then equipped with its canonical scalar product too.

 Consider an articulated arm of length~$k$ denoted by $(M_0,\dots, M_{k})$.
We def\/ine, on $(\mathbb{R}^{m+1})^{k+1}$, the vector f\/ields
\begin{gather*}
{\mathcal Z}_i = \sum\limits_{r=1}^{m+1} \big(x_{i+1}^r-x_{i}^r\big)\frac{\partial}{\partial x_{i}^r}
\qquad
\text{for}
\quad
i=0,\dots,k-1
\end{gather*}
Based on our previous assumptions, the kinematic evolution of the articulated arm is described by the controlled system
\begin{gather*}
\dot{q}= \sum\limits_{i=0}^{k-1} u_i{\mathcal Z}_i+ \sum\limits_{r=1}^{m+1}
u_{n+r} \frac{\partial}{\partial x_{k}^r}
\end{gather*}
with the following constraints:
$||x_{i}-x_{i+1}||=1$ for $i=0,\dots, k-1$ (see~\cite{LR} or~\cite{SP1}).

For $k\geq 1$, the {\it configuration space} ${\mathcal C}^k(m)$ is the set
\begin{gather*}
\{(x_0,\dots,x_{k})\,|\,   \Psi_i (x_0,\dots,x_{k})=0,\;\forall \, i=0,\dots, k-1 \},
\end{gather*}
where $\Psi_i (x_0,\dots,x_{k})=||x_{i}-x_{i+1}||^2-1$ and we set ${\mathcal C}^0(m):=\mathbb{R}^{m+1}$.

For $i=0,\dots,k-1$, the vector f\/ield:
\begin{gather*}
{\mathcal N}_i= \sum\limits_{r=1}^{m+1} \big(x_{i+1}^r-x_{i}^r\big)\left(\frac{\partial}{\partial x_{i+1}^r}-
\frac{\partial}{\partial x_{i}^r}\right)
\end{gather*}
is proportional to the gradient of $\Psi_i$.\footnote{In fact, we could use the dif\/ferential $d\Psi_i$ instead of
${\mathcal N}_i$; however, this choice is motivated by the use of the projection $\Pi_k$ whose kernel is generated by
$\{{\mathcal N}_i;\;i=0,\dots,k-1\}$}
It follows that the tangent space $T_q{\mathcal C}^k(m)$ is the subspace of $T_q(\mathbb{R}^{m+1})^{k+1}$ which is
orthogonal to ${\mathcal N}_i(q)$ for $i=0,\dots,k-1$.

Denote by ${\mathcal E}_k$ the distribution generated by the family of vector f\/ields
\begin{gather*}
\left\{{\mathcal Z}_0,\dots,{\mathcal Z}_{k-1}, \frac{\partial}{\partial
x_{k}^1},\dots, \frac{\partial}{\partial x_{k}^{m+1}}\right\}.
\end{gather*}
Let ${\mathcal D}_k$ be the distribution on ${\mathcal C}^k(m)$ def\/ined by ${\mathcal D}_k(q)=T_q{\mathcal C}^k(m)\cap
{\mathcal E}_k$.
Thus we have:
\begin{Lem}[see~\protect{\cite{SP1}}]
\label{distribinduite}
${\mathcal D}_k$ is a~distribution of constant rank $m+1$ generated by
\begin{gather*}
\big(x_{1}^r-x_{0}^r\big){\mathcal Z}_0+ \frac{\partial}{\partial x_{1}^r}
\qquad
\text{for}
\quad
r=1,\dots, m+1
\ \
\text{if}
\ \
k=1,
\\
\big(x_{k}^r-x_{k-1}^r\big)\!\left( \sum\limits_{i=0}^{k-2}
 \prod\limits_{j=i+1}^{k-1} {\mathcal A}_j{\mathcal Z}_i+{\mathcal Z}_{k-1}\right)\!
+ \frac{\partial}{\partial x_{k}^r}
\qquad
\text{for}
\quad
r=1, \dots, m+1
\
\text{if}
\
k\geq 2,
\end{gather*}
where ${\mathcal A}_j(q)=-\langle{\mathcal N}_{j}(q),{\mathcal N}_{j-1}(q)\rangle=\langle{\mathcal Z}_{j}(q),{\mathcal N}_{j-1}(q)\rangle$,
for $j=1,\dots, k-1$.
\end{Lem}
\begin{Not}
According to notations of Lemma~\ref{distribinduite}, we def\/ine
\begin{itemize}\itemsep=0pt
\item $Y_1={\mathcal Z}_0$ and $Y_0=0$,
\item for $k\geq 2$ $Y_k=\left( \sum\limits_{i=0}^{k-2} \prod\limits_{j=i+1}^{k-1} {\mathcal
A}_j {\mathcal Z}_i\right)+{\mathcal Z}_{k-1}$.
\end{itemize}
\end{Not}
\begin{Rem}\label{Y}
According to the previous notations,  the inductive relation holds:
\begin{gather}\label{recY}
Y_k={\mathcal A}_{k-1}Y_{k-1}+{\mathcal Z}_{k-1}.
\end{gather}
Consequently the distribution ${\mathcal D}_k$ is generated by the family
\begin{gather*}
\left\{\big(x_{k}^r-x_{k-1}^r\big)Y_k+ \frac{\partial}{\partial x_{k}^r}\,\Big|\, r=1,\dots, m+1\right\}.
\end{gather*}
\end{Rem}

The properties of ${\mathcal D}_k$ are summarized in the following result (see~\cite{LR, SP1}):
\begin{The}\label{drap}
The distribution ${\mathcal D}_k$ on ${\mathcal C}^k(m)$ satisfies the following properties:
\begin{enumerate}\itemsep=0pt
\item[$1)$] ${\mathcal D}_k$ is a~distribution of constant rank $m+1$,
\item[$2)$] ${\mathcal D}_k$ generates a~special multi-flag on ${\mathcal C}^k(m)$ of step~$m$ and length~$k$.
\end{enumerate}
\end{The}

\subsection{Articulated arm and spherical prolongation}\label{armsph}

The following canonical
tower of sphere bundles
\begin{gather}\label{towerC}
{\mathcal C}^k(m)\rightarrow \mathcal{C}^{k-1}(m)\rightarrow\dots\rightarrow \mathcal{C}^1(m)\rightarrow
\mathcal{C}^0(m):=\mathbb{R}^{m+1}
\end{gather}
is associated with an articulated arm on ${\mathbb R}^{m+1}$ ($m \geq 2$) of length $k \geq 1$,
where, for $j=1,\dots, k$, the projection ${\mathcal C}^j(m)\rightarrow \mathcal{C}^{j-1}(m)$ is the restriction of the
canonical projection
\begin{gather*}
\begin{matrix}
\mathbb{R}^{m+1}_0\times\dots\times\mathbb{R}^{m+1}_i\times\dots\times\mathbb{R}^{m+1}_{j}
&
\rightarrow
&
\mathbb{R}^{m+1}_0\times\dots\times\mathbb{R}^{m+1}_i\times\dots\times\mathbb{R}^{m+1}_{j-1}
\\
(x_0,\dots,x_{j-1},x_j) & \mapsto & (x_0,\dots,x_{j-1})
\end{matrix}
\end{gather*}

According to Theorems~\ref{drap} and~\ref{modelS}, we know that the dif\/ferential system
$(\mathcal{C}^k(m),\mathcal{D}_k))$ associated with an articulated arm of length~$k$ on $\mathbb{R}^{m+1}$ is locally
isomorphic to the canonical dif\/ferential system $(\hat{P}^k(m),\hat{\Delta}_k)$ at some appropriate points.
In fact, obtain a stronger result  (as stated in Theorem~\ref{1}(1)) stated as follows:

\begin{The}\label{artmod}
For each $m\geq 2$ and $1\leq j\leq k$, there exists a~diffeomorphism $F^j$ from $\hat{P}^j(m)$ on $\mathcal{C}^j(m)$
such that:
\begin{itemize}\itemsep=0pt
\item[$(i)$] $\rho^j\circ F^j =F^{j-1} \circ \hat{\pi}^j$,
where $\hat{\pi}^j: \hat{P}^j(m)\rightarrow \hat{P}^{j-1}(m)$
and
$\rho^j: \mathcal{C}^j(m)\rightarrow\mathcal{C}^{j-1}(m)$ are the canonical projections, and
\item[$(ii)$] $F^j_*(\hat{\Delta}_j)=\mathcal{D}_j$.
\end{itemize}
\end{The}

Therefore, according to Theorems~\ref{modelS} and~\ref{artmod} we can obtain:
\begin{The}
Let $\mathbb{D}:D=D_{k}\subset D_{k-1}\subset\dots\subset D_{j}\subset\dots\subset D_{1} \subset D_{0}= TM$ be a~special
multi-flag of step $m\geq 2$ and length $k\geq 1$.
Then, for any $x\in M$, there exists $y\in {\mathcal C}^k(m)$ for which the differential system
$(\mathcal{C}^k(m),\mathcal{D}_k,y)$ is locally equivalent to the differential system $(M,D,x) $.
\end{The}

The end of this subsection is devoted to the proof of Theorem~\ref{artmod} and thus the proof of Theorem~\ref{1}(1).
Before doing so, we need some auxiliary results.
\begin{Lem}\label{versionZ}
For $k\geq 1$, consider the natural decomposition:
\begin{gather*}
\big[T(\mathbb{R}^{m+1})^{k+1}\big]\Big|_{{\mathcal C}^k(m)}=T{\mathcal C}^k(m)\oplus\big[T{\mathcal C}^k(m)\big]^{\perp}
\quad
\text{on}
\quad
{\mathcal C}^k(m),
\end{gather*}
where $[T{\mathcal C}^k(m)]^{\perp}$ is the orthogonal of $T{\mathcal C}^k(m)$, and denote by $\Pi_k$ the orthogonal
projection of $[T(\mathbb{R}^{m+1})^{k+1}]_{|{\mathcal C}^k(m)}$ onto $T{\mathcal C}^k(m)$.
Finally, denote by ${\mathcal L}_k$ the vertical bundle defined by the natural fibration of ${\mathcal C}^k(m)$ over
${\mathcal C}^{k-1}(m)$.
We then obtain the following:
\begin{enumerate}\itemsep=0pt
\item[$1.$] The family of vector fields $\left\{\Pi_k\left( \frac{\partial}{\partial x_{k}^r}\right), r=1,\dots,m+1\right\}$
generates ${\mathcal L}_k$.
\item[$2.$] The distribution ${\mathcal D}_k$, is generated by ${\mathcal L}_k$ and the vector field $X_k=Y_k+{\mathcal V}_k$,
where
\begin{gather*}
{\mathcal V}_k= \sum\limits_{s=1}^{m+1}\big(x_k^s-x_{k-1}^s\big) \frac{\partial}{\partial x_{k}^s}.
\end{gather*}
\item[$3.$] The distribution ${\mathcal D}_k$ is also generated by the family of vector fields
\begin{gather*}
\big(x_{k}^r-x_{k-1}^r\big)X_k+\Pi_k\left( \frac{\partial}{\partial x_{k}^r}\right), \qquad r=1,\dots,m+1.
\end{gather*}
\end{enumerate}
\end{Lem}
\begin{proof}
Denote by $\mathfrak{H}_k$ the subdistribution of ${\mathcal E}_k$ generated by the family of vector f\/ields
\begin{gather*}
\left\{ \frac{\partial}{\partial x_{k}^r} \,\Big|\,  r=1,\dots, m+1\right\}.
\end{gather*}
It follows that $\mathfrak{H}_k\bigcap T{\mathcal C}^k(m)$ is a~distribution on ${\mathcal C}^k(m)$ which is contained
in ${\mathcal D}_k$.
In fact, we have
\begin{gather*}
{\mathcal L}_k=\ker d\Psi_{k-1}\bigcap \mathfrak{H}_k=\Pi_k(\mathfrak{H}_k).
\end{gather*}
Therefore, the distribution ${\mathcal L}_k$ is spanned by the family of vector f\/ields
\begin{gather*}
\left\{\Pi_k\left( \frac{\partial}{\partial x_{k}^r}\right) \,\Big|\, r=1,\dots, m+1\right\}.
\end{gather*}
On the other hand, $\mathfrak{H}_k$ is the vertical bundle of the canonical projection,
\begin{gather}\label{canoproj}
\begin{matrix}
\mathbb{R}^{m+1}_0\times\dots\times\mathbb{R}^{m+1}_{k} & \rightarrow & \mathbb{R}^{m+1}_0\times\dots\times\mathbb{R}^{m+1}_{k-1}
\\
(x_0,\dots,x_{k-1},x_k) & \mapsto & (x_0,\dots,x_{k-1})
\end{matrix}
\end{gather}
It follows that ${\mathcal L}_k$ is the vertical bundle of the induced projection of ${\mathcal C}^k(m)$ onto
${\mathcal C}^{k-1}(m)$.
Moreover, the f\/iber over $q\in {C}^{k-1}(m)$ of the previous f\/ibration is the unit sphere $S_q=\{(q,x_k)\, |\,
\Psi_{k-1}(q,x_k)=0\}$, which proves~(1).

Furthermore, the vector f\/ield ${\mathcal
V}_k= \sum\limits_{s=1}^{m+1}\left(x_k^s-x_{k-1}^s\right) \frac{\partial}{\partial x_{k}^s}$ is
vertical for the projection~\eqref{canoproj} and is orthogonal to each $S_q$.
Since $||{\mathcal V}_k||=1$, we thus have
\begin{gather}\label{projfibre}
\Pi_k\left( \frac{\partial}{\partial x_{k}^r}\right)=  \frac{\partial}{\partial x_{k}^r}-\big(x_k^r-x_{k-1}^r\big){\mathcal V}_k.
\end{gather}

From Remark~\ref{Y} the distribution ${\mathcal D}_k$ is generated by the family
\begin{gather*}
\left\{\big(x_{k}^r-x_{k-1}^r\big)Y_k+ \frac{\partial}{\partial x_{k}^r}  \,\Big|\, r=1,\dots, m+1\right\}.
\end{gather*}
Therefore the vector f\/ield
$X_k=Y_k+ \sum\limits_{r=1}^{m+1}(x_{k}^r-x_{k-1}^r) \frac{\partial}{\partial x_{k}^r}$ is
tangent to ${\mathcal D}_k$, but clearly this vector f\/ield is not tangent to ${\mathcal L}_k$.
Since ${\mathcal D}_k$ is a~distribution of constant rank $m+1$ and ${\mathcal L}_k$ is an (integrable) subdistribution of
rank~$m$, then ${\mathcal D}_k$ is generated by ${\mathcal L}_k$ and $X_k$, which proves~(2).

Finally, according to the relation~\eqref{projfibre}, each vector f\/ield
$(x_{k}^r-x_{k-1}^r)Y_k+ \frac{\partial}{\partial x_{k}^r}$ can be written~as
\begin{gather*}
\big(x_k^r-x_{k-1}^r\big)X_k+\Pi_k\left( \frac{\partial}{\partial x_{k}^r}\right).
\end{gather*}
Therefore from Lemma~\ref{distribinduite}, we obtain~(3).
\end{proof}

\begin{Pro}\label{sprolongarm}\quad
\begin{enumerate}\itemsep=0pt
\item[$1.$] There exists a~bundle isomorphism $\hat{\Psi}^k: {\mathcal D}_k\rightarrow {\mathcal C}^k(m)
\times
\mathbb{R}^{m+1}$ over the identity of ${\mathcal C}^k(m)$.

\item[$2.$] Let $\gamma_{k}$ be the Riemannian metric on the bundle ${\mathcal D}_k$ such that the morphism $\hat{\Psi}^k$ is an
isometry between ${\mathcal D}_k$ and ${\mathcal C}^k(m)\times\mathbb{R}^{m+1}$, where the bundle ${\mathcal C}^k(m)\times\mathbb{R}^{m+1}$
is given by the canonical Riemannian metric induced by the canonical Euclidean metric on each fiber $\{q\}\times\mathbb{R}^{m+1}$.
Then $\hat{\Psi}^k$ induces a~diffeomorphism $\Psi^k: S(\mathcal{D}_k,\mathcal{C}^k(m),\gamma_k)\rightarrow {\mathcal C}^{k+1}(m)$
which fulfills the following properties:
\begin{itemize}\itemsep=0pt
\item[$(i)$] $\Psi^k$ commutes with the canonical projections $S({\mathcal D}_k,\mathcal{C}^k(m), \gamma_k)
\rightarrow
\mathcal{C}^k(m)$ and $\mathcal{C}^{k+1}(m)$ $\rightarrow \mathcal{C}^k(m)$,
\item[$(ii)$] $\Psi^k_*[(\mathcal{D}_k)^{[1]}]=\mathcal{D}_{k+1}$.
\end{itemize}
\end{enumerate}
\end{Pro}

\begin{proof}
From Remark~\ref{Y}, the bundle ${\mathcal D}_k$ has $m+1$ non-zero global sections
\begin{gather*}
W_r=\big(x_{k}^r-x_{k-1}^r\big)Y_k+ \frac{\partial}{\partial x_{k}^r}
\qquad
\text{for}
\quad
r=1,\dots, m+1.
\end{gather*}
Thus ${\mathcal D}_k$ is a~trivial bundle.
This global trivialization gives rise to an isomorphism $\hat{\Psi}^k:\mathcal{D}_k \rightarrow \mathcal{C}^k(m)\times\mathbb{R}^{m+1}$
characterized by
\begin{gather*}
\hat{\Psi}^{k}\left(x_0,\dots,x_k, \sum\limits_{r=1}^{m+1}\nu^r W_r(x_0,\dots,x_k)\right)=(x_0,\dots,x_k,\nu),
\end{gather*}
where $\nu=(\nu^1,\dots,\nu^{m+1})$ and so $\hat{\Psi}^k$ satisf\/ies point (1).

Put on ${\mathcal D}_k$ the Riemannian metric ${\gamma}_k=\hat{\Psi}^*g$, where~$g$ is the canonical Euclidean metric on the
trivial bundle ${\mathcal C}^k(m)\times\mathbb{R}^{m+1}$.
It follows that the global basis $\left\{{W}_r\,|\,   r=1,\dots,m+1\right\}$ is an orthonormal family, and then the set
$S({\mathcal D}_k, {\mathcal C}^k(m),{\gamma}_k)$ is:
\begin{gather*}
\left\{\!\left(x_0,\dots,x_k, \sum\limits_{r=1}^{m+1}\nu^r W_r(x_0,\dots,x_k)\right) \!\Big|\, (x_0,\dots,x_k)\in {\mathcal C}^k(m),\;
\big(\nu^1\big)^2+\dots+\big(\nu^{m+1}\big)^2=1\right\},
\end{gather*}
and the f\/iber over $(x_0,\dots,x_k)\in {\mathcal C}^k(m)$ is equal to
\begin{gather*}
\left\{ \sum\limits_{r=1}^{m+1}\nu^r {W}_r(x_0,\dots,x_k) \,\Big|\, \big(\nu^1\big)^2+\dots+\big(\nu^{m+1}\big)^2=1\right\}.
\end{gather*}
From the choice of the metric ${\gamma}_k$, the restriction $\overline{\Psi}^k$ of $\hat{\Psi}^k$ to $S({\mathcal D}_k,
{\mathcal C}^k(m),{\gamma}_k)$ is then a~dif\-feomorphism
onto ${\mathcal C}^{k}(m)\times\mathbb{S}^m$.
Moreover, by construction, we have the following commutative diagram:
\begin{gather*}
\begin{matrix}
\hspace*{2mm}\overline{\Psi}^k 			
\\
S({\mathcal D}_k, {\mathcal C}^k(m),{\gamma}_k) \xrightarrow{\hspace*{1cm}} {\mathcal C}^{k}(m)
\times
\mathbb{S}^m
\qquad
\hfill
\\
\vspace*{0.3cm}\hfill \pi
\
\bigg\downarrow \hfill
\qquad
\hfill\bigg\downarrow
\
p_{{\mathcal C}^k}\hfill
\\
\hfill \mathcal{C}^{k}(m) \xrightarrow{\hspace*{1.2cm}}  	{\mathcal C}^{k}(m)\hfill
\\
\hspace*{-2mm} {\rm Id}
\end{matrix}
\end{gather*}
where the vertical maps are the canonical projections.
Consider now the map
\begin{gather*}
\widehat{\mathcal T}: \ {\mathcal C}^k(m)
\times
\mathbb{R}^{m+1}\rightarrow{\mathcal C}^k(m)
\times
\mathbb{R}^{m+1}
\end{gather*}
def\/ined by $\widehat{\mathcal T}(x_0,\dots,x_k,z)=(x_0,\dots,x_k,x_k+z)$.
This map induces a~dif\/feomorphism
\begin{gather*}
{\mathcal T}: \ {\mathcal C}^k(m)
\times
\mathbb{S}^m \rightarrow {\mathcal C}^{k+1}(m).
\end{gather*}
We set $ \Psi^k={\mathcal T}\circ\bar{\Psi}^k$
and obtain a~dif\/feomorphism $\Psi^k: S({\mathcal D}_k, {\mathcal C}^k(m),{\gamma}_k) \rightarrow {\mathcal
C}^{k+1}(m)$ which commutes with the canonical projections of the sphere bundles $\pi: S({\mathcal D}_k, {\mathcal
C}^k(m),{\gamma}_k) \rightarrow {\mathcal C}^k(m)$ and $\rho^{k+1} : {\mathcal C}^{k+1}(m)\rightarrow {\mathcal C}^k(m)$.

Given a~vector f\/ield~$W$ on ${\mathcal C}^{k}(m)$, for any lift $\overline{W}$ of~$W$ on ${\mathcal C}^{k+1}(m)$,
the vector f\/ield ${\mathcal W}=(T\Psi^k)^{-1}(\overline{W})$ is a~vector f\/ield on the manifold ${\mathcal D}_k$ such
that $T\pi({\mathcal W})=W$.
Moreover, $T\Psi^k$ sends $\ker (T\pi)$ onto $\ker (T\rho^{k+1})$.
Denote by $\overline{W}_r$ the vector f\/ield on ${\mathcal C}^{k}(m)
\times
{\mathbb{R}}^{m+1}$ def\/ined by
\begin{gather*}
\overline{W}_r(x_0,\dots,x_k,x_{k+1})=\big(x_{k}^r-x_{k-1}^r\big)Y_k+ \frac{\partial}{\partial x_{k}^r}
+ \frac{\partial}{\partial x_{k+1}^r}.
\end{gather*}
Notice that $\overline{W}_r(x_0,\dots,x_k,x_{k+1})$ is actually tangent to ${\mathcal C}^{k+1}(m)$, and we have
$Tp_{{\mathcal C}^k}(\overline{W}_r)={W}_r$.

We set ${\mathcal W}_r=(T\Psi^k)^{-1}(\overline{W}_r)$.

Now, given a~point $(x_0,\dots,x_k)\in{\mathcal C}^k(m)$, the spherical prolongation $({\mathcal D}_k)^{[1]}$ of
${\mathcal D}_k$ at a~point $\Big(x_0,\dots,x_k, \sum\limits_{r=1}^{m+1}\nu^r W_r(x_0,\dots,x_k)\Big)$ is generated
by $ \sum\limits_{r=1}^{m+1}\nu^r {\mathcal W}_r(x_0,\dots,x_k)$ and by the tangent space to the f\/iber
through this point.
This implies that the space $\Psi^k_*\big(({\mathcal D}_k)^{[1]}\big)$ over the point $(x_0,\dots,x_k,x_k+\nu)$, is generated by
$ \sum\limits_{r=1}^{m+1}\nu^r \overline {W}_r(x_0,\dots,x_k,x_k+\nu)$ and by the tangent space of the f\/iber
at this point.

But, according to Remark~\ref{Y}, this vector f\/ield $ \sum\limits_{r=1}^{m+1}\nu^r\overline
{W}_r(x_0,\dots,x_k,x_{k}+\nu)$ can be written~as
\begin{gather*}
 \sum\limits_{r=1}^{m+1}\left\{\nu^r\big(x_{k}^r-x_{k-1}^r\big)Y_k+\nu^r\left( \frac{\partial}{\partial
x_{k}^r}+ \frac{\partial}{\partial x_{k+1}^r}\right)\right \} ={\mathcal A}_k Y_k+{\mathcal Z}_k+{\mathcal
V}_{k+1},
\end{gather*}
where $x_{k+1}=x_k+\nu$.
According to~\eqref{recY} and Lemma~\ref{versionZ}(2), this last expression is exactly $X_{k+1}$.
We deduce that the space $\Psi^k_*\big(({\mathcal D}_k)^{[1]}\big)$ is generated, over the point $(x_0,\dots,x_k,x_{k+1})$,  by
$X_{k+1}$ and by the tangent space to the f\/iber at this point.
 Lemma~\ref{versionZ}(2) implies that $\Psi^k_*\big(({\mathcal D}_k)^{[1]}\big)$ and~${\mathcal D}_{k+1}$ coincide at
point $(x_0,\dots,x_k,x_{k+1})$.
\end{proof}

\begin{proof}[Proof of Theorem~\ref{artmod}] On the one hand, we have $\hat{P}^0(m)={\mathcal C}^0(m)=\mathbb{R}^{m+1}$ with
$\hat{\Delta}_0={\mathcal D}_0=T\mathbb{R}^{m+1}$ and, on the other hand, if $g_0$ is the canonical Riemannian metric on
$T\mathbb{R}^{m+1}$, we got $P^1(m)=S({\mathcal D}_0,{\mathcal C}^0(m),g_0)$ with $\hat{\Delta}_1=({\mathcal
D}_0)^{[1]}$.
Therefore, the result comes from Proposition~\ref{sprolongarm} for $k=1$.

Assume that there exists a~dif\/feomorphism $F^j:\hat{P}^j(m)\rightarrow {\mathcal C}^j(m)$ which satisf\/ies properties (i) and
(ii) of Theorem~\ref{artmod}.

From Proposition~\ref{sprolongarm}, we obtain a~dif\/feomorphism $\Psi^j: S({\mathcal D}_j, {\mathcal
C}^j(m),\gamma_{j})\rightarrow {\mathcal C}^{j+1}(m)$ such that $\Psi_*[({\mathcal D}_j)^{[1]}]={\mathcal D}_{j+1}$ and
satisfying the following commutative diagram
\begin{gather*}
\begin{matrix}
\hspace*{5mm} {\Psi}^j 			
\\
S({\mathcal D}_j, {\mathcal C}^j(m),{\gamma}_j) \xrightarrow{\hspace*{1cm}} {\mathcal C}^{j+1}(m)
\qquad
\hfill
\\
\vspace*{0.3cm}\hfill
\qquad
\bigg\downarrow \hfill		 		\hfill\bigg\downarrow \hfill
\\
\hfill
\qquad
\mathcal{C}^{j}(m)  \xrightarrow{\hspace*{1.3cm}}  	{\mathcal C}^{j}(m)\hfill
\\
\hspace*{3mm} {\rm Id}
\end{matrix}
\end{gather*}

According to previous induction at level~$j$, we  put on $\hat{P}^{j}(m)$ the Riemannian metric $\bar{\gamma}_j=(\Psi^j)^*(\gamma_{j})$.
From Lemma~\ref{prolongmap}, we  extend $F^j:\hat{P}^j(m)\rightarrow {\mathcal C}^j(m)$ to a~dif\/feomorphism
\begin{gather*}
\hat{\Theta}^{j}: \ S\big(\hat{\Delta}_j,\hat{P}^{j}(m),\bar{\gamma}_j\big)\rightarrow S\big({\mathcal D}_j,{\mathcal
C}^{j}(m),\gamma_j\big)
\end{gather*}
such that $\hat{\Theta}^{j}_*[(\hat{\Delta}_j)^{[1]}]=({\mathcal D}_j)^{[1]}$ and satisfying the following commutative diagram
\begin{gather*}
\begin{matrix}
\hspace*{-2mm}
\hat{\Theta}^j 			
\\
S\big(\hat{\Delta}_j, \hat{\mathcal P}^j(m),\bar{\gamma}_j\big) \xrightarrow{\hspace*{0.9cm}}
S\big({\mathcal D}_j,{\mathcal C}^j(m),\gamma_j\big)
\qquad
\hfill
\\
\vspace*{0.3cm}\hfill
\qquad
\bigg\downarrow \hfill		 		\hfill\bigg\downarrow
\qquad
\hfill
\\
\hfill \hat{P}^{j}(m)  \xrightarrow{\hspace*{1.3cm}}  	{\mathcal C}^{j}(m)\hfill
\\
\hspace*{-2mm}
F^j
\end{matrix}
\end{gather*}

We put on $\hat{P}^{j}(m)$ the Riemannian metric obtained by successive induction on the tower bundle~\eqref{towerS}
(see Section~\ref{cartsph} just after~\eqref{towerS}).
According to Lemma~\ref{changeg}, we also obtain a~dif\/feo\-mor\-phism $\Theta: \hat{P}^{j+1}(m)\rightarrow
S(\hat{\Delta}_j,\hat{P}^{j}(m),\bar{\gamma}_j)$ such that $\Theta_*(\hat{\Delta}_{j+1})=\hat{\Delta}_j^{[1]}$ and satisfying
the following commutative diagram
\begin{gather*}
\begin{matrix}
\hspace*{-5mm}\Theta
\\
\hfill
\qquad
\hat{P}^{j+1}(m) \xrightarrow{\hspace*{1cm}} S(\hat{\Delta}_j,\hat{P}^j(m), \hat{\gamma}_j)
\\
\vspace*{0.3cm}\hfill	\bigg\downarrow \hfill
\qquad
 		 \hfill\bigg\downarrow
\qquad
\hfill
\\
\hfill\hat{P}^{j}(m)  \xrightarrow{\hspace*{1.2cm}}  	\hat{P}^{j}(m)\hfill
\qquad
\\
\hspace*{-5mm}{\rm Id}
\end{matrix}
\end{gather*}

If we juxtapose the three last diagrams, we obtain the required dif\/feomorphism  $F^{j+1}=\Psi^j\circ \hat{\Theta}^j\circ\Theta$.
\end{proof}

\begin{Rem}
According to Theorems~\ref{modelS} and~\ref{artmod}, from towers~\eqref{towerS0} and~\eqref{towerC}, we obtain the
following diagram in which each vertical map is a~$2$-fold covering for $k\geq 1$:
\begin{gather*}
\begin{matrix}
& \hfill{\mathcal C}^k(m)\xrightarrow{\hspace*{0.4cm}}
\mathcal{C}^{k-1}(m)\xrightarrow{\hspace*{0.4cm}}\cdots\xrightarrow{\hspace*{0.4cm}}
\mathcal{C}^1(m)\xrightarrow{\hspace*{0.4cm}}\mathcal{C}^0(m):=\mathbb{R}^{m+1}\hfill
\\
\vspace*{0.15cm}&\hfill\!\!\!\!\!\!\!\!\!\Big\downarrow
\quad
\hfill\Big\downarrow
\hfill\ \ \cdots\hfill\Big\downarrow
\quad
\hfill\Big\downarrow
\qquad
\hfill
\\
& {P}^k(m)\xrightarrow{\hspace*{0.4cm}}
{P}^{k-1}(m)\xrightarrow{\hspace*{0.4cm}}\cdots\xrightarrow{\hspace*{0.4cm}}
{P}^{1}(m)\xrightarrow{\hspace*{0.4cm}}{P}^0(m):=\mathbb{R}^{m+1}
\end{matrix}
\end{gather*}
\end{Rem}

\subsection{Hyperspherical coordinates}\label{hypersphcoord}

Consider the natural global dif\/feomorphism ${\mathcal F}^k: {\mathcal C}^k(m)\rightarrow \mathbb{R}^{m+1}
\times
({\mathbb{S}}^m)^{k}$ given by
\begin{gather*}
{\mathcal F}^k(x_0,x_1,\dots,x_i,\dots,x_{k})= (x_0,x_1-x_0,\dots,x_i-x_{i-1},\dots,x_{k}-x_{k-1}).
\end{gather*}

Now, according to Theorem~\ref{artmod}, the map ${\mathcal F}^k\circ F^k$ is  a~global dif\/feomorphism from
$\hat{\mathcal P}^k(m)$ to ${\mathcal S}^k(m)= \mathbb{R}^{m+1}\times({\mathbb{S}}^m)^{k}$
and, if $\varrho^k:{\mathcal S}^k\rightarrow {\mathcal S}^{k-1}$ is the canonical projection, we
have the following commutative diagram:
\begin{gather}\label{diaghypersphere}
\begin{matrix}
F^k
\hspace{15mm} {\mathcal F}^k
\\
\hfill\hat {P}^k(m) \xrightarrow{\hspace*{0.4cm}} {\mathcal C}^k(m) \xrightarrow{\hspace*{0.4cm}} {\mathcal
S}^k(m) \hfill
\\
\Big\downarrow \hat{\pi}^k \hspace{10mm} \Big\downarrow \rho^{k}
\hspace{13mm} \Big\downarrow \varrho^{k}
\\
\hat{P}^{k-1}(m) \xrightarrow{\hspace*{0.4cm}} \mathcal{C}^{k-1}(m)
\xrightarrow{\hspace*{0.4cm}} {\mathcal S}^{k-1}(m)
\\
F^{k-1}
\hspace{15mm}
{\mathcal F}^{k -1}
\end{matrix}
\end{gather}

This global chart identif\/ies, each point $q=(x_0,x_1,\dots,x_i,\dots,x_{k})\in{\mathcal C}^k(m)$ with
\begin{gather*}
\zeta={\mathcal F}^k(q)=(x_0,z_1,\dots,z_i,\dots,z_{k})\in \mathbb{R}^{m+1}
\times
({\mathbb{S}}^m)^{k}.
\end{gather*}

We will put on each factor ${\mathbb{S}}^m$ the charts given by {\it hyperspherical coordinates} in $\mathbb{R}^{m+1}$
def\/ined as usually by the following relations:
\begin{gather*}
z^1=\rho\phi^1(\theta)= \rho\sin {\theta^1} \cdots\sin{\theta^{m-1}}\sin{\theta^m},
\\
z^2=\rho\phi^2(\theta)=\rho\sin
{\theta^1} \cdots\sin{\theta^{m-1}}\cos{\theta^m},
\\
z^3=\rho\phi^3(\theta)=\rho\sin {\theta^1}
\cdots\sin{\theta^{m-2}}\cos{\theta^{m-1}},
\\
\cdots\cdots\cdots\cdots\cdots\cdots\cdots\cdots\cdots\cdots\cdots\cdots\cdots\cdots
\\
z^{k}=\rho\phi^k(\theta)=\rho\sin{\theta^1} \cos{\theta^{2}},
\\
z^{k+1}=\rho\phi^{k+1}(\theta)=\rho\cos{\theta^1},
\end{gather*}
where $\rho^2=(z^1)^2+\dots+(z^{k+1})^2$, $0\leq \theta^m\leq 2\pi$ and $0\leq\theta^j\leq \pi$ for $1 \leq j \leq m-1$.

We denote by $\hat{\Phi}$ the map from $]0,+\infty[
\times
[0,\pi]
\times
\dots
\times
[0,\pi]
\times
[0,2\pi]$ to $\mathbb{R}^{m+1}$ def\/ined by $\hat{\Phi}(\rho,\theta)=\rho\Phi(\theta)=z$.
The jacobian matrix $D\hat{\Phi}$ of $\hat{\Phi}$ is
\begin{gather*}
D\hat{\Phi}=
\left(
\phi \ \rho\frac{\partial \phi}{\partial \theta^1} \ \cdots \ \rho\frac{\partial \phi}{\partial \theta^m}
\right),
\end{gather*}
where $\phi$ and $ \frac{\partial \phi}{\partial \theta^j}$ are the column vectors of components
$\left(\phi^1,\dots,\phi^{m+1}\right)$ and $\left(\frac{\partial \phi^1}{\partial\theta^j},\dots,\frac{\partial
\phi^{k+1}}{\partial \theta^j}\right)$, respectively.

The inverse of this matrix is then the transpose of the matrix
\begin{gather*}
\left(
\phi  \frac{1}{\rho||\frac{\partial \phi}{\partial \theta^1}||}\frac{\partial \phi}{\partial\theta^1} \  \cdots\
  \frac{1}{\rho||\frac{\partial\phi}{\partial \theta^m}||}\frac{\partial \phi}{\partial \theta^m}
\right).
\end{gather*}

For $i=1,\dots, k$, let ${\mathbb{S}}_{i}$ be the canonical sphere in the $i^{th}$ factor $\mathbb{R}^{m+1}_{i}$.
Given a~point $\alpha$ in the sphere ${\mathbb{S}}_{i}$, there exists hyperspherical coordinates
$z_{i}=\hat{\Phi}_{i}(\rho_{i},\theta_{i})=\rho_{i}\Phi_i(\theta_i^1,\dots\theta_i^k)$ def\/ined for $0\leq \theta_i^k\leq
2\pi$ and $0<\theta_i^j<\pi$, $j=1,\dots,m-1$, where $\Phi_i(0,\dots,0)=\alpha$.
Therefore, given a~point $\zeta=(x_0,z_1,\dots,z_i,\dots,z_{k})\in {\mathcal S}^k$, we obtain a~chart ${\mathcal
H}^k=({\rm Id}-x_0,(\hat{\Phi}_1)^{-1},\dots,(\hat{\Phi}_2)^{-1},\dots,(\hat{\Phi}_{k})^{-1})$ centered at $\zeta$ such that
its restriction to $\rho_i=1$, $i=1,\dots,k$, induces a~chart of ${\mathcal S}^k$ (centered at~$\zeta$).

Note that the map ${\mathcal H}^k=({\rm Id}-x_0,(\hat{\Phi}_1)^{-1},(\hat{\Phi}_2)^{-1},\dots,(\hat{\Phi}_{k})^{-1})$
is a~hyperspherical chart on ${\mathcal S}^k(m)$.
\begin{Def}\quad 
\begin{enumerate}\itemsep=0pt
\item For any $ \zeta\in {\mathcal S}^k(m)$, every map of type ${\mathcal H}^k$ around $\zeta$ is called
a~hyperspherical chart on ${\mathcal S}^k(m)$.
\item For any $q=({\mathcal F}^k)^{-1}(\zeta)$ in ${\mathcal C}^k(m)$, every map of type ${\mathcal H}^k\circ {\mathcal
F}^k$ around~$q$ is called a~hyperspherical chart on ${\mathcal C}^k(m)$.
\item For any $p=({\mathcal F}^k\circ{F^k})^{-1} (\zeta)$ on $\hat{\mathcal P}^k(m)$, every map of type ${\mathcal
H}^k\circ {\mathcal F}^k\circ F^k$ around~$p$ is called a~hyperspherical chart on $\hat{\mathcal P}^k(m)$.
\end{enumerate}
\end{Def}

\begin{Not}\label{nota2}\quad
\begin{itemize}\itemsep=0pt
\item $ A_i= \sum\limits_{r=1}^{m+1}\phi_{i-1}^r\phi_{i}^r$ for $i=1,\dots,k-1$ and $A_{k}=1$,
\item $Z_0= \sum\limits_{r=1}^{m+1}\phi_0^r\frac{\partial}{\partial x^r_0}$,
\item $Z_i= \sum\limits_{j=1}^m B_{i}^j\frac{\partial}{\partial{\theta_{i-1}^j}}$ for $i=1,\dots,k-1$,
with
$B^1_i= \sum\limits_{r=1}^{m+1}\frac{\partial{\phi_{i-1}^r}}{\partial {\theta_{i-1}^1}}\phi_{i}^r$
for $i=1,\dots,k-1$, and
$B_i^j={ \frac{1}{||\frac{\partial\phi_{i-1}}{\partial
\theta_{i-1}^j}||}} \sum\limits_{r=1}^{m+1} \frac{\partial{\phi_{i-1}^r}}{\partial
{\theta_{i-1}^j}}\phi_{i}^r$
for $i=1,\dots,k-1$ and $j=2,\dots,m$,
\item $X^i_l= \frac{\partial}{\partial \theta^i_l}$, for $i=1,\dots,m$ and $0\leq l\leq k-1$,
\item $X^0_l= \sum\limits_{i=0}^{l} f^i_l Z_i$ for $0\leq l\leq k-1$,
with $f^i_l= \prod\limits_{j=i+1}^{l}A_j$, for $i=0,\dots,l-1$, $0\leq l\leq k-1$ and $f^l_l=1$.
\end{itemize}
\end{Not}
\begin{Rem}\label{Aj}
In Lemma~\ref{distribinduite} we already def\/ined a~function ${\mathcal A}_j(q)=-\langle{\mathcal N}_{j}(q),{\mathcal N}_{j-1}(q)\rangle$.
It is clear that we have the relation $A_j\circ ({\mathcal H}^k\circ{\mathcal F}^k)={\mathcal A}_j$.
\end{Rem}

Theorem~\ref{1}(2) is obtained from (2) of the following result:
\begin{The}
\label{distribhyper}
For any $k\geq 1$, we have the following properties:
\begin{enumerate}\itemsep=0pt
\item[$1.$] In hyperspherical coordinates, on each manifold ${\mathcal S}^k(m)$, ${\mathcal C}^k(m)$ and $\hat{\mathcal P}^k(m)$,
the correspon\-ding typical distributions ${\mathcal F}^k_*({\mathcal D}_k)$, ${\mathcal D}_k$ and $\hat{\Delta}_k$
are generated by $\{X^0_{k-1},X^1_{k-1},\dots,X^m_{k-1}\}$, respectively.
\item[$2.$] We have a~net of commutative diagrams:
\begin{gather}\label{folds}
\begin{matrix}
\hfill \hat{P}^k(m) \xrightarrow{\hspace*{0.4cm}}\hat{P}^{k-1}(m) \xrightarrow{\hspace*{0.4cm}}\cdots
\xrightarrow{\hspace*{0.4cm}}\hat{P}^{1}(m) \xrightarrow{\hspace*{0.4cm}}\hat{P}^0(m):=\mathbb{R}^{m+1}\hfill
\\
\vspace*{0.3cm}
\hspace*{-11mm}\Big\downarrow
\hspace{21mm}\Big\downarrow
\hspace{11.5mm} \cdots\hspace{14mm}
\Big\downarrow
\hspace{17mm} \Big\downarrow
\\
\hfill {\mathcal C}^k(m) \xrightarrow{\hspace*{0.4cm}}\mathcal{C}^{k-1}(m)\xrightarrow{\hspace*{0.4cm}}\cdots
\xrightarrow{\hspace*{0.4cm}} \mathcal{C}^1(m) \xrightarrow{\hspace*{0.4cm}} \mathcal{C}^0(m):=\mathbb{R}^{m+1}\hfill
\\
\vspace*{0.3cm}
\hspace*{-11mm}\Big\downarrow
\hspace{21mm}\Big\downarrow
\hspace{11.5mm} \cdots\hspace{14mm}
\Big\downarrow
\hspace{17mm} \Big\downarrow
\\
\hfill {\mathcal S}^k(m) \xrightarrow{\hspace*{0.4cm}}\mathcal{S}^{k-1}(m) \xrightarrow{\hspace*{0.4cm}}\cdots
\xrightarrow{\hspace*{0.4cm}} \mathcal{S}^1(m) \xrightarrow{\hspace*{0.4cm}} \mathcal{S}^0(m):=\mathbb{R}^{m+1}\hfill
\\
\vspace*{0.3cm}
\hspace*{-11mm}\Big\downarrow
\hspace{21mm}\Big\downarrow
\hspace{11.5mm} \cdots\hspace{14mm}
\Big\downarrow
\hspace{17mm} \Big\downarrow
\\
\hfill{P}^k(m) \xrightarrow{\hspace*{0.4cm}} {P}^{k-1}(m) \xrightarrow{\hspace*{0.4cm}}\cdots
\xrightarrow{\hspace*{0.4cm}} {P}^{1}(m) \xrightarrow{\hspace*{0.4cm}}{P}^0(m):=\mathbb{R}^{m+1}
\end{matrix}
\end{gather}
with the following properties:
\begin{itemize}\itemsep=0pt
\item in each horizontal tower, the horizontal map between the space number~$l$ and the space number $l-1$
$(l\geq 1)$ is a~spherical fibration and a~projective space fibration in the first three lines and in the last line, respectively.
\item In each column number~$l$ each vertical map between two consecutive lines among the first three lines is
a~diffeomorphism which sends the typical distribution over the source space on the typical distribution over the image
space, and each vertical map between the two last lines is a~two-fold covering which have the same properties.
\end{itemize}
\end{enumerate}
\end{The}

\begin{proof}
In Section~5 of~\cite{SP1}, it is proved that, in hyperspherical coordinates, the distribution
$F^k_*({\mathcal D}_k)$ is precisely generated by $\{X^0_{k-1},X^1_{k-1},\dots,X^m_{k-1}\}$.
According to Theorem~\ref{artmod}, the dif\/feomorphism $F^k:\hat{\mathcal P}^k(m)\rightarrow {\mathcal C}^k(m)$ sends the
distribution $\hat{\Delta}_k$ onto ${\mathcal D}_k$.
This ends the proof of~(1).
 (2) is a~consequence of relation~\eqref{towerS0}, Theorem~\ref{artmod}, diagram~\eqref{diaghypersphere} and~(1).
\end{proof}
\begin{Rem}
\label{analytique}
All manifolds which appear in the towers~\eqref{folds} are analytic manifolds and all maps in these towers are also
analytic.
\end{Rem}

\section{{\bf RC} and {\bf RVT} codes and conf\/igurations of an articulated arm}

\subsection{{\bf RC} and {\bf RVT} codes according to~\cite{CMH}}\label{code}

In this subsection the theory of {\bf RC} and {\bf RVT} codes introduced in~\cite{CMH} will be adapted to the context
of spherical prolongation.

Consider a~distribution ${D}$ of constant rank on a~manifold~$M$ f\/itted with a~Riemannian metric~$g$.
We will denote, indistinctly, by ${\mathcal P}(M, D)$ the sphere bundle $S(M,D,g)$ or the projective bundle $P(M,D)$ and
by $D^{\{1\}}$ the spherical prolongation or the Cartan prolongation of~$D$ on ${\mathcal P}(M, D)$.
Therefore, we have the associated tower of bundles (see~\eqref{towerM} and~\eqref{towerS})
\begin{gather}\label{towerH}
{\mathcal P}^k(M) \xrightarrow{\hspace*{0.4cm}}{\mathcal P}^{k-1}(M) \xrightarrow{\hspace*{0.4cm}}\cdots
\xrightarrow{\hspace*{0.4cm}}{\mathcal P}^1(M) \xrightarrow{\hspace*{0.4cm}} {\mathcal P}^0(M):=M,
\end{gather}
where each manifold ${\mathcal P}^j$ is equipped with a~distribution denoted by $\mathfrak{D}_j$ such that
$\mathfrak{D}_0=TM$ and $\mathfrak{D}_j =(\mathfrak{D}_{j-1})^{\{1\}}$ for $1\leq j\leq k$.
When $M=\mathbb{R}^{m+1}$, ${\mathcal P}^k(m)$ denotes either $\hat{P}^k(m)$ or ${P}^k(m)$.

For any $1\leq j\leq k$, we denote by $\pi^j$ the natural projection of ${\mathcal P}^{j}(m)$ onto ${\mathcal
P}^{j-1}(m)$.
The tangent bundle to the f\/iber of $\pi^j:{\mathcal P}^j(m)\rightarrow {\mathcal P}^{j-1}(m)$ is the {\it vertical
bundle} denoted by $V_j$, and at any $p\in {\mathcal P}^j(m)$, $V_j(p)\subset\mathfrak{D}_j(p)$ by construction.

For any $p\in {\mathcal P}^{k-1}(m)$, the f\/iber $(\pi^{k})^{-1}(p)$ is denoted by $S^{k}(p)$.
Thus, for such a~point~$p$,~\eqref{towerH}, yields a~tower of f\/iber bundles
\begin{gather}\label{Ftower}
{\mathcal P}^{l}\big(S^k(p)\big) \xrightarrow{\hspace*{0.4cm}}{\mathcal P}^{l-1}\big(S^k(p)\big) \xrightarrow{\hspace*{0.4cm}} \cdots
\xrightarrow{\hspace*{0.4cm}} {\mathcal P}^{1}\big(S^k(p)\big) \xrightarrow{\hspace*{0.4cm}} {\mathcal P}^{0}\big(S^k(p)\big):=S^k(p)
\end{gather}
for any $1\leq l\leq k$.

Coming back to our general context, we have again a~distribution $\mathfrak{d}_j^k$ def\/ined inductively on each
${\mathcal P}^{j}(S^k(p))$, by $\mathfrak{d}_0^k=V_k$ and $\mathfrak{d}_j^k=[\mathfrak{d}_{j-1}^k]^{\{1\}}$ for $1\leq j\leq k$.
Such a~tower will be called {\it a~fiber prolongation tower}.
Of course, we have ${\mathcal P}^{l}(S^k(p))\subset {\mathcal P}^{k+l}(m)$ and~$\mathfrak{d}_j^k(q)$ is an hyperplane in~$\mathfrak{D}_{k+j}(q)$ for any~$q\in {\mathcal P}^j(S^k(p))$.
In particular, $\mathfrak{d}_0^k(q)$ is nothing else but the tangent space of~$S^k(p)$ at~$q$.

On the other hand, for any $k>l\geq 0$, let $\pi^{kl}$ be the natural projection of ${\mathcal P}^{k}(m)$ onto
${\mathcal P}^{l}(m)$ given by the composition $\pi^{k}\circ\pi^{k-1}\circ\dots\circ \pi^{l+1}$.
If~$p$ is a~point in ${\mathcal P}^{k}(M)$, we denote by $p_l$ its projection $p_l=\pi^{kl}(p)\in {\mathcal P}^{l}(m)$,
and we say that $ p_l$ {\it is under $p_k=p$}.
With these notations, for $k\geq 1$, each point $p_k\in{\mathcal P}^{k}(m)$ can be written $(p_{k-1},z)$ for some $z\in
S^{k}({p_{k-1}})$.

It follows that, at each level $k\geq 1$ we have the family of hyperplanes $\mathfrak{d}^i_j(p)$ inside the space~$\mathfrak{D}_k(p)$, with $i+j=k$.
In fact, each $\mathfrak{d}^i_j(p)$ comes from a~f\/iber prolongation of order~$j$ of the tangent space of the f\/iber
$S^i(p_{i-1})$ for $i=1,\dots,k$.

Recall that a~family $\left(E_i\right)_{i\in I}$ of hyperplanes of $\mathbb{R}^N$ is in general position if for every
f\/inite subset~$J$ of $ I$ the codimension of the intersection $\bigcap\limits_{i\in J}E_j$ is exactly equal to the
cardinal of~$J$.

\begin{The}[see~\protect{\cite{CMH}}] The family of hyperplanes $\mathfrak{d}^i_j(p)$ with $i+j=k$ is in general position inside the space $\mathfrak{D}_k(p)$.
\end{The}

\begin{proof}
This result is proved in~\cite{CMH} for Cartan prolongation (Theorem 6.1).
If $\hat{p}$ is a~point in $\hat{P}^k(m)$, we denote by~$p$ its projection $\tau^k(\hat{p})$ in $P^k(m)$ (see
Theorem~\ref{modelS}).
According to Theorem~\ref{modelS}, each hyperplane $\mathfrak{d}^i_j(\hat{p})$ in $\mathfrak{D}_k(\hat{p})$ projects, via
$\tau^k$, onto a~hyperplane $\mathfrak{d}^i_j({p})$ in ${\Delta}_k({p})$ corresponding to the previous process of f\/iber
prolongation in the equation~\eqref{towerS}.
The proof in the context of spherical prolongations is then a~consequence of Theorem~6.1 in~\cite{CMH}.
\end{proof}

According to~\cite{CH, CMH} we have the following def\/initions:

{\samepage \begin{Def}\quad 
\begin{enumerate}\itemsep=0pt
\item Any hyperplane $\mathfrak{d}^i_j(p)$ with $i+j=k$ in the space $\mathfrak{D}_k(p)$ is called a~{\it critical hyperplane} at~$p$.
A direction~$l$ or a~vector~$v$ in $\mathfrak{D}_k(p)$ is called {\it critical} if it lies in at least one critical hyperplane.
Otherwise~$l$ or~$v$ is called {\it regular}.
Moreover, a~critical direction~$l$ or a~vector~$v$ in $\mathfrak{D}_k(p)$ is called {\it vertical} or {\it tangency} if
the singular hyperplane containing this direction is $V_k(p)=\mathfrak{d}_0^k({p})$ or $\mathfrak{d}^i_j(\hat{p})$ for $j>0$, respectively.
\item A point $p=(p_{k-1},z)\in {\mathcal P}^k(m)$ is called {\it regular, critical, vertical or tangency} if $z\in
\mathfrak{D}_{k-1}(p_{k-1})$ is regular, critical, vertical or tangency respectively.
\end{enumerate}
\end{Def}}

\begin{Rem}\label{prolongprolongF}\quad
\begin{enumerate}\itemsep=0pt
\item Let $\hat{p}\in \hat{P}^k(m)$ and $p=\tau^k(\hat{p})\in P^k(m)$.
It follows from Theorem~\ref{modelS} that $\hat{p}$ is regular, critical, vertical or tangency if and only if~$p$ is
respectively regular, critical, vertical or tangency.
Conversely, for any $p\in P^k(m)$, each point in $\tau^k(p)\subset \hat{P}^k(m)$ has the same previous qualif\/ication
as~$p$.
\item We can consider, inside any f\/iber prolongations tower given by equation~\eqref{Ftower},  a~f\/iber prolongation
tower from some f\/iber of the projection ${\mathcal P}^{l}(S^k(p))\rightarrow {\mathcal P}^{l-1}(S^k(p))$ and look for
the corresponding critical hyperplane in $\mathfrak{d}^k_l(q)$.
Then such a~critical hyperplane is in fact an intersection of type $\mathfrak{d}^k_l(q)\cap \mathfrak{d}^i_j(q)$ with
$i>k$ and $k+l=i+j$ (see Proposition~6.2 in~\cite{CMH}).
\item If a~point $p=(p_{k-1},z)\in {\mathcal P}^k(m)$ is critical, then~$z$ may belong to the intersection of several
critical hyperplanes and {\it not only to one} critical hyperplane.
\end{enumerate}
\end{Rem}

The {\bf RC} code of a~point $p\in {\mathcal P}^k(m)$ is a~word $\sigma=\sigma_1\dots\sigma_l\dots\sigma_k$ whose letter
$\sigma_l$ is~$R$ or~$C$ if the point $p_l$ under $ p$ is regular or critical respectively.
Note that, by convention, the f\/irst letter is always~$R$.
Let $\sigma$ be the {\bf RC} code of a~point $p\in {\mathcal P}^k(m)$.
The {\bf RVT} code of~$p$ is a~word $\omega=\omega_1\dots\omega_l\dots\omega_k$ obtained from $\sigma$ in the following way:
\begin{itemize}\itemsep=0pt
\item $\omega_i=R$ if $\sigma_i=R$,
\item $\omega_i=V$ if $\sigma_i=C$ and the point $p_i$ under~$p$ is vertical,
\item $\omega_i=T$ if $\sigma_i=C$ and the point $p_i$ under~$p$ is tangency.
\end{itemize}

\begin{Rem}\label{propcode}\quad
\begin{enumerate}\itemsep=0pt
\item According to Remark~\ref{prolongprolongF}, the {\bf RC} or {\bf RVT} code of any point $\hat{p}\in \hat{P}^k(m)$
is the same as the {\bf RC} or {\bf RVT} code of its projection $p=\tau^k(\hat{p})\in P^k(m)$, respectively.
\item The {\bf RC} code gives rise to a~partition of ${\mathcal P}^k(m)$ into $2^{k-1}$ sets of points which have the
same {\bf RC} code $\sigma$.
Let $\hat{C}_\sigma$ or ${C}_\sigma$ be the set of point $\hat{p}\in \hat{P}^k(m)$ or $p\in P^k(m)$ whose {\bf RC} code
is $\sigma$ respectively.
Then $\tau^k(\hat{C}_\sigma)=C_\sigma$ and $(\tau^k)^{-1}(C_\sigma)=\hat{C}_\sigma$.
\item Remark~\ref{prolongprolongF}(2) implies that if~$p_i$ is tangency, then~$p_i$ must lie in a~f\/iber tower
prolongation for some~$p_j$ under~$p_i$.
Therefore, if~$i$ is the level at which the f\/irst letter~$C$ appears in a~{\bf RC} word, then the associated point $p_i$
must be vertical.
\item Each {\bf RC} code $\sigma$ generates theoretically $2^{n_\sigma}$ {\bf RVT} codes $\omega$ if $n_\sigma$ is the
number of letters~$C$ in $\sigma$.
However, from (3), a~letter~$T$ cannot immediately follow a~letter~$R$ in such a~code because each
tangency point must lies in a~prolongation tower of some point~$p_j$ under~$p_i$.
Consequently, after a~letter~$R$, there must appear at least a~letter~$V$ before any letter~$T$.
\item According to Remark~\ref{prolongprolongF}(3), for any critical point $p=(p_{k-1},z)\in {\mathcal P}^k(m)$,~$z$
may belong to the intersection of several critical hyperplanes.
Therefore, the {\bf RVT} code generated by the {\bf RC} code may not be well def\/ined.
In this case, when~$z$ is not vertical and belongs to only one of them, we need to be more clear in the code about the
possible letters ``T'' that may be written $T_1, T_2,\dots, T_\nu$.
Moreover, much more complicated codif\/ication is needed if~$z$ belongs to the intersection of several critical hyperplanes.
For instance, if this intersection is a~line, we can use a~codif\/ication by letters $L_1,L_2,\dots$ as is proposed
in~\cite{CH, CMH}.
\end{enumerate}
\end{Rem}

In the {\bf RVT} code, according to Remark~\ref{propcode}(5) and~\cite{CH, CMH}, we will use the
{\it following conventions}:
\begin{itemize}\itemsep=0pt
\item If~$z$
belongs to only one critical hyperplane we will use the letters $V, T_1, T_2, \dots$
in the {\bf RVT} code.
\item If~$z$ belongs to the intersection between exactly two critical hyperplanes referenced $T_i$ and $T_j$ we will use
letters of type $T_{ij}$ in the {\bf RVT} code.
Moreover, {\it we adopt the following convention}:
$\bf T_0$ is always relative to a~{\it vertical hyperplane} and $\bf T_i$, for $\bf i>0$, is relative  to a~{\it critical hyperplane
which is not vertical}.
\item More generally, if~$z$ belongs to the intersection of exactly~$n$ critical hyperplanes referenced $T_0=V$ and
$T_{i_1},\dots, T_{i_{n-1}}$ or $T_{i_1},\dots, T_{i_n}$ with $i_i\dots i_n\not=0$ we will use letters of type
$T_{0{i_1}\dots{i_n}}$ or $T_{i_1\dots {i_n}}$ in the {\bf RVT} code.
\end{itemize}

{\it Note that, in a~{\bf RVT} code, a~letter $T_0$ always means ``vertical'' and each letter $T_i$ with $i>0$ means
``tangency''}.
\begin{Def}
We will say that a~word $\omega$ (resp.
a~class ${\mathcal C}_\omega$) in {\bf RVT} code is of depth~$d$ if this word $\omega$ contains at least one letter of
type $T_{i_1\dots i_d}$ with $i_1\geq 0$.
\end{Def}

For instance, for $m=3$ and $1\leq k\leq 4$ the dif\/ferent {\bf RVT} codes which may appear are the following (compare
with~\cite{CMH} before Corollary 4.48):
\begin{alignat*}{3}
&\underline{k=1}: \quad && R,&\\
&\underline{k=2}: \quad && RR , \   RV, & \\
&\underline{k=3}: \quad && RRR, \ RRV, \ RVV, \ RVR, \ RVT, \ RT_0T_{01},& \\
&\underline{k=4}: \quad && RRRR, \ RRRV,  & \\
&&& RRVR , \ RRVV, \ RRVT, \ RRT_0T_{01} & \\
&&& RVRR, \ RVRV, \  RVVR, \ RVVV, \ RVVT, \ RVT_0T_{01}, & \\
&&& RVTR, \ RVTV, \ RVTT, \ RT_0T-1T_{01}, & \\
&&& RT_0T_{01}R, \ RT_0T_{01}V, \ RT_0T_{01}T_1, \ RT_0T_{01}T_2, \ RT_0T_{01}T_{01}, \ RT_0T_{01}T_{02}, \ RT_0T_{01}T_{12}. &
\end{alignat*}
All the words $RT_0T_{01}$, $RRT_0T_{01}$, $RVT_0T_{01}$, $RT_0T_1T_{01}$, $RT_0T_{01}R$, $RT_0T_{01}V$,
$RT_0T_{01}T_1$, $RT_0T_{01}T_2$, $RT_0T_{01}T_{01}$, $RT_0T_{01}T_{02}$ and $RT_0T_{01}T_{12}$ are of depth~$2$.
The other ones are of depth~$1$.

\subsection{Vertical points and conf\/igurations of articulated arms}\label{vertarm}

Consider a~point $q=(x_0,\dots, x_k)\in {\mathcal C}^k(m)$ and let $\hat{p}=(F^k)^{-1}(q)\in \hat{P}^k(m)$.
For $0\leq l\leq k$, we denote by $q_l=F^l(\hat{p}_l)\in{\mathcal C}^k(m)$, where $\hat{p}_l$ is any point under $\hat{p}$.
In fact, according to Theorem~\ref{artmod}, we have $q_l=(x_0,\dots,x_l)=\rho^k\circ\dots\circ
\rho^{l-1}(q)=\rho^{k,l-1}(q)$, where $\rho^k:{\mathcal C}^k(m)\rightarrow {\mathcal C}^{k-1}(m)$ is the natural
projection (see Theorem~\ref{artmod}).
We also say that $q_l$, $l=0,\dots, k-1$ are {\it points under}~$q$.
Moreover, we can write $q=(q_{k-1}, x_k)$ for $x_k\in (\rho^k)^{-1}(q_{k-1})$ and, by Theorem~\ref{artmod},
$\hat{p}=(F^k)^{-1}(q)$ is vertical if and only if the direction generated by $x_k-x_{k-1}$ is vertical according to the
projection $\rho^k:{\mathcal C}^k(m)\rightarrow {\mathcal C}^{k-1}(m)$.

More generally we can transpose the qualif\/ication of points of $\hat{P}^k(m)$ onto points of ${\mathcal C}^k(m)$:
\begin{Def}
Consider a~point $q=(F^k)(\hat{p})\in {\mathcal C}^k(m)$.
\begin{enumerate}\itemsep=0pt
\item~$q$ is called {\it regular, critical, vertical or tangency} if $\hat{p}$ is regular, critical, vertical or
tangency respectively.
\item The code of~$q$ will be the code of the corresponding point $\hat{p}$.
\end{enumerate}
\end{Def}

First of all, we have the following characterization of vertical points in ${\mathcal C}^k(m)$:
\begin{Pro}\label{caractV}
Fix some point $q\in {\mathcal C}^k(m)$.
\begin{enumerate}\itemsep=0pt
\item[$1.$] For all $2\leq l\leq k-1$, we have the following equivalent properties:
\begin{enumerate}\itemsep=0pt
\item[$(a)$] Consider the sandwich of rank~$l$:
\begin{gather*}
\begin{matrix}
[{\mathcal D}_{k}]_l&\subset & [{\mathcal D}_k]_{l-1}
\\
\hfill \cup\hfill&& \hfill\cup\hfill
\\
\hfill L([{\mathcal D}_k]_{l-1})&\subset &  L([{\mathcal D}_k]_{l-2})&
\end{matrix}
\end{gather*}
associated with ${\mathcal D}_k$. Then $[{\mathcal D}_{k}]_l(q)$ is contained in $ L([{\mathcal D}_k]_{l-2})(q)$.
\item[$(b)$] ${\mathcal A}_{l-1}(q)=0$.
\item[$(c)$] The configuration~$q$ of the articulated arm $(M_0,\dots,M_k)$ is such that the segments $[M_{l-2},M_{l-1}]$
and $[M_{l-1},M_{l}]$ are orthogonal at $M_{l-1}$.
\item[$(d)$] $q_l$ is vertical.
\end{enumerate}
\item[$2.$] $q$ is a~Cartan point if and only if each point $q_l$ under~$q$ is regular for $1\leq l\leq k-1$.
\end{enumerate}
\end{Pro}
\begin{Rem}
According to our def\/inition of a~{\it singular} points (see Def\/inition~\ref{Cart}), Proposition~\ref{caractV} implies
that a~point $q\in{\mathcal C}^k(m)$ is singular if and only if there exists a~point $q_l$ under~$q$ which is vertical.
\end{Rem}

A consequence of Proposition~\ref{caractV} is the following:
\begin{The}\label{verticalsets}\quad
\begin{enumerate}\itemsep=0pt
\item[$1.$] For $k\geq 2$, the set ${\mathcal C}_S$ of  singular points of $ {\mathcal C}^k(m)$ is a~subanalytic set of
codimension~$1$.
In particular, the set ${\mathcal C}^k_C(m)={\mathcal C}^k(m)\setminus {\mathcal C}_S$ of Cartan points is an open dense set.
\item[$2.$] Let $\omega$ be a~word of length~$k$ in letters~$R$ and~$V$ and denote by $\{i_1,\dots, i_\nu\}$ the set of index
$\{i\in \{1,\dots, k\}\,|\,  \omega_i=V\}$.
We have the following properties:
\begin{enumerate}
\item[$(i)$] The set ${\mathcal C}_{\omega}$ of points $q\in {\mathcal C}^k(m)$ whose {\bf RVT} code is $\omega$ is an
analytic submanifold of ${\mathcal C}^k(m)$ of codimension $\nu$.
\item[$(ii)$] The configuration of an articulated arm $(M_0,\dots,M_k)$ belongs to ${\mathcal C}_{\omega}$ if and only if
the unique consecutive segments $[M_{i-2},M_{i-1}]$ and $[M_{i-1},M_{i}]$ which are orthogonal at point $M_{i-1}$ occur
for $i=i_1,\dots, i_\nu$.
\end{enumerate}
\end{enumerate}
\end{The}
\begin{Rem}\quad 
\begin{enumerate}\itemsep=0pt
\item From the def\/inition of a~Cartan point (see Def\/inition~\ref{Cart}) it follows that for $k=1$ all points are Cartan
points and the set ${\mathcal C}_S$ is empty in this case.
\item Def\/inition~\ref{Cart} of a~Cartan point is somewhat dif\/ferent from the def\/inition given in~\cite{CH, CMH}.
However Proposition~\ref{caractV}(2) proves the equivalence of these def\/initions.
\item The result of  Theorem~\ref{verticalsets}(1) is well known (see~\cite{CH,CMH,LR,M1,S}).
\item  Theorem~\ref{verticalsets}(2) is also proved in~\cite{S} but with an another notation for this set.
\end{enumerate}
\end{Rem}

The proof of Proposition~\ref{caractV} needs the following lemma:
\begin{Lem}\label{sandwich}
For $2\leq l \leq k$ consider the sandwich of rank~$l$:
\begin{gather*}
\begin{matrix}
[{\mathcal D}_{k}]_l&\subset& [{\mathcal D}_k]_{l-1}
\\
\hfill\cup\hfill && \hfill\cup\hfill
\\
  L([{\mathcal D}_k]_{l-1})&\subset&  L([{\mathcal D}_k]_{l-2}).
\end{matrix}
\end{gather*}
Then $A_{l-1}(q)=0$ if and only if $[{\mathcal D}_{k}]_l(q)\subset L([{\mathcal D}_k]_{l-2})(q)$, and also $q_l$ is
vertical if and only if $A_{l-1}(q)=0$.
\end{Lem}
\begin{proof}
Let us use Notations~\ref{nota2}.

By a~simple computation (see the proof of Proposition~6.1 in~\cite{SP1}), we conclude that the
distribution $[{\mathcal D}_{k}]_l$ of the multi-f\/lag associated with ${\mathcal D}_k$ is generated, in hyperspherical
coordinates, by the union of the sets $\left\{X_{l-1}^0,X_{l-1}^1,\dots, X_{l-1}^m\right\}$
and $\big\{X_i^j\,|\,  j=1,\dots, m,\; l-2\leq i\leq k-1\big\}$,
and $L([{\mathcal D}_k]_{l-2})$ is generated
by $\big\{X_i^j \,|\,  j=1,\dots, m,\; l-2\leq i\leq k-1\big\}$.

We also get
\begin{gather*}
X_{l -1}^0=A_{l-1}X_{l-2}^0+Z_{l-1}.
\end{gather*}

By construction, each $Z_{l-1}$ belongs to $L([{\mathcal D}_k]_{l-2})$, and hence, if $A_{l-1}(q)=0$, it follows that
$[{\mathcal D}_{k}]_l(q)\subset L([{\mathcal D}_k]_{l-2})(q)$.
Since $\left\{X_{l-2}^0, X_{l-2}^1,\dots,X_{l-2}^{m+1}\right \}$ is a~basis of ${\mathcal D}_{l-1}$
at $q_{l-1}$ and $Z_{l-1}$ is a~linear combination of $\left\{X_{l-2}^1,\dots,X_{l-2}^{m+1}\right\}$ we then have
$X_{l-2}^0(q)\not=0$.
Thus, $[{\mathcal D}_{k}]_l(q)\subset L([{\mathcal D}_k]_{l-2})(q)$ if and only if $A_{l-1}(q)=0$.
According to Remark~\ref{Aj}, this ends the proof of the f\/irst equivalence in Lemma~\ref{sandwich}.

Consider now the dif\/feomorphism $\Psi^{l-1}$ from $S(\mathcal{D}_{l-1},\mathcal{C}^{l-1}(m),\gamma_{l-1})$ onto
${\mathcal C}^l(m)$ given in Proposition~\ref{sprolongarm}.
We can write $\Psi^{-1}(q_l)=(q_{l-1}, w_l)$ where $w_l$ is a~vector of norm $1$ in $\mathcal{D}_{l-1}(q_{l-1})$.
The family $\Big\{(x_{l-1}^r-x_{l-2}^r) Y_{l-1}+ \frac{\partial}{\partial x_{l-1}^r} \,\Big|\,
r=1,\dots,m+1\Big\}$ of vector f\/ields (see Lemma~\ref{versionZ}) spans ${\mathcal D}_{l-1}$, and is orthonormal
relative to the metric $\gamma_{l-1}$.
Therefore, we can write
\begin{gather*}
w_l=  \sum\limits_{r=1}^{m+1}z_l^r\left(\big(x_{l-1}^r-x_{l_2}^r\big) Y_{l-1}+ \frac{\partial}{\partial x_{l-1}^r}\right).
\end{gather*}
Moreover, according to this decomposition and from the def\/inition of $\Psi^{l-1}$ in the proof of
Proposition~\ref{sprolongarm}, we have
\begin{gather*}
\Psi^{l-1}(q_{l-1},w_l)=(x_0,\dots,x_{l-1}, x_{l-1}+z_l),
\end{gather*}
where $q_{l-1}=(x_0,\dots,x_{l-1})$ and $z_l=\big(z_l^1,\dots,z_l^{m+1}\big)$.

Since $\Big\{\Pi_{l-1}\big( \frac{\partial}{\partial x_{l-1}^r}\big), \, r=1,\dots,m+1\Big\}$
spans the tangent space to each f\/iber of the projection
${\mathcal C}^{l-1}(m)\rightarrow{\mathcal C}^{l-2}(m)$, the point $q_l$ is vertical if and only if
$ \sum\limits_{r=1}^{m+1}\big(x_l^r-x_{l-1}^r\big)\big(x_{l-1}^r-x_{l-2}^r\big)=0$.

But the f\/irst member of the previous relation is exactly ${\mathcal A}_{l-1}(q)$.
According to Remark~\ref{Aj}, this ends the proof of the Lemma~\ref{sandwich}.
\end{proof}
\begin{proof}[Proof of Proposition~\ref{caractV}]  Proposition~\ref{caractV}(1) is 
a~direct consequence of
Lemma~\ref{sandwich}.

Now,~$q$ is a~Cartan point if and only if $[{\mathcal D}_{k}]_l(q)$ is not contained in $L([{\mathcal D}_k]_{l-2})(q)$
for all $2\leq l\leq k$ (see the end of Section~\ref{multiflag}).
We claim that if~$q$ is a~Cartan point, no point $q_l$ under~$q$ is vertical.
If this was not true, there would exist some point $q_l$ under~$q$ which is tangency.
This would mean that $\hat{p}_l= (F^l)^{-1}(q_l)\in {P}^l(m)$ must be tangency.
Then, from  Remark~\ref{propcode}(3), there must exist a~point $\hat{p}_j$ under $\hat{p}_l$ which is vertical.
Therefore, from  Proposition~\ref{caractV}(1),~$q$ cannot be a~Cartan point.
We conclude that any point $q_l$ under~$q$ is regular.
The converse comes clearly from~(1) of the same proposition.
\end{proof}

\begin{proof}[Proof of Theorem~\ref{verticalsets}] If~$q$ is singular, from  Proposition~\ref{caractV}(2), there must exist
$q_l$ under~$q$ which is vertical.
It follows that the equation of the set ${\mathcal C}_S$ is $ \prod\limits_{l=1}^{k-1}{\mathcal A}_l=0$.

Note that at a~point~$q$ we have
\begin{gather*}
 \frac{\partial{\mathcal A}_l}{\partial x_{l+1}^r}=x_{l}^r-x_{l-1}^r
\qquad
\text{for}
\quad
r=1,\dots,m+1.
\end{gather*}
Taking into account the constraint $||x_l-x_{l-1}||^2=1$, we must have $ \frac{\partial{\mathcal
A}_l}{\partial x_{l+1}^r}(q)\not=0$ for some $1\leq r\leq m+1$.
According to Remark~\ref{analytique}, it follows that ${\mathcal C}_S$ is a~subanalytic subset of ${\mathcal C}^k(m)$ of
codimension $1$, which ends the proof of~(1).

According to the def\/inition of the set $\{i_1,\dots,i_\nu\}$ a~point~$q$ belongs to ${\mathcal C}_{\omega}$ if and only
if each point $q_{i_1},\dots,q_{i_\nu}$ under~$q$ is vertical.
From Proposition~\ref{caractV}(b), the equations of ${\mathcal C}({\omega})$ are then
\begin{gather}\label{equaCs}
{\mathcal A}_i(q)=0
\qquad
\text{for}
\quad
i=i_1-1,\dots, i_\nu-1,
\end{gather}
where each ${\mathcal A}_l$ depends only on the variables $x_{l-1}$, $x_l$ and $x_{l+1}$.
Thus, since $ \frac{\partial{\mathcal A}_l}{\partial x_{l+1}^r}(q)\not=0$ for some $1\leq r\leq m+1$, the
equations in~\eqref{equaCs} are independent.
According to Remark~\ref{analytique}, it follows that ${\mathcal C}_{\omega}$ is an analytic submanifold of ${\mathcal
C}^k(m)$ of codimension~$r$.

Theorem~\ref{verticalsets}(2) is a~direct consequence of Proposition~\ref{caractV}(c).
\end{proof}

\subsection{Tangency points and conf\/igurations of articulated arms}\label{tangearm}

We will prove the fundamental following results for tangency points $q\in {\mathcal C}^k(m)$
\begin{The}\label{tengency}\quad
\begin{enumerate}\itemsep=0pt
\item[$1.$] Assume that $q\in{\mathcal C}^k(m)$ is a~tangency point.
Then there exists $2\leq i\leq k-1$ such that the point $q_i$ under~$q$ is vertical.
We define
\begin{gather*}
l=\sup\{2\leq i\leq k-1
\
\text{such that}
\
q_i
\
\text{is vertical}\}.
\end{gather*}
Then, if $l<k$, for any $l<j\leq k$, the point $q_j$ under~$q$ must be tangency.
\item[$2.$] Denote by $R^hVT^l$ a~word of length $h+l+1\leq k$ in letters $R$, $T$, $V$, where $R^h$ denotes~$h$ consecutive letters~$R$ and $T^l$ denotes~$l$ consecutive letters~$T$.
Then the set ${\mathcal C}_{R^hVT^l}$ of points $q\in{\mathcal C}^{h+l+1}(m)$ whose {\bf RVT} code is ${R^hVT^l}$ is an
analytic submanifold of ${\mathcal C}^{h+l+1}(m)$ of codimension $l+1$.
The fiber of the projection of ${\mathcal C}_{R^hVT^l}$ onto $({\mathcal C}^h(m))_C$ over $q_h\in ({\mathcal C}^h(m))_C$
is the set $F^{h+l+1}(\hat{P}^{l}(\hat{S}^{h+1}(q_h)))$.
\item[$3.$] To each $q=(x_0,\dots, x_k)\in {\mathcal C}^k(m)$ and $0\leq h<k$ we associate a~field of directions $K_h(q)$ on
$\mathbb{R}^{m+1}$ defined by $K_h(q)$ generated by $x_{h+1}-x_h$.
Given a~configuration $q\in {\mathcal C}^k(m)$ of an articulated arm $ (M_0,\dots,M_k)$, the configuration $q_{h+l+1}$
of the induced articulated arm $(M_0,\dots,M_{h+l+1})$ belongs to ${\mathcal C}_{R^hVT^l}$ if and only if, each segment
$[M_{h+i},M_{h+i+1}]$ is orthogonal to the direction $K_h(q)$ for all $i=0,\dots,l$ and no other orthogonality constraint.
\end{enumerate}
\end{The}

For the proof 
we need  some intrinsic characterization of $F^{h+l+1}(\hat{P}^{l}(\hat{S}^{h+1}(q_h)))$
in ${\mathcal C}^{h+l+1}(m)$ and the critical hyperplane
$F^{h+i+1}(\mathfrak{d}^{h+1}_i)$ in ${\mathcal D}_{h+l+1}$.

Given any vertical point $q=(q_{h},w)\in {\mathcal C}^{h+1}(m)$, for $1\leq h< k$, denote by ${\mathcal
C}_q^{h+i+1}(m-1)$ the set $F^{h+i}(\hat{P}^{i}(\hat{S}^{h+1}(q_{h})))$
of $ {\mathcal C}^{h+i+1}(m)$ for
$i=1,\dots,k-h-1$, and call it {\it a~critical manifold}.
For $i=0$ we set ${\mathcal C}_q^{h+1}(m-1)=\hat{S}^{h+1}(q_{h})$.
Let $\delta^{h+1}_i$ be the singular hyperplane $F^{h+i+1}(\mathfrak{d}^{h+1}_i)\subset {\mathcal D}_{h+i+1}$ on
${\mathcal C}_q^{h+l}(m-1)$, for $i=0,\dots,k-h-1$.
We set ${\mathcal V}^{h}_{h+i+1}= \sum\limits_{s=1}^{m+1}(x_{h+1}^s-x_{h}^s) \frac{\partial}{\partial
x_{h+i+1}^s}$ for $i=0,\dots,k-h-1$.

\begin{Pro}
\label{delta}
Consider a~point $q\in{\mathcal C}^k(m)$ such that the point $q_{h+1}=(x_0,\dots,x_h, x_{h+1})$ under~$q$ is vertical.
For $i=0,\dots,k-h-1$,  the following properties hold:
\begin{enumerate}\itemsep=0pt
\item[$(i)$] the manifold ${\mathcal C}_{q_h}^{h+i+1}(m-1)$ is the subset of point $(q_{h},x_{h+1},\dots,x_{h+i+1})\in
{\mathcal C}^{h+i+1}(m)$ such that:
\begin{gather*}
||x_{h+1}-x_{h}||=1
\quad
\text{and}
\quad
\langle x_{j+1}-x_j,x_{h+1}-x_{h}\rangle=0
\quad
\text{for all}
\
j=h+2,\dots,h+i-1,
\end{gather*}
\item[$(ii)$] if $i>0$ the vertical space associated with the fibration ${\mathcal C}_{q_h}^{h+i+1}(m-1)\rightarrow
{\mathcal C}_{q_h}^{h+i}(m-1)$ is generated by
\begin{gather}\label{Tvert}
\left\{ \frac{\partial}{\partial x_{h+i+1}^r}-(x_{h+1}-x_{h}){\mathcal V}_{h+i+1}^{h}-(x_{h+i+1}^r-x_{h+i}^r){\mathcal V}_{h+i+1}
\,\Big|\,
r=1,\dots,m+1\right\},
\end{gather}
\item[$(iii)$] the vector field ${X}_{h+i+1}$ is tangent to ${\mathcal C}_{q_h}^{h+i+1}(m-1)$,
\item[$(iv)$] the distribution $\delta^{h+1}_i$ is the intersection between ${\mathcal D}_{h+i+1}$
and the tangent space to \linebreak ${\mathcal C}_{q_h}^{h+i+1}(m-1)$.
\end{enumerate}
\end{Pro}
\begin{proof}
For $i=0$, the set ${\mathcal C}_{q_h}^{h+1}(m-1)$ is the sphere $\hat{S}^{h+1}(q_{h})$ and the distribution
$\delta_0^{h+1}$ is the tangent space to this sphere.
Moreover, since $q_{h+1}$ is vertical, then ${\mathcal A}_h(q_{h+1})=0$.
Therefore, we have $X_{h+1}={\mathcal Z}_{h}+{\mathcal V}_{h+1}$ which is tangent to $\hat{S}^{h+1}(q_{h})$ in
$q_{h+1}$.
Thus the proper\-ties~(i), (iii) and (iv) of Proposition~\ref{delta} are true for $i=0$.
Assume now that for all $0\leq j<i$ these last properties are true.
According to our assumption, $\delta^{h+1}_{i-1}$ is generated by the family:
\begin{gather*}
\left\{X_{h+1},  \frac{\partial}{\partial x_{h+1}^r}-\big(x_{h+1}^r-x_{h}\big){\mathcal V}_{h+1}\,\Big|\,
r=1,\dots,m+1\right\}
\quad
\text{for}
\quad
i=1,
\\
\left\{X_{h+i},  \frac{\partial}{\partial x_{h+i}^r}-(x_{h+1}-x_{h}){\mathcal V}_{h+i}^{h}-\big(x_{h+i}^r-x_{h+i-1}\big){\mathcal V}_{h+i}
\,\Big|\,
r=1,\dots,m+1\right\}
\quad
\text{for}
\ \
i>1.
\end{gather*}

Recall that ${\mathcal D}_{h+i}$ is generated by:
\begin{gather*}
\left\{X_{h+i},  \frac{\partial}{\partial x_{h+i}^r}-\big(x_{h+i}^r-x_{h+i-1}\big){\mathcal V}_{h+i}\,\Big|\,
r=1,\dots,m+1\right\}.
\end{gather*}

Therefore on ${\mathcal C}_{q_h}^{h+i}(m-1)$ the distribution $\delta^{h+1}_{i-1}$ is the intersection
between ${\mathcal D}_{h+i}$ and the kernel
of the dif\/ferential form $ \sum\limits_{s=1}^{m+1}(x_{h+1}^s-x_{h}^s)(dx_{h+i}^s-dx^s_{h+i-1})$.

According to our assumption, we can see that $\delta^{h+1}_{i-1}$ is generated by the family
\begin{gather*}
\left\{U_r=(x_{h+i}^r-x_{h+i-1}^r)Y_{h+i}+ \frac{\partial}{\partial x_{h+i}^r}-\big(x_{h+i}^r-x_{h+i-1}^r\big){\mathcal V}_{h+i}^h
\,\Big|\,
r=1,\dots,m+1\right\}.
\end{gather*}

According to the proof of Proposition~\ref{sprolongarm}, a~point
$\Big(q_h,x_{h+1},\dots,x_{h+i}, \sum\limits_{r=1}^{m+1}\nu^r{W}_r\Big)$ belongs to the manifold
$\delta_{i-1}^{h+1}$ if and only if $(q_h,x_{h+1},\dots,x_{h+i})$ belongs to ${\mathcal C}_{q_h}^{h+i}(m-1)$ and
$ \sum\limits_{r=1}^{m+1}\nu^r(x_{h+1}^r-x_{h}^r)=0$.

Note that at such a~point we have:
\begin{gather}\label{deltai-1}
 \sum\limits_{r=1}^{m+1}\nu^rU_r
= \sum\limits_{r=1}^{m+1}\nu^r{W}_r- \sum\limits_{r=1}^{m+1}\nu^r\big(x_{h+1}^r-x_h^r\big){\mathcal
V}_{h+1}^h= \sum\limits_{r=1}^{m+1}\nu^r{W}_r.
\end{gather}

Consider the submanifold $S (\delta_{i-1}^{h+1},{\mathcal C}_{q_h}^{h+i}(m-1),\gamma_{h+i})$ of $S({\mathcal
D}_{h+i},{\mathcal C}^{h+i}(m),\gamma_{h+i})$.
According to the proof of Proposition~\ref{sprolongarm} $\Psi^{h+i}(S (\delta_{i-1}^{h+1},{\mathcal
C}_{q_h}^{h+i},\gamma_{h+i}))$ is a~submanifold of ${\mathcal C}^{h+i+1}(m)$ which is f\/ibered on ${\mathcal
C}_{q_h}^{h+i}(m-1)$.
Moreover the restriction of $\Psi^{h+i}$ to $S (\delta_{i-1}^{h+1},{\mathcal C}_{q_h}^{h+i}(m-1),\gamma_{h+i})$ commutes
with the canonical projections of $S \delta_{i-1}^{h+1},{\mathcal C}_{q_h}^{h+i}(m-1),\gamma_{h+i})$ and $\Psi^{h+i}(S
(\delta_{i-1}^{h+1},{\mathcal C}_{q_h}^{h+i}(m-1),\gamma_{h+i}))$ onto ${\mathcal C}_{q_h}^{h+i}(m-1)$ respectively.

Now according to~\eqref{deltai-1}, the manifold $S (\delta_{i-1}^{h+1},{\mathcal C}_{q_h}^{h+i}(m-1),\gamma_{h+i})$ is
the set of points
\begin{gather*}
\left(q_h,x_{h+1},\dots,x_{h+i}, \sum\limits_{r=1}^{m=1}\nu^rW_r\right)\in S\big({\mathcal D}_{h+i},{\mathcal
C}^{h+i}(m),\gamma_{h+i}\big)
\end{gather*}
with the following constraints:
\begin{gather*}
(q_h,x_{h+1},\dots,x_{h+i})\in {\mathcal C}_{q_h}^{h+i}(m-1),
\qquad
 \sum\limits_{r=1}^{m+1}\big(\nu^r\big)^2=1,
\qquad
 \sum\limits_{r=1}^{m+1}\nu^r\big(x_{h+1}^r-x_{h}^r\big)=0.
\end{gather*}

Since ${\mathcal C}^{h+i+1}(m)={\mathcal C}^{h+i}(m)\times\mathbb{S}^m$, then
  $\Psi^{h+i}(S (\delta_{i-1}^{h+1},{\mathcal C}_{q_h}^{h+i}(m-1),\gamma_{h+i}))$ is
a~submanifold of ${\mathcal C}_{q_h}^{h+i}(m-1)\times\mathbb{S}^m$ def\/ined by the equation
\begin{gather}\label{eqCqhi1}
 \sum\limits_{r=1}^{m+1}\big(x^r_{h+i+1}-x^r_{h+i}\big)\big(x_{h+1}^r-x_{h}^r\big)=0.
\end{gather}

But from the construction of $F^{h+i+1}$, $F^{h+i+1}(\hat{P}^{i}(\hat{S}^{h+1}(q_{h}))$ is equal to
$\Psi^{h+i}(S (\delta_{i-1}^{h+1},{\mathcal C}_{q_h}^{h+i}(m-1),\gamma_{h+i}))$, which is precisely
the set ${\mathcal C}^{h+i+1}_{q_h}(m-1)$.
Therefore is proved~(i).
This implies that the vertical bundle associated with the f\/ibration
${\mathcal C}_{q_h}^{h+i+1}(m-1)\rightarrow {\mathcal C}_{q_h}^{h+i}(m-1)$ is generated by the family~\eqref{Tvert}.

Let $\overline{U}_r$ be the vector f\/ield on ${\mathcal C}^{h+i}(m)\times{\mathbb{R}}^{m+1}$ def\/ined by
\begin{gather*}
\overline{U}_r=\big(x_{h+i}^r-x_{h+i-1}^r\big)Y_{h+i}+ \frac{\partial}{\partial
x_{h+i}^r}-\big(x_{h+i}^r-x_{h+i-1}^r\big){\mathcal V}_{h+i}^h + \frac{\partial}{\partial x_{k+i+1}^r}.
\end{gather*}

The vector f\/ield $\overline{U}_r(x_0,\dots,x_{h+i},x_{h+i+1})$ is actually tangent to ${\mathcal C}^{h+i+1}_{q_h}(m)$
and $\overline{U}_r$ projects onto ${U}_r$.

Therefore as in the proof of Proposition~\ref{sprolongarm}, we set ${\mathcal U}_r =
(T\Psi^{h+i})^{-1}(\overline{U}_r)$.
The distribution $[\delta_{i-1}^{h+1}]^{[1]}$ is generated by the vertical bundle of the f\/ibration $S
(\delta_{i-1}^{h+1},{\mathcal C}_{q_h}^{h+i}(m-1),\gamma_{h+i})\rightarrow {\mathcal C}_{q_h}^{h+i}(m-1)$ and the vector
f\/ield $ \sum\limits_{r=1}^{m+1}\nu^r{\mathcal U}_r $.
According to equation~\eqref{eqCqhi1}, the vector f\/ield $T\Psi^{h+i}( \sum\limits_{r=1}^{m+1}\nu^r{\mathcal U}_r)$
is nothing else but $X_{h+i+1}$ on ${\mathcal C}^{h+i+1}_{q_h}(m)$.
Moreover, since $\bar{U}_r$ is tangent to ${\mathcal C}^{h+i+1}_{q_h}(m)$, then $X_{h+i+1}$ is also tangent to
${\mathcal C}^{h+i+1}_{q_h}(m)$.
Finally the distribution $\Psi^{h+i}_*([\delta_{i-1}^{h+1}]^{[1]})$ is generated by $X_{h+i+1}$ and the vertical bundle
of the f\/ibration ${\mathcal C}_{q_h}^{h+i+1}(m-1)\rightarrow {\mathcal C}_{q_h}^{h+i}(m-1)$.
Since this vertical bundle is the intersection between the vertical bundle of the f\/ibration ${\mathcal
C}^{k+i+1}(m)\rightarrow {\mathcal C}^{h+i}(m)$ and the tangent space to ${\mathcal C}_{q_h}^{h+i+1}(m-1)$ then $\Psi^{h+i}_*([\delta_{i-1}^{h+1}]^{[1]})$ is the intersection between ${\mathcal D}_{h+i+1}$ and the tangent space to
${\mathcal C}_{q_h}^{h+i+1}(m-1)$.
But according to the def\/inition of the family of distributions $\{\delta_i ^{h+1}\}_{i\geq 0}$ and the construction of
$F^{h+i+1}$ we get $\Psi^{h+i}_*([\delta_{i-1}^{h+1}]^{[1]})=\delta_i^{h+1}$.
This ends the proof.
\end{proof}

\begin{proof}[Proof of Theorem~\ref{tengency}] Fix some tangency point $q\in {\mathcal C}^k(m)$.
According to Remark~\ref{propcode}(3), there must exist a~vertical point~$q_i$ under~$q$.
Let $q_l$ be the last vertical point under~$q$ then ${\mathcal A}_{l-1}(q)=0$.
In fact, given $q_{l-1}=(x_0,\dots, x_{l-1})$, this relation characterizes the vertical points $q_l=(x_0,\dots,x_{l-1},x_l)\in
{\mathcal C}^l(m)$.
$F^l(\hat{S}^l(q_{l-1}))$ is exactly the set of points $(q_{l-1},x_l)\in {\mathcal C}^l(m)$ such
that ${\mathcal A}_{l-1}(q_{l-1},x_l)=0$.

Now, if $l+1\leq k$, then $q_{l+1}$ cannot be vertical and the point $q_{l+1}$ is no more
regular, since otherwise the existence of a~vertical point between $q_{l+1}$ and~$q$ would contradict the def\/inition of~$l$.
It follows that $q_{l+1}$ must be a~tangency point.

We proceed by induction.
Assume that for $l\leq i< k$ the point $q_i$ is tangency.
By the same arguments as previously, $q_{i+1}$ must also be a~tangency point.
It follows that, by induction, we obtain~(1).

We shall now prove~(2).
For each $0\leq i\leq j\leq k-1$ we set
\begin{gather*}
{\mathcal A}_{j,i}(x_0,\dots,x_k)=\langle x_{j+1}-x_j,x_{i+1}-x_i\rangle.
\end{gather*}

Note that ${\mathcal A}_{j,i}$ is def\/ined on any ${\mathcal C}^{h+l+1}(m)$ for $0\leq i\leq j\leq h+l$.

For $l=0$, the point $q_{h+1}$ is vertical if and only if ${\mathcal A}_{{h+1},h}(q_{h+1})=0$, i.e., if and only if
$q_{h+1}$ belongs to ${\mathcal C}_{q_h}^{h+1}(m-1)$ from Proposition~\ref{delta}(1).
Therefore~(2) is true for $l=0$.
Assume that for all $0\leq i<l$ the point $q_{h+i+1}$ belongs to ${\mathcal C}_{R^hVT^i}$ if and only if ${\mathcal
A}_{j,h}(q_{h+j+1})=0$ for all $j=h+1,\dots,h+i$.
By def\/inition, the point $q_{h+l+1}$ belongs to ${\mathcal C}_{R^hVT^l}$ if and only if $q_{h+l+1}$ is tangency.
From~(1) each point $q_{h+i+1}$ under $q_{h+l+1}$ must be tangency for $i=0,\dots,l-1$.
In particular this means that $q_{h+l}$ belongs to ${\mathcal C}_{q_h}^{h+l}(m-1)$.
It follows that $q_{h+l+1}$ is tangency if and only if the direction generated by $x_{h+l+1}-x_{h+l}$ belongs to
$\delta_{l}^{h}(q_{h+l})$ and according to the proof of Proposition~\ref{delta} we get the relation
$\langle x_{h+l+1}-x_{h+l},x_{h+1}-x_{h}\rangle=0$.
Therefore if $q_{h+l}$ belongs to ${\mathcal C}_{R^hVT^{l-1}}$ then $q_{h+l+1}$ belongs to ${\mathcal C}_{R^hVT^{l}}$ if and only if
 ${\mathcal A}_{h+l,h}(q_{h+l+1})=0$.

${\mathcal C}_{R^hVT^{l}}$ then is def\/ined by the equations
\begin{gather}\label{eqRVT}
{\mathcal A}_{j,h}=0,\qquad j=h+1,\dots,h+l,
\end{gather}
in ${\mathcal C}^{h+l+1}(m)$.
It is clear that these equations are independent.
In particular, $q_{h+i+1}$ belongs to ${\mathcal
C}_{q_h}^{h+i+1}(m-1)=F^{h+i+1}(\hat{P}^i(\hat{S}^{h+1}(q_h))$ for all $i=0,\dots, k-h-l-1$.
This ends the proof of~(2).

 (3) is an interpretation of the equations~\eqref{eqRVT} in terms of orthogonality.
\end{proof}

\section{Relation between EKR classes of depth at most~1,\\ {\bf RVT} codes and articulated arms}\label{relatiomEKRRVT}

\subsection{Mormul EKR coding according to~\cite{M1,M2}}\label{EKR}

In~\cite{M1,M2}, P.~Mormul has constructed a~coding system for labeling singularity classes of germs of special
multi-f\/lag which he called ``extended Kumpera--Ruiz'' (``EKR'' in short).
Mormul's codes are f\/inite sequences in $\mathbb{N}$.
We now summarize how P.~Mormul def\/ines this coding system.

Given a~coordinate system $(y^1,\dots, y^s)$ on $\mathbb{R}^s$ consider a~distribution $\mathcal D$ def\/ined on
a~neighbourhood of $0\in \mathbb{R}^s$ by $m+1$ vector f\/ields $Z_1,\dots,Z_{m+1}$.
A~new distribution~${\mathcal D}'$ is associated with~${\mathcal D}$ on a~neighbourhood of $0\in \mathbb{R}^{s+m}$
relatively to a~coordinate system $(y^1,\dots,y^s,x^1,\dots,x^m)$ by an operation denoted ${\bf j}$ where~$\bf j$ takes
values in $ \{1,2,\dots,m+1\}$ in the following way:
for a~f\/ixed value~$j$ of $\bf j$, the distribution ${\mathcal D}'$ is generated by{\samepage
\begin{gather*}
\begin{split}
& \Bigg\{Z^{'}_{1}=x^{1}Z_{1}+\dots+ x^{j-1}Z_{j-1}+Z_{j}+\big(x^j+c^j\big)Z_{j+1}+\dots+\big(x^m+c^m\big)Z_{m+1},
\\
& \phantom{\Bigg\{ }
Z'_2= \frac{\partial}{\partial x^1},
\
\dots,
\
Z'_{m+1}= \frac{\partial}{\partial x^m}\Bigg\},
\end{split}
\end{gather*}
where $c^j,\dots,c^m$ are constants that may or may not be equal to zero.}

For instance, when $m=2$, ${\mathcal D}'$ is a~distribution of rank $3$ generated by
\begin{gather*}
Z'_1=
\begin{cases}
Z_1+\big(x^1+c^1\big)Z_2+\big(x^2+c^2\big) Z_3& \text{for}\quad {\bf j}=1,
\\
x^1 Z_1+Z_2+\big(x^2+c^2\big)Z_3            &\text{for}\quad {\bf j}=2,
\\
x^1 Z_1+x^2 Z_2+Z_3                 &\text{for}\quad {\bf j}=3,
\end{cases}
\end{gather*}
and $ Z'_2= \frac{\partial}{\partial x^1}$, $Z'_{3}= \frac{\partial}{\partial x^2}$.

This procedure is initiated for $\mathcal D^{(0)}$ generated by
$\left\{Z_0^{(0)}= \frac{\partial}{\partial y^1}, \dots,
Z^{(0)}_{m+1}= \frac{\partial}{\partial y^{m+1}}\right\}$
on $\mathbb{R}^{m+1}$ and we obtain by
a~f\/irst operation $\mathbf{j}_1$ a~new distribution ${\mathcal D}^{(1)}$ of rank $m+1$ on
a~neighborhood of $0\in \mathbb{R}^{2m+1}$ generated by the produced vector f\/ields $\big\{Z^{(1)}_1,\dots,Z^{(1)}_{m+1}\big\}$.
By induction on the composition of consecutive operations ${\bf j_1, j_2,\dots, j_{k}}$
for each word ${\bf j_1 j_2\dots j_{k}}$, we obtain a~distribution $\Delta_{\bf j_1\dots j_{k}}$ on a~neighborhood of $0\in \mathbb{R}^{(k+1)m+1}$
generated by the associated produced $(m+1)$ vector f\/ields $\big\{Z^{(k)}_1,Z_2^{(k)},\dots,Z_{m+1}^{(k)}\big\}$.

We have then the following result:
\begin{The}[see \protect{\cite{M1, M2}}]
Consider a~differential system ${D}$ which spans a~special multi-flag on a~manifold~$M$ of dimension ${(k+1)m+1}$.
Every point $x\in M$ the differential system $(M,D,x)$ is locally equivalent to some differential system $(\Delta_{\bf
j_{1}j_{2} \dots j_{k}}, \mathbb{R}^{(k+1)m+1},0)$.
Moreover, the value of ${\bf j_{1}j_{2} \dots j_{k}}$ can be chosen such that ${\bf j_1}=1$ and, in the case where ${\bf
j_{l+1}}$ $>$ $\max ({\bf j_{1},j_{2}, \dots, j_{l}})$ then we have ${\bf j_{l+1}}$ = $ 1+ \sup ({\bf j_{1},j_{2}, \dots,
j_{l}})$ for $l=1,\dots,k-1$ $($the rule of the least possible new jumps upwards$)$.
\end{The}

Therefore, with a~given germ of distribution $(M,D,x)$ on a~manifold~$M$ it is associated a~well def\/ined sequence of values
${j_{1}j_{2} \dots j_{k}}$ which satisf\/ies the rule of least upward jumps.
Conversely, a~germ of distribution~$D$ determines an unique sequence of operations ${\bf j_1, j_2,\dots, j_{k}}$
satisfying the rule of least upward jumps (see footnote~6 of~\cite{M2}).
This sequence is called a~{\it singularity class of special multi-flags} in~\cite{M1}.
We will say that this a~sequence is an~{\bf EKR} {\it class} of germ of distributions which is encoded by the unique associated
sequence of integers ${j_{1}j_{2} \dots j_{k}}$.
\begin{figure}[t]\centering
\includegraphics{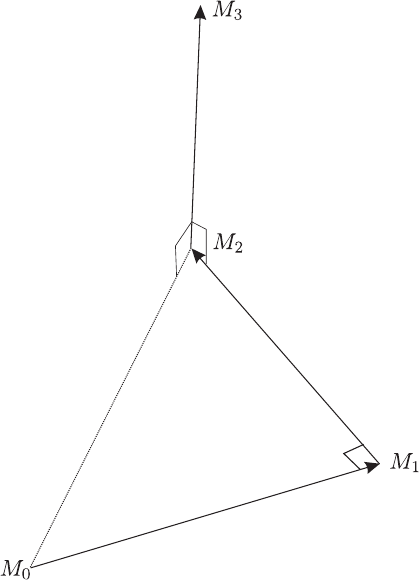}

\caption{EKR set $\Sigma_{123}$ of depth $2$ of an articulated arm $(M_0,M_1,M_2,M_3)$.}\label{sigma123}
\end{figure}

Since ${\mathcal D}_k$ generates a~special multi-f\/lag of step $m \geq 2$ and length $k \geq 1$ on ${\mathcal C}^k(m)$,
it is associated with  any point $q\in {\mathcal C}^k(m)$ a~word ${j_{1}j_{2} \dots j_{k}}$ def\/ined by the germ
$({\mathcal D}_k,{\mathcal C}^k(m), q)$.
We will denote by $\Sigma_{j_{1}j_{2} \dots j_{k}}$ the set of conf\/igurations $q\in {\mathcal C}^k(m)$ such that this
associated word is ${j_{1}j_{2} \dots j_{k}}$ and called it an {\bf EKR} set.
The integer $d=\sup\{j_1,\dots, j_k\}-1$ will be called the {\it depth} of $\Sigma_{j_{1}j_{2} \dots j_{k}}$ (see
Fig.~\ref{sigma123}).

\subsection{Stratif\/ication of EKR sets of depth at most~1 by {\bf RVT} codes}\label{stratEKR1}

According to the notations introduced in Section~\ref{intro}, the following result gives a~complete description of EKR
sets of depth at most $1$ in terms of {\bf RVT} sets.
This result gives a~proof of Theorems~\ref{2} and~\ref{3}(2):
\begin{The}\label{stratEKR}\quad
\begin{enumerate}\itemsep=0pt
\item[$1.$] The EKR set $\Sigma_{1\dots 1}$ is the set of Cartan points.
In particular, $\Sigma_{1\dots 1}$ is an open dense set whose complementary is a~subanalytic set of ${\mathcal C}^k(m)$
of codimension $1$.
\item[$2.$] Let $\Sigma_{j_1\dots j_k}$ be an~EKR set of depth $1$ and $\{i_1,\dots,i_\nu\}$ be the set of indices~$i$ such that $j_i=2$.
Then $\Sigma_{j_1\dots j_k}$ is an analytic submanifold of ${\mathcal C}^k(m)$ of codimension $\nu$.
Moreover,~$q$ belongs to $\Sigma_{j_1\dots j_k}$ if and only if the configuration at~$q$ of the articulated arm
$(M_0,\dots, M_k)$, the segments $[M_{i-2},M_{i-1}]$ and $[M_{i-1},M_{i}]$ are orthogonal at $M_{i-1}$ for
all index $i\in \{i_1,\dots,i_\nu\}$.
\item[$3.$] In the previous situation we have
\begin{enumerate}\itemsep=0pt
\item[$(i$)] A~point~$q$ belongs to $\Sigma_{j_1\dots j_k}$ if and only if, in its ${\bf RVT}$ code the
only letters~$V$ are at rank $i_1, \dots, i_\nu$.
\item[$(ii)$] A {\bf RVT} set ${\mathcal C}_\omega$ is contained in $\Sigma_{j_1{\dots} j_k}$ if and only if $\omega$ is of
type $R^{h_0}VT^{l_1}\!R^{h_1}{\dots}V T^{l_\nu}\!R^{h_{\nu}}$ and each letter~$V$ is exactly at rank $i_1,\dots, i_\nu$.
More precisely, $l_{\lambda}+h_{\lambda}=i_{\lambda+1}-i_{\lambda}-1$ for $\lambda=1,\dots,\nu-1$ and
$l_\nu+h_\nu= k-i_\nu-1$.
This set is an analytic submanifold of $\Sigma_{j_1\dots j_k}$ of codimension $l_1+\dots+l_\nu$.
In particular ${\mathcal C}_{R^{h_0}VR^{h_1}\dots V R^{h_{\nu}}}$ is an open dense set of $\Sigma_{j_1\dots j_k}$ for
$h_0=i_1-1$, $h_{\lambda}=i_{\lambda+1}-i_{\lambda}-1$ for $\lambda=1,\dots,\nu-1$ and $h_\nu=k-i_\nu-1$.
\item[$(iii)$] $\Sigma_{j_1\dots j_k}$ is the union of all sets of type ${\mathcal C}_{R^{h_0}VT^{l_1}R^{h_1}\dots V
T^{l_\nu}R^{h_{\nu}}}$ which satisfies~$(ii)$.
\end{enumerate}
\end{enumerate}
\end{The}
\begin{Rem}
The decomposition of $\Sigma_{j_1\dots j_k}$ given in (iii) above into {\bf RVT} sets is in agreement with the
decomposition of such EKR sets for $k =3$ described by Howard in the appendix of~\cite{CMH}.
Therefore, the description in (iii) can be seen as a~generalization of Howard's result.
\end{Rem}

For the proof of this theorem we need the following proposition which will be used in Section~\ref{2k4}.
\begin{Pro}\label{cactEKR3}
Consider a~configuration~$q$ which belongs to some EKR set $\Sigma_{{j_1}\dots{j_k}}$ of depth~$d$.
\begin{enumerate}\itemsep=0pt
\item[$1.$] If $d\geq 2$, we have $j_h\geq 2$ if and only if ${\mathcal A}_{h-1}=0$.
\item[$2.$] If $d\geq 2$ there exists a~rank~$h$ such that $j_h=3$ and the letter of rank~$h$ in the {\bf RVT} code of~$q$ is
of type $T_{0i}$ with $i\geq 1$.
\item[$3.$] Assume that the {\bf RVT} code of~$q$ is a~word which contains a~letter of type $T_{rs}$.
Then this word contains also a~letter of type $T_{0i}$ for some rank $l<h$ with $i\geq 1$.
If the letter of rank~$h$ is $T_{0i}$ with $i\geq 1$, then $2\leq j_h\leq 3$.
If $j_h=2$ there exists a~letter of type $T_{01}$ at some rank $l<h$ with $j_l=3$.
In particular $d\geq 2$.
\item[$4.$] If $d=2$, assume that the letter of rank~$h$ is $T_{0i}$ and~$h$ is the first index~$h$ such that $j_h=3$.
Then $i=1$ and we have one and only one of the following situations:
\begin{itemize}\itemsep=0pt
\item  $j_{h+1}=2$ if and only if the letter of rank $h+1$ is~$V$ or $T_{02}$,
\item  $j_{h+1}=3$ if and only if the letter of rank $h+1$ is $T_{01}$.
\end{itemize}
\end{enumerate}
\end{Pro}
\begin{proof}
According to the def\/inition of the operation {\bf j}, a~point $q_l$ is regular if and only if $j_l=1$, and according to
the proof of Lemma~\ref{sandwich}, $q_l$ is regular if and only if ${\mathcal A}_{l-1}(q_l)\not=0$.
Therefore, since the depth of $\Sigma_{j_1\dots j_k}$ is at most $2$, then $j_l\geq 2$ if and only if ${\mathcal A}_{l-1}=0$.
This ends the proof of~(1).

Now based on Proposition~\ref{delta} and the convention in the {\bf RVT} code, note that if the letter of rank~$h$ is of type
$T_{rs}$, there exists $\alpha<\beta$ such that the direction generated by $x_h-x_{h-1}$ belongs only the critical hyperplanes
$\delta_\alpha^{h-\alpha-1}$ and $\delta_\beta^{h-\beta-1}$.
Therefore from Proposition~\ref{delta} the point $q_h$ satisf\/ies the following constraints:
\begin{gather}\label{Tab}
{\mathcal A}_{l,h-\alpha-2}=0
\quad
\text{for}
\quad
h-\alpha\leq l<h
\qquad
\text{and}
\qquad
{\mathcal A}_{l,h-\beta-2}=0
\quad
\text{for}
\quad
h-\beta\leq l<h.
\end{gather}

Assume that there exists an index~$h$ such that $j_h \geq 3$.
The rule of least upward jumps assumes that there exists a~rank $l<h$ such that $j_l=3$.
Therefore we may assume that $j_h=3$.
In this case, the point $q_h$ is critical and based on Lemma~\ref{sandwich} we get ${\mathcal A}_{h-1}(q_h)=0$.
The def\/inition of the operation ${\bf j}$ implies that the projection of $X_h(q_h)$ on ${\mathcal C}^{h-1}(m)$ is
contained in an hyperplane of~$\delta_0^{h-1}$ at~$q_{h-1}$.
It follows that $q_h$ belongs to some critical manifold ${\mathcal C}_{q_l}^h(m-1)$ so that the point $q_{l+1}$
under $q_{h}$ is vertical.
In particular we have $l\leq h-2$.
Based on our our convention on the {\bf RVT} code, the letter of rank~$h$ is of type $T_{0i}$ with $i\geq 1$.
This ends the proof of~(2).

Assume now that the letter of rank~$h$ is $T_{rs}$.
As we have already seen, there exists $\alpha<\beta$ such that the direction generated by $x_h-x_{h-1}$ belongs to
$\delta_\alpha^{h-\alpha-1}$ and $\delta_\beta^{h-\beta-1}$.
Equations~\eqref{Tab} imply that the direction generated by $x_{h-\alpha}-x_{h-\alpha-1}$ belongs to
$\delta_0^{h-\alpha-1}$ and $\delta_{\beta-\alpha}^{h-\beta-\alpha-1}$.
Thus the letter of rank $h-\alpha$ is of type $T_{0i}$ with $i\geq 1$.
Let $\nu$ be the f\/irst rank such that the corresponding letter in the {\bf RVT} code is of type $T_{0i}$.
From our conventions on the {\bf RVT} code this letter must be $T_{01}$.
Assume that the direction generated by $x_\nu-x_{\nu-1}$ belongs to $\delta_0^{\nu-1}$ and
$\delta_{\gamma}^{\nu-\gamma-1}$.

From equations~\eqref{Tab} for $\nu-\gamma < j\leq \nu$ and $q_{\nu-\gamma}$, the point $q_j$ must be tangency and
$q_{\nu-\gamma}$ must be vertical.
This implies that $j_{\nu-\gamma}\geq 2$.
But from the choice of $\nu$ it follows that each letter of rank smaller than $\nu-1$ is of type $R$, $V$, $T$.
It follows that $j_{\nu-\gamma}=2$.
Therefore $q_j$ is tangency for $\nu-\gamma < j \leq\nu$ and the two points $q_{\nu-\gamma}$ and $q_\nu$ are vertical.
According to the choice of $\nu$, the letter of rank $\nu-\gamma$ is~$V$ and if $\gamma >1$ all letters of rank
$\nu-\gamma<l<\nu$ are equal to~$T$.
Again from the choice of $\nu$, $j_l=1$ for $h-\nu<l<\nu$ if $\nu>1$.
We get now,  at point $q_{\nu-1}$ the two relations
\begin{gather*}
\langle x_\nu-x_{\nu-1},x_{\nu-1}-x_{\nu-2}\rangle=0
\qquad
\text{and}
\qquad
\langle x_\nu-x_{\nu-1},x_{\nu-\gamma}-x_{\nu-\gamma-1}\rangle=0.
\end{gather*}
For a~given f\/ixed point $q_{h-1}$ we get two independent linear relations in $x_\nu-x_{\nu-1}$.
Now from the choice of $\nu$, each $j_l$ belongs to $\{1,2\}$ for $1\leq l<\nu$.
Therefore from the def\/inition of the operation $\bf j$ and the rule of least upward jumps it follows that $j_\nu=3$.
This ends the proof of~(3).

Finally assume that the letter of rank~$h$ is $T_{01}$ and~$h$ is the f\/irst index~$h$ such that $j_h=3$.
Note that if ${\mathcal A}_h(q_{h+1})=0$, again from Theorem~\ref{tengency}(1) this implies that there exists
$\alpha>0$ such that the point $q_{h-\alpha}$ is vertical and $q_{j}$ is tangency for $h-\alpha<j\leq h$.
In particular $q_h$~belongs to ${\mathcal C}_{q_{h-\alpha-1}}^{h}(m-1)$.

On the one hand assume that $j_{h+1}=2$.
This condition is equivalent to ${\mathcal A}_h(q_{h+1})=0$ and the projection of $X_{h+1}(q_{h+1})$ onto ${\mathcal
C}^{h}(m)$ must not be tangent to ${\mathcal C}_{q_{h-\alpha-1}}^{h}(m-1)$ at $q_h$ otherwise $j_{h+1}=3$.
This implies that $q_{h+1}$ is vertical but the letter of rank $h+1$ can not be $T_{01}$ otherwise the projection of
$X_{h+1}(q_{h+1})$ onto ${\mathcal C}^{h}(m)$ must be tangent to ${\mathcal C}_{q_{h-\alpha-1}}^{h}(m-1)$ at $q_h$.
Assume that the letter of rank $h+1$ is of type $T_{rs}$, of course we have $T_{rs}\not=T_{01}$.
Now from the proof of~(2), the letter of rank $h +1$ and the previous must be of type $T_{0i}$ with $i\geq 1$.
Taking into account our convention on the {\bf RVT} code we get $T_{rs}=T_{02}$

Conversely if the letter of rank $h+1$ is~$V$ then we get $j_{h+1}=2$.
If this letter is $T_{02}$, this means that ${\mathcal A}_{h, h-2}(q_{h+1})\not=0$ but ${\mathcal A}_h(q_{h+1})=0$.
It follows that the projection of $X_{h+1}(q_{h+1})$ onto ${\mathcal C}^{h}(m)$ is not tangent to ${\mathcal
C}_{q_{h-\alpha-1}}^{h}(m-1)$ and so we have $j_{h+1}=2$.

On the other hand assume that $j_{h+1}=3$.
This implies that ${\mathcal A}_h(q_{h+1})=0$ and the letter of rank $h+1$ must be of type $T_{rs}$.
Suppose that $T_{rs}\not=T_{01}$, by the same argument as previously, we obtain $T_{rs}=T_{02}$ which implies
that $j_{h+1}=2$ and gives rise to a~contradiction.
Therefore we must have $T_{rs}=T_{01}$.
Conversely if the letter of rank $h+1$ is $T_{01}$, as we have already seen, we have $j_{h+1}=3$.
\end{proof}

\begin{proof}[Proof of Theorem~\ref{stratEKR}]
(1) and (2) are consequences of  Propositions~\ref{cactEKR3}(1),~\ref{caractV} and  Theorem~\ref{verticalsets}(1).
Now from the convention on the {\bf RVT} code and  Proposition~\ref{cactEKR3}(2) the {\bf RVT} code of a~point
in an EKR set of depth $1$ contains only letters in $\{R,V,T\}$.
Therefore property~(i) in~(3) is a~consequence of Propositions~\ref{caractV} and~\ref{cactEKR3}(1).

We now focus on property~(ii).
If $\omega$ is of type $R^{h_0}VT^{l_1}R^{h_1}\dots V T^{l_\nu}R^{h_{\nu}}$ then the letters~$V$ are exactly at rank
$i_1,\dots, i_\nu$ with the relations $l_{\lambda}+h_{\lambda}=i_{\lambda+1}-i_{\lambda}-1$ for $\lambda=1,\dots,\nu-1$
and $l_\nu+h_\nu= k-i_\nu-1$, so the set ${\mathcal C}_{R^{h_0}VT^{l_1}R^{h_1}\dots V T^{l_\nu}R^{h_{\nu}}}$ must be
contained in the set $\Sigma_{j_1\dots j_k}$ from~(i).
Consider any word $\omega$ in a~{\bf RVT} code such that each letter of rank $i_1,\dots, i_\nu$
is~$V$, and take any $q\in{\mathcal C}_\omega\subset {\mathcal C}^k(m)$.
If $i_2=i_1+1$ the {\bf RVT} code of $q_{i_2-1}$ is of type $R^{h_0}VT^{l_1}R^{h_1}$ with $l_1=h_1=0$.
Assume now that $i_2-i_1>1$.
If for all $i_1<j< i_2$ each point $q_j$ is regular and then the {\bf RVT} code of $q_{i_2-1}$ is of type
$R^{h_0}VT^{l_1}R^{h_1}$ with $l_1=0$.
Now suppose that there exists some $q_j$ under~$q$ with $i_1<j< i_2$ which is critical.
The point $ q_j$ must be tangency by Theorem~\ref{tengency}(1), and, moreover, for $i_1<j'\leq j$, the point
$q_{j'}$ must also be tangency.
We set $l_1=\max\{j-i_1\,|\,   q_j\; \textrm{tangency}\}$.
Then, for $i_1+l_1<j<i_2$, the point $q_j$ must be regular otherwise from the previous argument $q_j$ must be tangency
which contradicts the def\/intion of $l_1$.
It follows that the {\bf RVT} code of $q_{i_2-1}$ is of type $R^{h_0}VT^{l_1}R^{h_1}$.
By induction on $1\leq i\leq \nu$, the same arguments shows that $\omega$ must be of type $R^{h_0}VT^{l_1}R^{h_1}\dots V
T^{l_\nu}R^{h_{\nu}}$.

Finally, from the proof of Theorem~\ref{verticalsets}(2), it follows that equations of ${\mathcal
C}_{R^{h_0}VT^{l_1}R^{h_1}\dots V T^{l_\nu}R^{h_{\nu}}}$ is the union of $\nu$ systems of type~\eqref{eqRVT}, and so we
get a~set of $\nu+l_1+\dots+l_\nu$ independent equations.
It~follows that ${\mathcal C}_{R^{h_0}VT^{l_1}R^{h_1}\dots V T^{l_\nu}R^{h_{\nu}}}$ is an analytic submanifold of
${\mathcal C}^k(m)$ of codimension $\nu+l_1+\dots+l_\nu$.
On the other hand, the equations of $\Sigma_{j_1\dots j_k}$ are ${\mathcal A}_{i_{\lambda}-1}=0$ for $\lambda=1,\dots,\nu$.
These equations are exactly the f\/irst equations of the $\nu$ systems of type~\eqref{eqRVT}
which def\/ine ${\mathcal C}_{R^{h_0}VT^{l_1}R^{h_1}\dots V T^{l_\nu}R^{h_{\nu}}}$.
This ends the proof of property~(ii).

First of all from Proposition~\ref{cactEKR3}(3), if $d=1$ then the depth of any word in {\bf RVT}
code is~$1$.
Therefore property~(iii) is a~direct consequence of properties~(i) and~(ii).
\end{proof}

\subsection{EKR sets of depth 1, {\bf RVT} codes and conf\/igurations\\ of an articulated arms}\label{depth1}

We will now give a~complete interpretation of the previous result in terms of conf\/igurations of an articulated arm as
stated in Theorem~\ref{3}(2):
\begin{The}\label{RVTconf1}
Let $\Sigma_{j_1\dots j_k}$ be an EKR set of depth $1$ in ${\mathcal C}^k(m)$ and $\{i_1,\dots, i_\nu\}$ the set $\{i
\in \{2,\dots,k\}\,|\,   j_i=2\}$.

A configuration $q\in \Sigma_{j_1\dots j_k}$ belongs to the RVT set ${\mathcal C}_{R^{h_0}VT^{l_1}R^{h_1}\dots V
T^{l_\nu}R^{h_{\nu}}}\subset \Sigma_{j_1\dots j_k}$ if and only, at~$q$, the only orthogonality constraint is that each segment
$[M_{{i_\lambda}+l-1},M_{{i_\lambda}+l}]$ is orthogonal to the direction on $\mathbb{R}^{m+1}$ generated by
$\overrightarrow{M_{i_\lambda-2}M_{i_\lambda-1}}$ for all $l=0,\dots,l_{\lambda}$ and $\lambda=1,\dots,\nu$.
\end{The}

\begin{Rem}
The property (ii) of Theorem~\ref{tengency} is a~particular case of Theorem~\ref{RVTconf}.
Note that we can f\/ind a~similar result in~\cite{S} with more restricted context.
\end{Rem}

For the proof of this result, we need the notion of ``induced articulated arm''.

Given two integers~$r$ and~$s$ such that $0\leq r<s\leq k$, we can look for the motion of an {\it induced articulated
arm}, which consists of segments of the original articulated arm joining $M_{r}$ to $M_{s}$.
We can then study the {\it induced articulated arm} $(M_{r},\dots, M_s)$.
We def\/ine $\kappa=s-r$, and we denote by ${\mathcal C}^{rs}(m)$ the image of ${\mathcal C}^k(m)$ by the canonical
projection $\varrho^{rs}$ from $\mathbb{R}^{m+1}_0\times\dots\times\mathbb{R}^{m+1}_k$ onto $\mathbb{R}^{m+1}_r\times\dots\times\mathbb{R}^{m+1}_s$.

In fact, we have: ${\mathcal C}^{rs}(m)=\{q_{rs}=(x_r,x_{r+1},\dots x_s)\,|\,   q=(x_0,\dots, x_k)\in {\mathcal C}^k(m)\}$.

Taking into account Section~\ref{arm}, let ${\mathcal E}_{rs}$ be the distribution on $(\mathbb{R}^{m+1})^{\kappa+1}$
spanned by
\begin{gather*}
{\mathcal Z}_{r},
\quad
\dots,
\quad
{\mathcal Z}_{s-1},
\quad
 \frac{\partial}{\partial x_{s}^1},
\quad
\dots,
\quad
 \frac{\partial}{\partial x_{s}^{m+1}}
\end{gather*}
and let ${\mathcal D}_{rs}$ be the distribution induced by ${\mathcal E}_{rs}$ on ${\mathcal C}^{rs}(m)$.

In terms of Notations~\ref{nota2}, the mechanical system describing the evolution of an induced arm $(M_{r},\dots, M_s)$
is a~controlled system on $\mathbb{R}^{m+1}\times(\mathbb{S}^m)^{\kappa}\equiv{\mathcal C}^\kappa(m)$
naturally associated with the distribution ${\mathcal D}_{rs}$.

Consider a~word $R^{h_0}VT^{l_1}R^{h_1}\dots V T^{l_\nu}R^{h_{\nu}}$ of~$k$ letters in a~{\bf RVT} code, and associate
with this word the sequences $r_0,\dots,r_{\nu}$ and $s_0,\dots, s_\nu$ def\/ined by
\begin{itemize}\itemsep=0pt
\item ${s}_0=h_0$ and $r_0=0$,
\item ${s}_{i}={s}_{i-1}+h_i+l_i+1=h_0+h_1+l_1+1+\dots+h_i+l_i+1$ and $r_i=s_{i-1}-1 $ for $i=1,\dots,\nu$.
\end{itemize}

We get then the following characterization:
\begin{Lem}\label{RVTconf}
The configuration $q\in{\mathcal C}^k(m)$ of an articulated arm $(M_0,\dots, M_k)$ belongs to
the set ${\mathcal C}_{R^{h_0}VT^{l_1}R^{h_1}\dots V T^{l_\nu}R^{h_{\nu}}}$ if and only if the induced articulated arm associated
with the pair of integers $({r}_{i},{s}_{i})$ is such that $\varrho^{{r}_{i}{s}_{i}}(q)$ belongs
to the set ${\mathcal C}_{R^{h_0}}\subset {\mathcal C}^{{r}_0{s}_0}(m)={\mathcal C}^{s_0}(m)$ for $i=0$
and belongs  to the set ${\mathcal C}_{RVT^{h_i}R^{h_i}}\subset {\mathcal C}^{{r}_{i}{s}_{i}}(m)$ for all $i=1,\dots, \nu$.
\end{Lem}

\begin{proof}
For any $q\in{\mathcal C}^k(m)$ we denote as usual by $q_l$ any point of ${\mathcal C}^l(m)$ under~$q$, and we f\/ix
a~conf\/iguration $q\in {\mathcal C}^k(m)$ of the articulated arm $(M_0,\dots, M_k)$.
First of all, for $i=0$, the induced articulated arm associated with $(r_0,s_0)$ has the induced conf\/iguration~$q_{s_0}$.
The {\bf RVT} code of~$q_{s_0}$ consists of~$h_0$ f\/irst letters of the {\bf RVT} code of~$q$.
Therefore, these f\/irst $h_0$ letters are~$R^{h_0}$ if and only if the {\bf RVT} code of~$q_{s_0}$ is~$R^{h_0}$.

Assume that the f\/irst~$s_i$ letters of the {\bf RVT} code of~$q$ are $R^{h_0}VT^{l_1}R^{h_1}\dots$ $V
T^{l_i}R^{h_{i}}$ if and only if the {\bf RVT} code of the conf\/iguration $q_{r_is_i}=\varrho^{{r}_{i}{s}_{i}}(q)$ of the
associated induced articulated arm is ${R^{h_0}}$ for $i=0$ and ${RVT^{l_i}R^{h_i}}$ for all $1 \leq i\leq \mu-1<\nu$.
Consider the conf\/iguration $q_{r_\mu s_\mu}=\varrho^{{r}_{\mu}{s}_{\mu}}(q)\in{\mathcal C}^{r_\mu s_\mu}(m)$ of the
associated induced articulated arm.
Denote by~$q_{{r_\mu}l}$ the conf\/iguration under $q_{r_\mu s_\mu}$ for $r_\mu\leq l\leq s_\mu$, and set $\kappa_\mu=s_\mu-r_\mu$.
By convention, the {\bf RVT} code of $q_{{r_\mu}{r_\mu+1}}$ is~$R$.
Now, according to Proposition~\ref{caractV}, $q_{{r_\mu}{r_\mu+2}}$
is vertical in ${\mathcal C}^{r_\mu s_\mu}(m)\equiv{\mathcal C}^{\kappa_\mu}(m)$ if and only if
\begin{gather*}
 \sum\limits_{j=1}^{m+1}\big(x_{r_{\mu}+2}^j-x_{r_{\mu}+1}^j\big)\big(x_{r_{\mu}+1}^j-x_{r_{\mu}}^j\big)=0.
\end{gather*}
This is equivalent to ${\mathcal A}_{r_\mu+1}(q)=0$.
It follows that $q_{{r_\mu}{r_{\mu+2}}}$ is vertical if and only if~$q_{r_{\mu+2}}$ is vertical.
Finally, the letter of rank~$2$ in the {\bf RVT} code of $q_{r_\mu s_\mu}$ is~$V$ if and only if the letter of rank
$s_\mu+1=h_0+h_1+l_1+1+\dots+h_{\mu}+l_{\mu}+2$ is also~$V$.
Consider now an integer $r_{\mu}+2+l$ with $0\leq l\leq l_{\mu}$.
Either $q_{r_\mu r_{\mu}+l+2}$ is critical or it is regular.
If $q_{r_\mu r_{\mu}+l+2}$ were vertical, then, from the previous argument, $q_{r_{\mu}+l+2}$ would also be vertical,
which contradicts the def\/inition of the set $\{i_1,\dots, i_\nu\}$.
Assume that $q_{r_\mu r_{\mu}+l+2}$ is tangency.
Then $q_{r_\mu r_{\mu}+l'+2}$ is also tangency for all $0\leq l'\leq l$.
We have this property if and only if the following relations hold (see~\ref{eqRVT}):
\begin{gather*}
\langle x_{r_\mu+l'+3}-x_{r_\mu+l'+2}, x_{r_\mu+1}-x_{r_\mu}\rangle=0
\quad
\text{for all}
\quad
0\leq l'\leq l.
\end{gather*}
According to our assumption and the equations of ${\mathcal C}_{R^{h_0}VT^{l_1}R^{h_1}\dots
VT^{l_{\mu-1}}R^{h_{\mu-1}}}$ (see proof of property~(ii) in Theorem~\ref{stratEKR}), $q_{r_\mu r_{\mu+l+2}}$ is
tangency if and only if $q_{r_{\mu}+l+2}$ is also tangency.

Now, $q_{r_\mu r_{\mu+l+2}}$ is regular if and only if $\langle x_{r_\mu+l+3}-x_{r_\mu+l+2}, x_{r_\mu+1}-x_{\mu}\rangle\not=0$.
On the one hand, if $q_{r_{\mu}+l+2}$ is regular we must have $\langle x_{r_\mu+l+3}-x_{r_\mu+l+2},
x_{r_\mu+1}-x_{r_\mu}\rangle\not=0$ and so $q_{r_\mu r_{\mu+l+2}}$ is regular.
On the other hand, according to the choice of the set $\{i_1,\dots, i_\nu\}$, if $q_{r_{\mu}+l+2}$ is critical then it must be
tangency.
From the def\/inition of the {\bf RVT} code, since $x_{r_\mu+l+3}-x_{r_\mu+l+2}$ belongs to one and only one critical
hyperplane then $\langle x_{r_\mu+l+3}-x_{r_\mu+l+2}, x_{r_\mu+1}-x_{r_\mu}\rangle=0$ and therefore $q_{r_\mu
r_{\mu}+l+2}$ can not be regular.
$q_{r_\mu r_{\mu+l+2}}$ is regular if and only if $q_{r_{\mu+l+2}}$ is regular.

It follows that our assumption is then true for the integer $\mu$.
\end{proof}

\begin{proof}[Proof of Theorem~\ref{RVTconf1}]\sloppy
Based on Theorem~\ref{tengency}(3), for each induced articulated arm associated
with a~pair $(r_i,s_i)$, there exists a~family of directions $\{K_i(q)\}$ in $\mathbb{R}^{m+1}$ for
$q\in {\mathcal C}^{r_i s_i}(m)$ generated by $x_{r_i+1}-x_{r_i}$ such that the conf\/iguration $q_{r_i,s_i}$ belongs to
the set ${\mathcal C}_{R^VT^{l_i}R^{h_i}}\subset {\mathcal C}^{r_i s_i}(m)$ if and only if this conf\/iguration fulf\/ills
the following property: each segment $[M_{{r_i}+1+l},M_{{r_i}+2+l}]$ is orthogonal at $M_{{r_i}+1+l}$ to $K_i(q_{r_i
s_i})$, for $l=0,\dots, l_i$ and there is no other orthogonality constraint.

The theorem is then a~consequence of Lemma~\ref{RVTconf}.
\end{proof}

\subsection[EKR sets of depth 2, {\bf RVT} codes and conf\/igurations of an articulated arm for $1\leq k\leq 4$]{EKR sets of depth 2, {\bf RVT} codes and conf\/igurations\\ of an articulated arm for $\boldsymbol{1\leq k\leq 4}$}
\label{2k4}

The combination of all possible {\bf RVT} codes of depth $2$ has an exponential growth relatively to the
length~$k$ of a~special multi-f\/lag.
Therefore, in this subsection we only describe the relations between EKR sets of $2$-depth, {\bf RVT} codes and
conf\/igurations of articulated arms for $k=4$.
In fact, this situation corresponds to the results of~\cite {CH,CMH,MP}.

First of all, for $k=3$, we have only $\Sigma_{123}$, which is an EKR set of depth $2$, and for $k=4$ we have fourteen
EKR sets (of depth $2$) whose numerical codes are (see for instance~\cite{MP}) $1111$, $1112$, $1121$, $1122$,
$1123$, $1211$, $1212$, $1213$, $1221$, $1222$, $1223$, $1231$, $1232$, $1233$.

Therefore, for $1\leq k\leq 4$, the other EKR sets for $1\leq k\leq 4$ are of depth $1$.

Recall that at the end of Section~\ref{code} we have seen that for $k=3$ we have only one {\bf RVT} set of depth $2$
(i.e.\  $RT_0T_{01}$) but we have ten {\bf RVT} sets for $k=4$.
All other {\bf RVT} sets are of depth at most $1$.

Notice that the decomposition of EKR sets of depth $1$ into {\bf RVT} sets are of depth $1$ can be found in
Theorem~\ref{stratEKR}, and the corresponding interpretation in terms of conf\/igurations of an articulated arm can also
be found in Theorem~\ref{RVTconf}.
This is why we have given such results only for EKR sets of depth $2$ previously enumerated.

For this purpose, we need the following characterizations of some EKR sets of depth $2$ which is an easy consequence of
Proposition~\ref{cactEKR3}:
\begin{Pro}
\label{cactEKR3bis}
Let $\Sigma_{{j_1}\dots{j_k}}$ be an EKR set of depth $2$ in ${\mathcal C}^{k}(m)$ with $k\geq 3$.
Consider an integer $2\leq h<k$.
Assume that $j_l\in\{1,2\}$ for all $1\leq l\leq h-1$ and denote by $\{i_1,\dots, i_\nu\} $ the set of indexes
$i\in\{1,\dots,h-1\} $ such that $ j_i=2$.
\begin{enumerate}\itemsep=0pt
\item[$1.$] $j_h=3$ if and only if the letter of rank~$h$ of the {\bf RVT} code of any $q\in\Sigma_{{j_1}\dots{j_k}}$ is
$T_{01}$ where $T_1$ refers to the singular hyperplane $\delta_{i_\nu}^{h-i_\nu-1}$.
\item[$2.$] if $j_h=3$ then we have one and only one of the following situations:
\begin{itemize}\itemsep=0pt
\item $j_{h+1}=1$ if and only if the letter of rank $h+1$ belongs to the set $\{R,T_1, T_{2},T_{12}\}$,
\item $j_{h+1}=2$ if and only if the letter of rank $h+1$ belongs to the set $\{V,T_{02}\}$,
\item $j_{h+1}=3$ if and only if the letter of rank $h+1$ is $T_{01}$.
\end{itemize}
\end{enumerate}
\end{Pro}

For $k=4$, we can easily get the decomposition of an EKR set of depth at most $2$ into {\bf RVT} sets as given in the
following table by application of the previous proposition.
For $k\leq 3$, the results are particular cases of Theorem~\ref{stratEKR}.
\begin{center}
{\it Decomposition of EKR classes into {\bf RVT} classes}
\vspace{1mm}

\begin{tabular}{|c|c|rl}
\hline
EKR class & {\bf RVT} class
\\
\hline
$1111$ & $RRRR$
\\
\hline
$1112$ & $RRRV$
\\
\hline
$1121$ & $RRVR$, $RRVT$
\\
\hline
$1122$ & $RRVV$
\\
\hline
$1123$ & $RRT_0T_{01}$
\\
\hline
$1211$ & $RVRR$, $RVTR$, $RVTT$
\\
\hline
$1212 $ & $RVRV$, $RVTV$
\\
\hline
$1213$ & $RT_0T_1T_{01}$
\\
\hline
$1221$ & $RVVR$, $RVVT$
\\
\hline
$1222$ & $RVVV$
\\
\hline
$1223$ & $RVT_0T_{01}$
\\
\hline
$1231$ & $RT_0T_{01}R$, $RT_0T_{01}T_1$, $RVT_0T_{01}T_2$, $RVT_0T_{01}T_{12}$
\\
\hline
$1232$ & $RT_0T_{01}V$, $RT_0T_{01}T_{02}$,
\\
\hline
$1233$ & $RT_0T_{01}T_{01}$.
\\
\hline
\end{tabular}
\end{center}

For $1\leq k\leq 4$, we only give an interpretation of {\bf RVT} sets in terms of conf\/igurations of articulated arm when
the {\bf RVT} code contains a~letter of type $T_{ij}$.
The other cases are particular cases of Theorem~\ref{RVTconf}.
The proof of the following descriptions are obtained from the decomposition of each EKR set in {\bf RVT} sets given in
the previous table and by an easy interpretation in terms of orthogonality of the equations of type ${\mathcal
A}_{j,i}=0$ of each such sets (see the proof of Theorem~\ref{tengency}).
These proofs are left to the reader.

{\it Let $q=(x_0,\dots, x_k)$ be a~configuration of an articulated arm $(M_0,\dots,M_k)$ with $k\leq 4$.
We have the following characterizations:}
\begin{itemize}\itemsep=0pt
\item $q$ belongs to $\Sigma_{123}={\mathcal C}_{RT_0T_{01}}$ if and only if $[M_{i-2},M_{i-1}]$ and
$[M_{i-1},M_i]$ are orthogonal at $M_{i-1}$ for $i=2,3$, $[M_{2},M_3]$ is  orthogonal to the direction generated
by $\overrightarrow{M_0M_1}$ and no other orthogonality in the conf\/iguration~$q$.
\item $q$ belongs to $\Sigma_{1123}={\mathcal C}_{RRT_0T_{01}}$ if and only if $[M_{i-2},M_{i-1}]$ and
$[M_{i-1},M_i]$ are orthogonal at $M_{i-1}$ for $i=3,4$, $[M_{3},M_4]$ is orthogonal to the direction generated
by $\overrightarrow{M_1M_2}$ and no other orthogonality in the conf\/iguration~$q$.
\item $q$ belongs to ${\mathcal C}_{RT_0T_1T_{01}}=\Sigma_{1213}$ if and only if $[M_{i-2},M_{i-1}]$ and
$[M_{i-1},M_i]$ are orthogonal at $M_{i-1}$ for $i=2,4$, $[M_{3},M_4]$ is orthogonal to the direction generated
by $\overrightarrow{M_0M_1}$ and no other orthogonality in the conf\/iguration~$q$.
\item $q$ belongs to $\Sigma_{1223}={\mathcal C}_{RVT_0T_{01}}$ if and only if $[M_{i-2},M_{i-1}]$ and
$[M_{i-1},M_i]$ are orthogonal at~$M_{i-1}$ for $i=2,3,4$, $[M_{3},M_4]$ is orthogonal to the direction
generated by $\overrightarrow{M_0M_1}$ and no other orthogonality in the conf\/iguration~$q$.
\item in $\Sigma_{1231}$:
\begin{itemize}\itemsep=0pt
\item[(i)] $q$ belongs to ${\mathcal C}_{RT_0T_{01}R} $ if and only if $[M_{i-2},M_{i-1}]$ and $[M_{i-1},M_i]$ are orthogonal
at~$M_{i-1}$ for $i=2,3$, $[M_{2},M_3]$ is orthogonal to the direction generated by $\overrightarrow{M_0M_1}$
and no other orthogonality in the conf\/iguration~$q$.
\item[(ii)] $q$ belongs to ${\mathcal C}_{RT_0T_{01}T_1}$ (resp.\
${\mathcal C}_{RT_0T_{01}T_2}$) if and only if the previous constraints hold, $[M_{3},M_4]$ is
orthogonal to the direction generated by $\overrightarrow{M_1M_2}$ (resp.
$\overrightarrow{M_0M_1}$) and no other orthogonality in the conf\/iguration~$q$.
\item[(iii)] $q$ belongs to ${\mathcal C}_{RT_0T_{01}T_{12}}$ if and only if we have the previous constraints of (ii) hold,
 $[M_{3},M_4]$ is orthogonal at the directions generated by $\overrightarrow{M_0M_1}$ and by~$\overrightarrow{M_1M_2}$ and no other orthogonality in the conf\/iguration~$q$.
\end{itemize}
\item in $\Sigma_{1232}$:
\begin{itemize}\itemsep=0pt
\item[(i)] $q$ belongs to ${\mathcal C}_{RT_0T_{01}V}$ if and only if $[M_{i-2},M_{i-1}]$ and $[M_{i-1},M_i]$ are orthogonal
at~$M_{i-1}$ for $i=2,3,4$, $[M_{2},M_3]$ is orthogonal to the directions generated $\overrightarrow{M_0M_1}$
and no other orthogonality in the conf\/iguration~$q$.
\item[(ii)] $q$ belongs to ${\mathcal C}_{RT_0T_{01}T_{02}}$ if and only if $[M_{i-2},M_{i-1}]$ and $[M_{i-1},M_i]$ are
orthogonal at $M_{i-1}$ for $i=2,3,4$, $[M_2,M_3]$ and $[M_{3},M_4]$ are orthogonal at the direction generated
by $\overrightarrow{M_0M_1}$ and $\overrightarrow{M_1M_2}$ respectively and no other orthogonality in the
conf\/iguration~$q$.
\end{itemize}
\item The point $q$ belongs to ${\mathcal C}_{RT_0T_{01}T_{01}}=\Sigma_{1233}$ if and only if $[M_{i-2},M_{i-1}]$
and $[M_{i-1},M_i]$ are orthogonal at $M_{i-1}$ for $i=2,3,4$, $[M_2,M_3]$ and $[M_{3},M_4]$ are orthogonal at
the direction generated by $\overrightarrow{M_0M_1}$ and no other orthogonality in the conf\/iguration~$q$.
\end{itemize}
\begin{proof}[Proof of Proposition~\ref{cactEKR3bis}]
(1) is established in the proof of Proposition~\ref{cactEKR3}(3).
In~(2) the last two situations correspond to Proposition~\ref{cactEKR3}(4).
The f\/irst situation is an elementary computation in terms of critical hyperplane and is left to the reader.
\end{proof}

\section*{Main notations}
\begin{itemize}\itemsep=0pt
\item $\mathbb{D}:D=D_{k}\subset D_{k-1}\subset\dots\subset D_{j}\subset\dots\subset D_{1} \subset D_{0}= TM $:
special multi-f\/lag of step $m\geq 2$ and
length $k\geq 1$ (Section~\ref{multiflag}).

\item $\begin{matrix}
D_{j}&\subset& D_{j-1}
\\
\hfill\cup\hfill&& \hfill\cup\hfill
\\
L(D_{j-1})& \subset& L(D_{j-2})
\end{matrix}$: sandwich of rank~$j$ (Section~\ref{multiflag}).

\item $ P^k(m)\rightarrow P^{k-1}(m)\rightarrow\dots\rightarrow P^{1}(m)\rightarrow
P^0(m):=\mathbb{R}^{m+1}$:
tower of projective bundles (Section~\ref{projtower}).

\item $\Delta_j$: typical distribution on $P^{j}(m)$ which is the Cartan prolongation of $\Delta_{j-1}$
(Section~\ref{projtower}).

\item $S(D,M,g)$: sphere bundle in the distribution~$D$ associated with Riemannain metric~$g$
(Section~\ref{cartsph}).

\item $ \hat{P}^k(m)\rightarrow \hat{P}^{k-1}(m)\rightarrow\dots\rightarrow \hat{P}^{1}(m)\rightarrow
\hat{P}^0(m):=\mathbb{R}^{m+1}$:
tower of sphere bundles (Section~\ref{cartsph}).

\item $\hat{\Delta}_j$: typical distribution on $ \hat{P}^j(m)$ which is the spherical prolongation of
$\hat{\Delta}_{j-1}$ (Section~\ref{cartsph}).

\item $(M_0,\dots,M_k)$ articulated arm or system of rigid bars in $\mathbb{R}^{m+1}$ of length~$k$
(Section~\ref{arm}).

\item ${\mathcal C}^k(m)$: conf\/iguration space of an articulated arm in $\mathbb{R}^{m+1}$ of length~$k$
(Section~\ref{arm}).

\item $q=(x_1,\dots,x_k)\in {\mathcal C}^k(m)$: conf\/iguration of an articulated arm where
$x_i=\big(x_i^1,\dots, x_i^r,\dots$, $x_i^{m+1}\big)$ for $i=0,\dots,k$.

\item ${\mathcal D}_k$: typical distribution on ${\mathcal C}^k(m)$ associated with an articulated arm of
length~$k$ (Section~\ref{arm}).

\item ${\mathcal Z}_i= \sum\limits_{r=1}^{m+1}\big(x_{i+1}^r-x_i^r\big)\frac{\partial}{\partial x_i^r}$
for $i=0,\dots,k-1$ (Section~\ref{arm}).

\item ${\mathcal A}_i= \sum\limits_{r=1}^{m+1}\big(x_{i+1}^r-x_i^r\big)\big(x_i^r-x_{i-1}^r\big)$ for
$i=0,\dots,k-1$ (Section~\ref{arm}).

\item ${\mathcal A}_{j,i}= \sum\limits_{r=1}^{m+1}\big(x_{j+1}^r-x_j^r\big)\big(x_{i+1}^r-x_{i}^r\big)$ for
$i=0,\dots,k-1$ and $i< j <k$ (in proof of Theorem~\ref{tengency}).

\item $Y_k=\Big( \sum\limits_{i=0}^{k-2} \prod\limits_{j=i+1}^{k-1} {\mathcal A}_j
{\mathcal Z}_i\Big)+{\mathcal Z}_{k-1}={\mathcal A}_{k-1}Y_{k-1}+{\mathcal Z}_{k-1}$ (Section~\ref{arm}).

\item $X_k=Y_k+ \sum\limits_{r=1}^{m+1}(x_k^r-x_{k-1}^r)\frac{\partial}{\partial x_k^r}$
(Section~\ref{arm}).

\item ${\mathcal D}_k$ is spanned by (Section~\ref{arm}):

\begin{itemize}\itemsep=0pt
\item[$\circ$] $\left\{\big(x_{k}^r-x_{k-1}^r\big)Y_k+ \frac{\partial}{\partial x_k^r},
\,
r=1,\dots m+1\right\}$,
\item[$\circ$] $\left\{\big(x_{k}^r-x_{k-1}^r\big)X_k+\Pi_k\left( \frac{\partial}{\partial x_k^r}\right)
,\,
r=1,\dots m+1\right\}$,
\item[$\circ$] $\left\{X_k, \Pi_k\left( \frac{\partial}{\partial x_k^r}\right):\; r=1,\dots,m+1\right\}$,
where $\Pi_k:T(\mathbb{R}^{m+1})^{k+1}\rightarrow T{\mathcal C}^k(m)$ is
orthogonal projection.
\end{itemize}

\item $\Psi^j: S({\mathcal D}_j,{\mathcal C}^j(m),\gamma_j)\rightarrow {\mathcal C}^{j+1}(m)$ such that
$\Psi^j_*(({\mathcal D}_j)^{[1]})={\mathcal D}^{j+1}$ (Section~\ref{armsph}).

\item $F^j:\hat{P}^j(m)\rightarrow {\mathcal C}^j(m)$ such that $F^j_*(\hat{\Delta}_j)={\mathcal D}_j$
(Section~\ref{armsph}).

\item ${\mathcal F}^k:{\mathcal C}^k(m)\rightarrow {\mathcal S}^k(m)\equiv \mathbb{R}^{m+1}
\times
(\mathbb{S}^m)^k$ with ${\mathcal F}^k(x_0,\dots,x_k)=(x_0,x_1-x_0,\dots,x_k-x_{k-1})$
(Section~\ref{armsph}).

\item $
\begin{cases}
z^1=\rho\phi^1(\theta)= \rho\sin {\theta^1} \dots\sin{\theta^{m-1}}\sin{\theta^m},
\\
z^2=\rho\phi^2(\theta)=\rho\sin
{\theta^1} \cdots\sin{\theta^{m-1}}\cos{\theta^m},
\\
z^3=\rho\phi^3(\theta)=\rho\sin {\theta^1}
\cdots\sin{\theta^{m-2}}\cos{\theta^{m-1}},
\\
\cdots\cdots\cdots\cdots\cdots\cdots\cdots\cdots\cdots\cdots\cdots\cdots\cdots
\\ z^{k}=\rho\phi^k(\theta)=\rho\sin{\theta^1} \cos{\theta^{2}},
\\
z^{k+1}=\rho\phi^{k+1}(\theta)=\rho\cos{\theta^1},
\end{cases}
$ hyperspherical coordinates
(Section~\ref{armsph}).

\item ${\mathcal P}^j(m)$: either $P^j(m)$ or $\hat{P}^j(m)$ (Section~\ref{code}).

\item $\mathfrak{D}_j$: the typical distribution on ${\mathcal P}^j(m)$ (i.e.\  either $\Delta_j$ or
$\hat{\Delta}_j$) (Section~\ref{code}).

\item $\mathfrak{d}_j^i$ with $j+i=k$: family of singular hyperplanes inside $\mathfrak{D}_k$
(Section~\ref{code}).

\item {\bf RVT} code  (Section~\ref{code}).

\item ${\mathcal C}_\omega$: set of conf\/igurations $q\in {\mathcal C}^k(m)$ whose {\bf RVT} code is the
word $\omega$ (Section\ref{code}).

\item operation {\bf j}  (Section~\ref{EKR}).

\item EKR class (Section~\ref{EKR}).

\item $\Sigma_{j_1\dots j_k}$: set of conf\/iguration $q\in {\mathcal C}^k(m)$ for which the germ of the
distribution ${\mathcal D}_k$ in~$q$ belongs to the EKR class coded by ${j_1\dots j_k}$.
\end{itemize}

\subsection*{Acknowledgments}

We would like to thank warmly the anonymous referees for the care and time they spent in the detailed reading of
dif\/ferent versions and for their questions, suggestions and comments that helped us to signif\/icantly improve the initial
version.

\pdfbookmark[1]{References}{ref}
\LastPageEnding


\begin{thebibliography}{99}
\footnotesize\itemsep=0pt

\bibitem{Ad1}
Adachi J., Global stability of distributions of higher corank of derived length
  one, \href{http://dx.doi.org/10.1155/S1073792803130735}{\textit{Int. Math. Res. Not.}} \textbf{2003} (2003), 2621--2638.

\bibitem{Ad2}
Adachi J., Global stability of special multi-f\/lags, \href{http://dx.doi.org/10.1007/s11856-010-0072-3}{\textit{Israel~J. Math.}}
  \textbf{179} (2010), 29--56.

\bibitem{CH}
Castro A.L., Howard W.C., A {M}onster {T}ower approach to {G}oursat
  multi-f\/lags, \href{http://dx.doi.org/10.1016/j.difgeo.2012.06.005}{\textit{Differential Geom. Appl.}} \textbf{30} (2012), 405--427.

\bibitem{CMH}
Castro A.L., Montgomery R., Spatial curve singularities and the
  {M}onster/{S}emple tower, \href{http://dx.doi.org/10.1007/s11856-012-0031-2}{\textit{Israel~J. Math.}} \textbf{192} (2012),
  381--427.

\bibitem{KuRu}
Kumpera A., Rubin J.L., Multi-f\/lag systems and ordinary dif\/ferential equations,
  \textit{Nagoya Math.~J.} \textbf{166} (2002), 1--27.

\bibitem{LR}
Li S.J., Respondek W., The geometry, controllability, and f\/latness property of
  the {$n$}-bar system, \href{http://dx.doi.org/10.1080/00207179.2011.569955}{\textit{Internat.~J. Control}} \textbf{84} (2011),
  834--850.

\bibitem{MZ1}
Montgomery R., Zhitomirskii M., Geometric approach to {G}oursat f\/lags,
  \href{http://dx.doi.org/10.1016/S0294-1449(01)00076-2}{\textit{Ann. Inst. H. Poincar\'e Anal. Non Lin\'eaire}} \textbf{18} (2001),
  459--493.

\bibitem{M1}
  Mormul P., Geometric singularity classes for special $k$-flags, $k \geq 2$, of arbitrary length,
in Singularity Theory Seminar, Editor S.~Janeczko, Warsaw University of Technology, Vol.~8, 2003, 87--100.

\bibitem{M2}
Mormul P., Multi-dimensional {C}artan prolongation and special {$k$}-f\/lags, in
  Geometric Singularity Theory, \href{http://dx.doi.org/10.4064/bc65-0-12}{\textit{Banach Center Publ.}}, Vol.~65, Polish
  Acad. Sci., Warsaw, 2004, 157--178.

\bibitem{MP}
Mormul P., Pelletier F., Special 2-f\/lags in lengths not exceeding four: a study
  in strong nilpotency of distributions, \href{http://arxiv.org/abs/1011.1763}{arXiv:1011.1763}.

\bibitem{PR1}
Pasillas-L{\'e}pine W., Respondek W., Contact systems and corank one involutive
  subdistributions, \href{http://dx.doi.org/10.1023/A:1013015602007}{\textit{Acta Appl. Math.}} \textbf{69} (2001), 105--128,
  \href{http://arxiv.org/abs/math.DG/0004124}{math.DG/0004124}.

\bibitem{PLR}
Pasillas-L{\'e}pine W., Respondek W., On the geometry of {G}oursat structures,
  \href{http://dx.doi.org/10.1051/cocv:2001106}{\textit{ESAIM Control Optim. Calc. Var.}} \textbf{6} (2001), 119--181,
  \href{http://arxiv.org/abs/math.DG/9911101}{math.DG/9911101}.

\bibitem{P}
Pelletier F., Espace de conf\/iguration d'un syst\`eme m\'ecanique et tours de
  f\/ibr\'es associ\'ees \`a un multi-drapeau sp\'ecial, \href{http://dx.doi.org/10.1016/j.crma.2011.12.007}{\textit{C.~R.~Math.
  Acad. Sci. Paris}} \textbf{350} (2012), 71--76.

\bibitem{SY}
Shibuya K., Yamaguchi K., Drapeau theorem for dif\/ferential systems,
  \href{http://dx.doi.org/10.1016/j.difgeo.2009.03.017}{\textit{Differential Geom. Appl.}} \textbf{27} (2009), 793--808.

\bibitem{S}
Slayman M., Bras articul\'e et distributions multi-drapeaux, Ph.D. Thesis,
  Universit\'e de Savoie, Laboratoire de Math\'ematiques (LAMA), 2008.

\bibitem{SP1}
Slayman M., Pelletier F., Articulated arm and special multi-f\/lags,
  \textit{J.~Math. Sci. Adv. Appl.} \textbf{8} (2011), 9--41,
  \href{http://arxiv.org/abs/1205.2990}{arXiv:1205.2990}.

\end{thebibliography}
\end{document}